\documentclass[12pt, a4paper]{article}
\usepackage{amssymb,latexsym, amsmath,amscd,amsfonts,amsthm,mathbbol, mathrsfs}
\usepackage[all]{xy}
\newtheorem{prop}{Proposition}[section]
\newtheorem{defi}[prop]{Definition}

\newtheorem{lem}[prop]{Lemma}
\newtheorem{thm}[prop]{Theorem}

\newtheorem{remar}[prop]{Remark}
\newtheorem{cor}[prop]{Corollary}

\DeclareMathAlphabet{\mathpzc}{OT1}{pzc}{m}{it}
\DeclareMathOperator{\Aut}{Aut}
\DeclareMathOperator{\End}{End}
\DeclareMathOperator{\Hom}{Hom}
\DeclareMathOperator{\Ind}{Ind}
\DeclareMathOperator{\cInd}{c-Ind}

\DeclareMathOperator{\PGL}{PGL}
\DeclareMathOperator{\GL}{GL}
\DeclareMathOperator{\SL}{SL}
\DeclareMathOperator{\Ker}{Ker}

\DeclareMathOperator{\Image}{Im}

\DeclareMathOperator{\Stab}{Stab}

\DeclareMathOperator{\Inj}{Inj}
\DeclareMathOperator{\inj}{inj}
\DeclareMathOperator{\Irr}{Irr}
\DeclareMathOperator{\supp}{Supp}
\DeclareMathOperator{\id}{id}
\DeclareMathOperator{\Mod}{Mod}

\DeclareMathOperator{\ev}{ev}
\DeclareMathOperator{\Fr}{Fr}
\DeclareMathOperator{\St}{St}
\DeclareMathOperator{\soc}{soc}
\DeclareMathOperator{\rad}{rad}
\DeclareMathOperator{\Rep}{Rep}
\newcommand{\cIndu}[3]{\cInd_{#1}^{#2}{#3}}

\newcommand{\Indu}[3]{\Ind_{#1}^{#2}{#3}}
\newcommand{\KK}{\mathfrak{K}}

\newcommand{\FF}{\mathcal F}

\newcommand{\Fq}{\mathbf{F}_{q}}
\newcommand{\Fbar}{\overline{\mathbf{F}}_{p}}

\newcommand{\oo}{\mathfrak o}
\newcommand{\oF}{\oo_{F}}
\newcommand{\pF}{\mathfrak{p}_{F}}

\newcommand{\pif}{\varpi_{F}}

\newcommand{\ZZ}{\mathbb{Z}}
\newcommand{\VV}{\mathcal{V}}
\newcommand{\LLL}{\mathcal {L}}
\newcommand{\DD}{\mathcal{D}}

\newcommand{\kzr}{\mathfrak K(\sigma_{0})}
\newcommand{\kon}{\mathfrak K(\sigma_{1})}
\newcommand{\CC}{\mathcal C}

\newcommand{\QK}[1]{I_{1}}
\newcommand{\coeff}{\mathcal{COEF}_G}
\newcommand{\diag}{\mathcal{DIAG}}
\newcommand{\WW}{\mathcal W}
\newcommand{\const}{\mathcal K}
\newcommand{\HH}{\mathcal H}
\newcommand{\vs}[2]{V_{#1, #2}}

\newcommand{\rr}{\mathbf{r}}
\newcommand{\lie}{\mathfrak{g}}
\newcommand{\UU}{\mathcal U}
\newcommand{\CCC}{\mathbb C}
\newcommand{\II}{\mathcal I}

\newcommand{\val}{val}
\newcommand{\BB}{\mathcal B_{\gamma}}
\newcommand{\QQ}{\mathbf{Q}}
\newcommand{\PP}{\mathcal P}
\newcommand{\kK}{\mathcal K}
\newcommand{\GG}{\mathcal G}
\newcommand{\red}{\mathrm red}
\title{Coefficient systems and supersingular representations of $\GL_2(F)$}
\author{Vytautas Paskunas}
\date{\today.}
\begin{document} 
\maketitle 
\tableofcontents
\pagebreak
\section{Introduction}
Recently Breuil in \cite{breuil} has determined the isomorphism classes of the irreducible smooth $\Fbar$-representations 
of $\GL_2(\QQ_p)$. This allowed him to define a ``correspondance semi-simple modulo $p$ pour $\GL_2(\QQ_p)$''.  Under this correspondence
 the isomorphism classes of irreducible smooth $2$-dimensional  $\Fbar$-representations of the Weil group of $\QQ_p$ are in bijection with 
 the isomorphism classes of ``supersingular'' irreducible smooth $\Fbar$-representations of $\GL_2(\QQ_p)$. The term ``supersingular'' was coined 
by Barthel and Livn\'e. Roughly speaking a supersingular representation is the $\Fbar$ analogue of a supercuspidal representation over $\CCC$, 
see Definition \ref{superdefi}. Let $F$ be a non-Archimedean local field, with a residue class field  $\Fq$ of the characteristic $p$. All 
the irreducible smooth $\Fbar$-representations of $G=\GL_2(F)$, which are \underline{not} supersingular,  have  been determined by Barthel and Livn\'e 
in \cite{bl1} and \cite{bl2}, and also by Vign\'eras in \cite{vig}, with no restrictions 
on $F$. However, if $F\neq\QQ_p$ then the method of Breuil fails and relatively little is known about the supersingular representations of $G$.   

This paper is an attempt to shed some light on this question. We fix a uniformiser $\pif$  of $F$ and  we construct $q(q-1)/2$ pairwise 
non-isomorphic, irreducible, supersingular, admissible (in the usual smooth sense) representations of $G$, which admit a central character, such that 
$\pif$ acts trivially. If $F=\QQ_p$ then using the results of Breuil we may show that our construction gives all the supersingular representations 
up to a twist by an unramified quasi-character.  We conjecture that this is true for arbitrary $F$. If $\rho$ is an irreducible smooth 
$\Fbar$-representation of the Weil group $W_F$ of $F$, then the wild inertia subgroup of $W_F$ acts trivially on $\rho$, since it is pro-$p$ and
normal in $W_F$. This implies that  there are only $q(q-1)/2$ isomorphism classes of irreducible smooth $2$-dimensional 
$\Fbar$-representations $\rho$ of the Weil group  of $F$ such that $(\det\rho)(\Fr)=1$. Here, $\Fr$ is the Frobenius automorphism corresponding
to $\pif$ via the local class field theory. So the conjecture would be true if there was a Langlands type of correspondence. 

The starting point in this theory is that every pro-$p$ group acting smoothly on an $\Fbar$-vector has a  non-zero invariant vector. Let 
$I_1$ be the unique maximal pro-$p$ subgroup of the standard Iwahori subgroup $I$ of $G$. Given a smooth representation $\pi$ of $G$ the 
Hecke algebra $\HH=\End_G(\cIndu{I_1}{G}{\Eins})$ acts on the $I_1$-invariants $\pi^{I_1}$. It is expected 
that this functor induces a bijection between the irreducible smooth representations of $G$ and the irreducible modules of $\HH$. This happens if 
$F=\QQ_p$. Moreover, if $F$ is arbitrary and  $\pi$ is an irreducible smooth representation of $G$, which is not supersingular, then $\pi^{I_1}$ 
is an irreducible $\HH$-module. All the irreducible modules of $\HH$ that do not arise this way are called supersingular. They have been determined 
by Vign\'eras and we give a list of them in the Definition \ref{list}. There are $q(q-1)/2$ isomorphism classes of irreducible supersingular 
modules of $\HH$ up to a twist by an unramified quasi-character.

Given a supersingular module $M$ of $\HH$ we construct two $G$-equivariant coefficient systems $\VV$ and $\II$ 
on the Bruhat-Tits tree $X$ of $\PGL_2(F)$ and a morphism of $G$-equivariant coefficient systems between them. Once we pass to the $0$-th homology, 
this induces a homomorphism of $G$-representations. We show that the image of this homomorphism 
$$\pi=\Image(H_0(X,\VV)\rightarrow H_0(X,\II))$$
is a smooth irreducible representation of $G$, which is supersingular, since $\pi^{I_1}$ contains a supersingular module $M$. 
Moreover, we show that two non-isomorphic irreducible supersingular modules give rise to non-isomorphic representations. However, the question 
of determining all smooth irreducible representations $\pi$ of $G$, such that $\pi^{I_1}$ contains $M$, remains open.

We will describe the contents of this  paper in more detail. In Section \ref{one} we recall the algebra structure of $\HH$ and 
the definition of  supersingular modules. 

Sections \ref{GLFq} and \ref{princindrep} deal with some aspects of the $\Fbar$-representation theory of $\Gamma=\GL_2(\Fq)$. 
In Section \ref{GLFq} we give two different descriptions of the irreducible $\Fbar$-representations of $\Gamma$, one of them 
due to Carter and Lusztig \cite{cl} and the other one due to Brauer and Nesbitt \cite{brauer}, and a dictionary between them. Let $U$ be the subgroup
of unipotent upper-triangular matrices in $\Gamma$, then $U$ is a $p$-Sylow subgroup of $\Gamma$. If $\rho$ is a representation
of $\Gamma$, then the Hecke algebra $\HH_{\Gamma}=\End_{\Gamma}(\Indu{U}{\Gamma}{\Eins})$ acts on the $U$-invariants $\rho^{U}$. This functor
induces a bijection between the irreducible representations of $\Gamma$ and the irreducible right modules of $\HH_{\Gamma}$. 

Every  representation $\rho$ of $\Gamma$  has an injective envelope $\iota:\rho\hookrightarrow \inj \rho$. By this we mean, a representation $\inj \rho$ of 
$\Gamma$ and an injection $\iota$, such that $\inj \rho$ is an injective object in the category of $\Fbar$-representations of $\Gamma$ and every non-zero 
$\Gamma$-invariant subspace of $\inj \rho$ intersects $\iota(\rho)$ non-trivially. Injective envelopes are unique up to isomorphism. 
In Section \ref{princindrep} we determine the $\HH_{\Gamma}$-module structure of $(\inj \rho)^{U}$, for an irreducible representation 
$\rho$ of $\Gamma$. This is important to us, so we give  two ways of doing it. If $p=q$ then the dimension of $(\inj \rho)^{U}$ is small and 
this enables us to give an elementary argument. In general we use the results of Jeyakumar \cite{j}, where he describes explicitly 
injective envelopes of irreducible representations of $\SL_2(\Fq)$. 

Let $\oF$ be the ring of integers of $F$, let $K=\GL_2(\oF)$. The reduction modulo the prime ideal of $\oF$ induces a  surjection
$K\rightarrow \Gamma$, let $K_1$ be the kernel of this map. The Hecke algebra $\HH_{K}=\End_{K}(\Indu{I_1}{K}{\Eins})$ is naturally a 
subalgebra of $\HH$. Let $M$ be a supersingular module of $\HH$, 
then the restriction of $M$ to $\HH_{K}$ is isomorphic to a direct sum of two irreducible modules of $\HH_{K}$. Since $K/K_1\cong \Gamma$ 
we may identify representations of $K$ on which $K_1$ acts trivially with the representations of $\Gamma$. This induces an identification 
$\HH_K=\HH_{\Gamma}$. Since the irreducible modules of $\HH_{\Gamma}$
are in bijection with the irreducible representations of $\Gamma$, there exists a unique representation $\rho=\rho_M$ of $\Gamma$, such that 
$\rho$ is isomorphic to a direct sum of two irreducible representations of $\Gamma$, and $\rho^{U}\cong M|_{\HH_{\Gamma}}$.
Let $\rho\hookrightarrow \inj \rho$ be an injective 
envelope of $\rho$ in the category of $\Fbar$-representations of $\Gamma$. We consider now both $\rho$ and $\inj \rho$ as representations of $K$.
 We have an exact sequence 
\begin{displaymath}
\xymatrix{0\ar[r]&\rho^{I_1}\ar[r]& (\inj \rho)^{I_1}}
\end{displaymath}
of $\HH_K$-modules. The main result of Section \ref{princindrep} are Propositions \ref{injqp} ($p=q$), Propositions \ref{injiwahori} and  
\ref{injregular} (general case),
which say that there exists an action of $\HH$, extending the action of $\HH_K$, on $(\inj \rho_{\gamma})^{I_1}$, such that the above exact sequence
yields an exact sequence
\begin{displaymath}\tag{E}\label{EA}
 \xymatrix{0\ar[r]&M\ar[r]& (\inj \rho)^{I_1}}
\end{displaymath}
of $\HH$-modules. The fact that we can extend the action and obtain (\ref{EA}) implies the existence of a certain $G$-equivariant coefficient 
system $\II$ on the tree $X$. 

The inspiration to use coefficient systems comes from the works of Schneider and Stuhler \cite{ss1} and \cite{ss2}, where the authors work over
the complex numbers, and Ronan and Smith \cite{rs}, where the $\Fbar$ coefficient systems are studied for finite Chevalley groups. We introduce coefficient
systems  in Section \ref{coeffsys}. Let $\sigma_1$ be an edge on $X$ containing a vertex $\sigma_0$. Since, $G$ acts transitively on the vertices of 
the tree $X$, the category of $G$-equivariant coefficient systems
is  equivalent to a category of diagrams $\diag$. The  objects of $\diag$ are triples $(\rho_0,\rho_1,\phi)$, where $\rho_0$ is a smooth representation 
of $\kzr$, $\rho_1$ is a smooth representation of $\kon$ and $\phi$ is a $\kon\cap \kzr$-equivariant homomorphism, $\phi:\rho_1\rightarrow \rho_0$,
where $\kzr$ and $\kon$ are the $G$-stabilisers of $\sigma_0$ and $\sigma_1$.
 The proof of equivalence between the two categories is the main result of Section \ref{coeffsys}. As a corollary  we obtain a nice way of passing 
from ``local'' to ``global'' information, see Corollary \ref{thepoint}, and we use this in the construction of $\II$.  

More precisely, we start with a supersingular $\HH$-module $M$ and find the unique smooth representation $\rho=\rho_M$ of $K$, such that $\rho$ is 
isomorphic to a direct sum of two irreducible representations of $K$, and $\rho^{I_1}\cong M|_{\HH_K}$, as above. We then consider an injective envelope
 $\rho\hookrightarrow \Inj \rho$ of $\rho$ in the category of smooth $\Fbar$-representations of $K$. Let $\sigma_1$ be an edge on $X$ fixed by $I$ 
and let $\sigma_0$ be a vertex fixed by $K$. We extend the action of $K$ on $\Inj \rho$ 
to the action of $F^{\times}K=\kzr$, so that a fixed uniformiser acts trivially. We denote this  representation  by $Y_{0}$. Let us assume that
we may extend the action of $F^{\times}I=\kon\cap\kzr$ on $Y_{0}|_{F^{\times}I}$ to the action of $\kon$. We denote the corresponding 
representation of $\kon$ by $Y_{1}$. The triple $Y=(Y_{0}, Y_{1}, \id)$ is an object in a category $\diag$, which is 
equivalent to the category of $G$-equivariant coefficient systems on the tree $X$, by the main result of  
Section \ref{coeffsys}. So $Y$ gives us a $G$-equivariant coefficient system $\II$. Moreover, the restriction maps of $\II$ are all isomorphisms. This
implies that 
$$H_0(X,\II)|_K\cong \Inj \rho.$$ 
In particular, we have an injection
$$\rho\hookrightarrow \Inj \rho\cong H_0(X,\II_{\gamma})|_K,$$
which gives us an exact sequence of vector spaces 
\begin{displaymath}
\xymatrix{0\ar[r]&\rho^{I_1}\ar[r]& H_0(X,\II)^{I_1}}.
\end{displaymath}
We show in Subsection \ref{II} that using (\ref{EA}) we may extend the action of $F^{\times}I$ on $Y_{0}|_{F^{\times}I}$ to the action of 
$\kon$, so that the image of $\rho^{I_1}$ in $H_0(X,\II)^{I_1}$ is an $\HH$-invariant subspace, isomorphic to $M$ as an $\HH$-module. We let 
$\pi$ be the $G$-invariant subspace of $H_0(X,\II)$ generated by the image of $\rho$. 
In Theorem \ref{irrconst} we prove that $\pi$ is irreducible and supersingular. We also show that $\pi$ is the socle of $H_0(X,\II)$. 
The space $H_0(X,\II)^{I_1}$ is always finite dimensional, we determine the $\HH$-module structure in  Proposition \ref{Injmod}. 
The proofs rely on some general properties of injective envelopes, which we recall in Subsection \ref{injectenve}. Using injective envelopes 
we also give a new proof of the criterion for admissibility of a smooth representation of $G$, which works in a very general context, see Subsection 
\ref{admissibility}. 

We would like to explain the thinking behind the construction of the coefficient system $\VV$ in Subsection \ref{VV}. Let $\pi$ be a 
smooth representation of $G$, generated by its $I_1$-invariant vectors. We may associate to $\pi$ a $G$-equivariant 
coefficient system $\FF_{\pi}$ as follows.
Given a simplex $\sigma$ on $X$, we let $U^{1}_{\sigma}$ be the maximal normal pro-$p$ subgroup of the $G$-stabiliser of $\sigma$. With this notation
$U^{1}_{\sigma_1}=I_1$ and $U^{1}_{\sigma_0}=K_1$. We may 
consider the coefficient system of invariants $\FF_{\pi}=(\pi^{U^{1}_{\sigma}})_{\sigma}$, where the restriction maps are inclusions. Since 
$\pi$ is generated by its $I_1$-invariants the natural map
$$H_0(X,\FF_{\pi})\rightarrow \pi$$
is surjective. If we are working over the complex numbers then a theorem of Schneider and Stuhler in \cite{ss1}, says that the above homomorphism 
is in fact an isomorphism. If we are working over $\Fbar$, then $H_0(X,\FF_{\pi})$ can be much bigger than $\pi$.

The construction of $\VV$ is motivated by the following question. Let $M$ be a supersingular module of $\HH$ and suppose that there 
exists a smooth irreducible $\Fbar$-representation $\pi$ of $G$ such that $\pi^{I_1}\cong M$. What can be said about the corresponding coefficient system 
$\FF_{\pi}$? It is enough  to understand the action of $K$ on $\pi^{K_1}$. 
This reduces the question to the representation theory of $\GL_2(\Fq)$.
In Corollary \ref{shouldbeall} we show that there exists an injection $\VV\hookrightarrow \FF_{\pi}$  
and hence every $\pi$ as above is a quotient of $H_0(X,\VV)$. We would 
like to point out that although in most cases we do not know whether such $\pi$ exists, the coefficient system $\VV$ is always
well defined. Moreover, if $\pi$ is any non-zero irreducible quotient of $H_0(X,\VV)$, then we show that $\pi$ is supersingular, 
since $\pi^{I_1}$ contains a supersingular $\HH$-module $M$. This implies that $H_0(X,\VV)$ is a quotient of one of 
the spaces considered by Barthel and Livn\'e in \cite{bl2}. Corollary \ref{allnice} implies that at least in some cases the quotient map is an isomorphism.
Now the Remarque 4.2.6 in \cite{breuil} shows that in general $\dim H_0(X,\VV)^{I_1}$ is infinite. The irreducible representation $\pi$, which we construct
in this paper, is a quotient of $H_0(X,\VV)$, moreover the space $\pi^{I_1}$ is finite dimensional. Hence, in contrast to 
the situation over $\CCC$, in general $H_0(X,\VV)$ is very far away from being irreducible. 

We believe that our construction of irreducible representations will work for other groups. Our strategy could be applied most directly to
the group $G=\GL_N(F)$, where $N$ is a prime number. If $N$ is prime then the maximal open, compact-mod-centre subgroups of $G$ are 
the $G$-stabilisers of chambers (simplices of maximal dimension) and vertices in the Bruhat-Tits building of $G$ and if 
we had the equivalent of (\ref{EA}) then the construction of the coefficient system $\II$ and our proofs would carry through. 
However, in order to do this one needs to understand the $\HH_{\Gamma}$-module structure of $(\inj \rho)^U$, (or at least the action of $B$
on $(\inj \rho)^{U}$, at the cost of not knowing $\HH$-module structure of $H_0(X,\II)^{I_1}$), where $\rho$ is 
an irreducible $\Fbar$-representation of $\Gamma=\GL_N(\Fq)$, $B$ is the subgroup of upper-triangular matrices, and $U$ is the subgroup of unipotent 
upper-triangular matrices of $\Gamma$. This might be quite a difficult problem, since already for $N=2$ the dimension of $(\inj \rho)^U$ can be as 
big as $2^n-1$, if $q=p^n$.
  
\textit{Acknowledgements.} I would like to thank Michael Spiess and Thomas Zink for a number of useful discussions and for
looking after me in general. I would like to thank Marie-France Vign\'eras for her encouragement and her comments on this work.
\subsection{Notation}
Let $F$ be a non-Archimedean local field, $\oF$ its ring of integers, $\pF$ the maximal ideal of $\oF$. Let $p$ be the characteristic and
let $q$ be the number of elements of the residue class field of $F$. We fix a uniformiser $\pif$ of $F$. 

Let $G=\GL_2(F)$ and  $K=\GL_2(\oF)$. Reduction modulo $\pF$ induces a surjective homomorphism 
$$\red: K\rightarrow \Gamma=\GL_2(\Fq).$$
Let $K_1$ be the kernel of $\red$. Let $B$ be the subgroup of $\Gamma$ of upper triangular matrices. Then
$$B=HU$$
where $H$ is the subgroup of diagonal matrices and $U$ is the subgroup of unipotent matrices in $B$. It is of importance, that the order of
$H$ is prime to $p$ and $U$ is a $p$-Sylow subgroup of $\Gamma$. Let $I$ and $I_1$ be the subgroups of $K$, given by
$$I=\red^{-1}(B), \quad I_1=\red^{-1}(U).$$  
Then $I$ is the Iwahori subgroup of $G$ and $I_1$ is the unique maximal pro-$p$ subgroup of $I$. Let $T$ be the subgroup of diagonal matrices
in $K$, and let $T_1=T\cap K_1=T\cap I_1$. Let $N$ be the normaliser of $T$ in $G$. We introduce some special elements of $N$. Let
$$\Pi=\begin{pmatrix}  0 & 1\\ \pif & 0 \end{pmatrix},\quad  n_s=\begin{pmatrix}  0 & -1\\ 1 & 0 \end{pmatrix},
\quad s=\begin{pmatrix}  0 & 1\\ 1 & 0 \end{pmatrix}. $$
The images of $\Pi$ and $n_s$ in $N/T$, generate it as a group. The normaliser $N$ acts on $T$ by conjugation, and hence it acts on the group of characters
of $T$. This action factors through $T$, so if $w\in N/T$ and $\chi$ is a character of $T$, we will write $\chi^{w}$ for the 
character, given by
$$\chi^{w}(t)= \chi(w^{-1}t w), \quad \forall t\in T.$$ 
Let $\tilde{B}$ be the group of upper-triangular matrices in $G$, then
$\tilde{B}= \tilde{T}\tilde{U}$
where $\tilde{T}$ is the group of diagonal matrices in $G$ and $\tilde{U}$ is the group of unipotent matrices in $\tilde{B}$.

\begin{defi}\label{superdefi} Let $\pi$ be a smooth irreducible $\Fbar$-representation of $G$, such that  $\pi$ admits
 a central character, then $\pi$ is called supersingular if $\pi$ is not a subquotient of 
$\Indu{\tilde{B}}{G}{\chi}$, for any smooth quasi-character $\chi: \tilde{B}\rightarrow \tilde{B}/\tilde{U}\cong \tilde{T}\rightarrow \Fbar^{\times}$.
\end{defi}

All the representations considered in this paper are over $\Fbar$, unless it is stated otherwise.

\section{Hecke algebra}\label{one}

\begin{lem}\label{propinvariants} Let $\mathcal P$ be a pro-$p$ group and let $\pi$ be a smooth non-zero representation of $\mathcal P$, then 
the space $\pi^{\mathcal P}$ of $\mathcal P$-invariants is non-zero.
\end{lem}
\begin{proof} We choose a non-zero vector $v$ in $\pi$. Let 
$\rho=\langle \mathcal P v\rangle_{\Fbar}$
 be a subspace of $\pi$ generated by $\mathcal P$ and $v$. Since the action of $\mathcal P$ on $\pi$ is smooth, the stabiliser 
$\Stab_{\mathcal P}(v)$ has finite index in $\mathcal P$, hence $ \rho$ is finite dimensional. Let $v_1,\ldots, v_{d}$ be an 
$\Fbar$ basis of $\rho$. The group $\mathcal P$ acts on $\rho$ and the kernel of this action is given by
$$\Ker \rho= \bigcap_{i=1}^{d} \Stab_{\mathcal P}(v_i).$$ 
In particular, $\Ker \rho$ is an open subgroup of $\mathcal P$. Hence, $\mathcal P/ \Ker \rho $ is a finite group, whose order is a 
power of $p$. Now,
$$\rho^{\mathcal P}=\rho^{\mathcal P/\Ker \rho}\neq 0$$
since $\mathcal P/\Ker \rho$ is a finite $p$-group, see \cite{serre}, \S 8, Proposition 26.  
\end{proof}

Let $\pi$ be a smooth representation of $G$, then 
$$\pi^{I_1}\cong \Hom_{I_1}(\Eins, \pi)\cong \Hom_{G}(\cIndu{I_1}{G}{\Eins},\pi)$$
by Frobenius reciprocity. Let $\HH$ be the Hecke algebra
$$\HH=\End_{G}(\cIndu{I_1}{G}{\Eins})$$
then via the above isomorphism $\pi^{I_1}$ becomes naturally a right $\HH$-module.  We obtain a functor
$$\Rep_G\rightarrow \Mod-\HH,\quad \pi\mapsto \pi^{I_1},$$
where $\Rep_{G}$ is a category of smooth $\Fbar$-representations of 
$G$ and $\Mod-\HH$ is the category of right $\HH$-modules. Since $I_1$ is an open pro-$p$ subgroup of $G$, Lemma \ref{propinvariants} implies 
that $\pi^{I_1}=0$ if and only if $\pi=0$. This functor is our basic tool. We want to study the structure of $\HH$. We follow \cite{cl}, 
where finite groups with split $BN$-pair are treated, a lot of the proofs just carry over formally.

\begin{defi}\label{defTn} Let $g\in G$ and $f\in \cIndu{I_1}{G}{\Eins}$ we define  $T_{g}\in \HH$ by  
$$(T_{g}f)(I_1 g_1)=\sum_{I_1 g_2 \subseteq I_1 g^{-1} I_1  g_{1}} f(I_1 g_2).$$
\end{defi}
\begin{lem}\label{iwahoridecomp} We may write $G$ as a disjoint union
$$G=\mathop{\/\dot{\bigcup}}\limits_{n\in N/T_1}I_1 n I_1$$
of double cosets.
\end{lem}
\begin{proof}This follows from the Iwahori decomposition.
\end{proof}

It is immediate that the definition of $T_g$ depends only on the double coset $I_1gI_1$. The Lemma above implies that
it is enough to consider $T_n$, where $n\in N$ is a representative of a coset in $N/T_1$.

\begin{defi} Let $\varphi\in \cIndu{I_1}{G}{\Eins}$ be the unique function such that
$$\supp \varphi = I_1\quad \textrm{and}\quad \varphi(u)=1, \quad \forall u\in I_1.$$ 
\end{defi}

\begin{lem}\label{basisTinv} 
\begin{itemize} 
\item[(i)] The function $\varphi$ generates $\cIndu{I_1}{G}{\Eins}$ as a $G$-rep\-re\-sen\-ta\-tion. 
\item[(ii)] $\supp T_n \varphi=I_1n I_1$  and $(T_n\varphi)(g)=1$, for every $g\in I_1 n I_1$. In particular,
$$T_n\varphi=\sum_{u\in I_1/(I_1\cap n^{-1}I_1 n)} u n^{-1}\varphi.$$ 
\item[(iii)] The set $\{T_n \varphi: n\in N/T_1\}$ is an $\Fbar$-basis of $(\cIndu{I_1}{G}{\Eins})^{I_1}$.
\item[(iv)] The set $\{T_n: n\in N/T_1\}$ is an $\Fbar$-basis of $\HH$.
\end{itemize}
\end{lem}
\begin{proof} Let $g\in G$, then $\supp (g^{-1}\varphi) =I_1 g$ and $(g^{-1}\varphi)(I_1 g)=1$. Part (i) follows immediately. 

Let $f\in \cIndu{I_1}{G}{\Eins}$, then by examining the definition of $T_n$, one obtains that $\supp (T_n f)\subseteq I_1 n\supp f$. Hence, 
$\supp (T_n \varphi)\subseteq I_1 n I_1$. Since $T_n$ is a $G$-equivariant homomorphism and $I_1$ acts trivially on $\varphi$, it is enough to prove that 
$(T_n \varphi)(n)=1$. Since $\supp \varphi =I_1$, it is immediate from  Definition \ref{defTn} that $(T_n \varphi)(I_1n)=\varphi(I_1)=1$. The last part 
follows from decomposing $I_1 n I_1$ into right cosets and applying the argument used in Part (i). 

Let $n, n'\in N$, and suppose that $nT_1\neq n' T_1$, then Lemma \ref{iwahoridecomp} implies that $I_1 n I_1\neq I_1 n' I_1$. By Part (ii) the functions 
$T_n\varphi$ and $T_{n'}\varphi$ have disjoint support. This implies that the set $\{T_n \varphi: n\in N/T_1\}$ is linearly independent.
Any $f\in (\cIndu{I_1}{G}{\Eins})^{I_1}$, is constant on the double cosets $I_1 n I_1$, for $n\in N$, and since $\supp f$ is compact, $f$ is supported
only on finitely many such, hence Lemma \ref{iwahoridecomp} and Part (ii) imply that $\{T_n \varphi: n\in N/T_1\}$ is also a spanning set.  
Hence we get Part (iii).

Let $\psi\in \HH$, Part (i) implies that $\psi=0$ if and only if $\psi(\varphi)=0$. This observation coupled with Part (iii) implies Part (iv).
\end{proof} 
\begin{cor}\label{T} Let $\pi$ be a smooth representation of $G$ and let $v\in \pi^{I_1}$, then the action of $T_n$ on $\pi^{I_1}$ is given by 
$$v T_n= \sum_{u\in I_1/(I_1\cap n^{-1}I_1 n)} u n^{-1}v.$$
\end{cor}
\begin{proof}The isomorphism $\Hom_{G}(\cIndu{I_1}{G}{\Eins},\pi)\cong \pi^{I_1}$ is given explicitly by $\psi \mapsto \psi(\varphi)$. 
Let $\psi$ be the 
unique $G$-invariant homomorphism, such that $\psi(\varphi)=v$, then 
$$vT_{n_s}=(\psi\circ T_{n_s})(\varphi)= \psi(T_{n_s}\varphi)=\psi(\sum_{u\in I_1/(I_1\cap n^{-1}I_1 n)} u n^{-1}\varphi).$$
The last equality follows from Lemma \ref{basisTinv} (ii). Since, $\psi$ is $G$-invariant, we obtain the Lemma.
\end{proof}  

\begin{lem}\label{noco} Let $n', n\in N$ and suppose that $n$ normalises $I_1$, then $$T_{n'}T_n= T_{n' n}, \quad T_n T_{n'}= T_{n n'}.$$
\end{lem}
\begin{proof} Lemma \ref{basisTinv} (i) implies that it is enough to show that the homomorphisms map $\varphi$ to the same function. Let 
$f\in \cIndu{I_1}{G}{\Eins}$ then since $n$ normalises $I_1$ we have $(T_n(f))(g)= f(ng)$ and $T_n\varphi= n^{-1}\varphi$. Now the Lemma follows from 
Lemma \ref{basisTinv} (ii).
\end{proof}

Let $t\in T$ and let $h$ be the image of $t$ in $H$, via $T/T_1\cong H$, we will write $T_h$ for the homomorphism $T_t$. 
\begin{defi}\label{simide}Let $\chi: H\rightarrow \Fbar^{\times}$ be a character, we define 
$$e_{\chi}=\frac{1}{|H|}\sum_{h\in H} \chi(h)T_h.$$
Let $$\varphi_{\chi}=e_{\chi}\varphi,$$ then $\varphi_{\chi}$ is the unique function in $\cIndu{I_1}{G}{\Eins}$ such that 
$$\supp \varphi_{\chi}=I,\quad \varphi_{\chi}(g)= \chi(gI_1),\quad \forall g\in I,$$
 via the isomorphism $I/I_1\cong H$.
\end{defi}
\begin{lem}\label{idem}
\begin{itemize}
\item[(i)] $e_{\chi}^{2}=e_{\chi}$ and $e_{\chi}e_{\chi'}=0$, if $\chi\neq \chi'$.
\item[(ii)]$\id =\sum_{\chi} e_{\chi}$, where the sum is taken over all characters $\chi :H\rightarrow \Fbar^{\times}$.
\item[(iii)]$e_{\chi}(\cIndu{I_1}{G}{\Eins})\cong \cIndu{I}{G}{\chi}$.
\end{itemize}
\end{lem}
\begin{proof} We note that $H$ is abelian and the order of $H$ is prime to $p$. Parts (i) and (ii) follow from the orthogonality relations of characters.
Lemma \ref{basisTinv} (i) implies that $e_{\chi}(\cIndu{I_1}{G}{\Eins})$ is generated by $\varphi_{\chi}$ and this implies Part (iii).
\end{proof}
\begin{cor}Let $\pi$ be a smooth representation of $G$, then $I$ acts on $(\pi^{I_1})e_{\chi}$ by a character $\chi$. Moreover,
$$\pi^{I_1}\cong\oplus_{\chi}(\pi^{I_1})e_{\chi}.$$
\end{cor}
\begin{proof}The group $I$ acts on $\pi^{I_1}$. Since $I_1$ acts trivially and $I/I_1\cong H$, which is abelian and of order  prime to $p$, 
the space $\pi^{I_1}$ decomposes into one dimensional $I$ invariant subspaces. Corollary \ref{T} implies that $e_{\chi}$ cuts out the $\chi$-isotypical 
subspace. The last part follows from Lemma \ref{idem} (ii).   
\end{proof}
\begin{lem}\label{relation}
\begin{itemize}
\item[(i)] $T_{n_s}e_{\chi}= e_{\chi^{s}}T_{n_s}$,\quad $T_{\Pi}e_{\chi}=e_{\chi^{s}}T_{\Pi}$.
\item[(ii)] If $\chi=\chi^{s}$ then 
$$T^2_{n_s} e_{\chi}=- T_{n_s}e_{\chi}.$$
If $\chi\neq \chi^{s}$ then
$$T^{2}_{n_s}e_{\chi}=0.$$
\end{itemize}
\end{lem}
\begin{proof}Part (i) follows from Lemma \ref{noco}. Lemma \ref{basisTinv} (i) implies that it is enough to calculate 
$T_{n_s}^{2}e_{\chi} \varphi=T_{n_s}^2 \varphi_{\chi}$. Applying Lemma \ref{basisTinv} (ii) twice we obtain
$\supp T_{n_s}^2\varphi_{\chi}\subseteq K$. Hence it is enough to do the calculation in the space $\Indu{I_1}{K}{\Eins}$. Since $K_1$ 
acts trivially on this
space, it is enough to do the calculation in the space $\Indu{U}{\Gamma}{\Eins}$. Then the Lemma is a special case of \cite{cl} Theorem 4.4.  
\end{proof} 
\begin{lem}\label{pinsm}Let $m\ge 0$  and let $w=\Pi n_s$ then the following hold:
\begin{itemize}
\item[(i)] $I_1 w I_1 w^m I_1 = I_1 w^{m+1} I_1$,
\item[(ii)] $I_1w^{-1}I_1w^{m+1} \cap I_1 w^m I_1= I_1 w^m$,
\item[(iii)] $T_{w^m} = (T_w)^m = (T_{\Pi} T_{n_s})^m$.
\end{itemize}
\end{lem} 
\begin{proof}The first two parts can be checked by a direct calculation. For Part (iii) we observe that 
$$ \supp T_{w} T_{w^m}\varphi \subseteq I_1 w \supp T_{w^m}\varphi=I_1 w I_1 w^{m} I_1= I_1 w^{m+1}I_1,$$ 
where the last equality is Part (i). Part (ii) and Lemma \ref{basisTinv} (ii) imply that
$$(T_w T_{w^{m}}\varphi)(w^{m+1})=1.$$
 Since $I_1$ acts trivially on $\varphi$ and all the homomorphisms are $G$-equivariant, we may apply Lemma \ref{basisTinv} (ii) again to obtain
$$T_{w} T_{w^m}\varphi= T_{w^{m+1}}\varphi.$$
Lemma \ref{basisTinv} (i) implies that $T_w T_{w^{m}}= T_{w^{m+1}}$. Induction and Lemma \ref{noco} gives us Part (iii).
\end{proof}
\begin{lem}\label{generelat}
\begin{itemize}
\item[(i)]Let $n\in N$, then there exists $h\in H$ and integers $a\in\{0,1\}$, $m\ge 0$ and $b\in \ZZ$ such that
$$T_n= T_{\Pi}^{a}(T_{\Pi}T_{n_s})^{m} T_{\Pi}^{b}T_h$$
where $T_{\Pi}^{-1}=T_{\Pi^{-1}}$.
\item[(ii)]The elements $T_{n_s}$, $T_{\Pi}$, $T_{\Pi^{-1}}$ and $e_{\chi}$, for every character $\chi:H\rightarrow \Fbar^{\times}$,  generate
$\HH$ as an algebra.
\end{itemize}
\end{lem} 
\begin{proof}We note that Lemma \ref{noco} implies that $T_{\Pi}$ is invertible with $T_{\Pi}^{-1}= T_{\Pi^{-1}}$ and $T_{\Pi}^{2}$ is central in $\HH$.
Every $n\in N$ maybe written as $n=\Pi^{a}(\Pi n_s)^{m}\Pi^b t$, where $t\in T$. Lemma \ref{noco} and Lemma \ref{pinsm}(iii) imply Part (i). 
Hence $T_{n_s}$, $T_{\Pi}$, $T_{\Pi^{-1}}$ and  $T_h$, for $h\in H$ generate $\HH$ as an algebra. Lemma \ref{noco} implies that 
$T_h e_{\chi}= \chi(h^{-1})e_{\chi}$ and hence Lemma \ref{idem} (ii) implies that $T_h$ can be expressed as a linear combination of idempotents 
$e_{\chi}$.
This gives us Part (ii).
\end{proof}

\begin{lem}\label{basisechi} 
\begin{itemize} 
\item[(i)] The set $\{e_{\chi}T_n \varphi: n\in N/T, \chi:H\rightarrow \Fbar^{\times}\}$ is an $\Fbar$-basis of $(\cIndu{I_1}{G}{\Eins})^{I_1}$.
\item[(ii)] The set $\{e_{\chi}T_n: n\in N/T, \chi:H\rightarrow \Fbar^{\times}\}$ is an $\Fbar$-basis of $\HH$.
\end{itemize}
\end{lem}
\begin{proof}Since $e_{\chi}T_h = \chi(h^{-1})e_{\chi}$ Lemma \ref{basisTinv} (iii) implies that  the set 
$\{e_{\chi}T_n \varphi: n\in N/T, \chi:H\rightarrow \Fbar^{\times}\}$ is a spanning set. Since the elements $e_{\chi}$ are orthogonal idempotents 
it is enough to show that the set $\{e_{\chi}T_n\varphi: n\in N/T\}$ is linearly independent for a fixed character $\chi$. Lemma \ref{basisTinv} (ii) 
implies that $\supp e_{\chi}T_n \varphi =I n I$. Lemma \ref{iwahoridecomp} implies that if $nT\neq n'T$, then $e_{\chi}T_n\varphi$ and 
$e_{\chi}T_{n'}\varphi$ have disjoint support and hence the set is linearly independent. Part (ii) follows from Part (i) and Lemma \ref{basisTinv} (i).
\end{proof}

\subsection{Supersingular modules }
All the irreducible modules of $\HH$ have been determined by Vign\'eras in \cite{vig}. They naturally split up into two classes.

\begin{prop}  
 Let $\pi$ be a smooth irreducible representation of $G$, which admits a central character. Suppose that $\pi$ is not supersingular, 
then $\pi^{I_1}$ is an irreducible $\HH$-module.
\end{prop}
\begin{proof}See \cite{vig} E.5.1.
\end{proof}
The modules as above could be called non-supersingular, we are interested in all the rest. 
\begin{defi}\label{list}Let $\chi:H\rightarrow \Fbar^{\times}$ be a character, let $\gamma=\{\chi, \chi^{s}\}$ and let $\lambda\in \Fbar^{\times}$.
We define a standard  supersingular module $M_{\gamma}^{\lambda}$ to be a right $\HH$-module such that its underlying vector space is $2$ dimensional
$$M_{\gamma}^{\lambda}=\langle v_1, v_2\rangle_{\Fbar}$$
and  the action of $\HH$ is determined by the following:
\begin{itemize}
\item[(i)]If $\chi=\chi^{s}$ then
$$v_1 e_{\chi}=v_1, \quad v_1 T_{n_{s}}=-v_1,\quad v_1 T_{\Pi}= v_2$$ 
and 
$$v_2 e_{\chi}=v_2, \quad v_2 T_{n_{s}}=0,\quad v_2 T_{\Pi}= \lambda v_1.$$
\item[(ii)]If $\chi\neq \chi^{s}$ then
$$v_1 e_{\chi}=v_1, \quad v_1 T_{n_{s}}=0,\quad v_1 T_{\Pi}= v_2$$ 
and 
$$v_2 e_{\chi^{s}}=v_2, \quad v_2 T_{n_{s}}=0,\quad v_2 T_{\Pi}= \lambda v_1.$$
\end{itemize}
\end{defi}

To show that these relations define an action of $\HH$ requires some work, this is done in \cite{vig}.

\begin{lem}\label{gammagamma}The modules $M_{\gamma}^{\lambda}$ are irreducible and  
$$M_{\gamma'}^{\lambda'}\cong M_{\gamma}^{\lambda}$$
if and only if  $\gamma'= \gamma$ and $\lambda'=\lambda$.
\end{lem}
\begin{proof}The definition immediately gives that $M_{\gamma}^{\lambda}$ does not have a $1$ dimensional submodule, hence it is irreducible.
 If $\chi':H\rightarrow \Fbar^{\times}$ is a character, such that $\chi'\not\in \gamma$ then 
$$M_{\gamma}^{\lambda}e_{\chi'}=0.$$
Hence, $\gamma=\gamma'$. The element $T_{\Pi}^{2}$ acts on $M_{\gamma}^{\lambda}$ by a scalar $\lambda$. Hence, $\lambda=\lambda'$. 
\end{proof}

 The following Proposition explains why $M_{\gamma}^{\lambda}$ are called supersingular.
\begin{prop}
Let $M$ be an irreducible $\HH$ module, such that $M\not\cong \pi^{I_1}$ for any non-\-su\-per\-sin\-gu\-lar irreducible representation $\pi$, then
$$M\cong M_{\gamma}^{\lambda}$$
for some $\gamma$ and $\lambda$. 
\end{prop}
\begin{proof}See \cite{vig} C.2 and E.5.1.  
\end{proof}

\begin{cor}\label{supersingularity} Let $\pi$ be a smooth irreducible representation of $G$, admitting a central character. Suppose that 
$\pi^{I_1}$ contains a  submodule isomorphic to $M_{\gamma}^{\lambda}$ for some $\gamma$ and $\lambda$, then $\pi$ is supersingular.
\end{cor}

We will also need to consider the following extension of supersingular modules.
\begin{defi}  Let $\chi:H\rightarrow \Fbar^{\times}$ be a character, such that $\chi\neq \chi^{s}$, let $\gamma=\{\chi, \chi^{s}\}$ and let 
$\lambda\in \Fbar^{\times}$. Let
$$\HH^{\lambda}= \HH/(T_{\Pi}^{2}-\lambda)\HH$$
then we define a right $\HH$-module $L^{\lambda}_{\gamma}$ to be
$$L_{\gamma}^{\lambda}= e_{\chi}\HH^{\lambda}/e_{\chi}(T_{\Pi}T_{n_s}-T_{n_s}T_{\Pi}) \HH^{\lambda}.$$
\end{defi} 
The definition seems to be asymmetric in $\chi$ and $\chi^{s}$, however the multiplication from the left by $T_{\Pi}$ induces an isomorphism
$$e_{\chi}\HH^{\lambda}/e_{\chi}(T_{\Pi}T_{n_s}-T_{n_s}T_{\Pi}) \HH^{\lambda}\cong e_{\chi^{s}}\HH^{\lambda}/
e_{\chi^{s}}(T_{\Pi}T_{n_s}-T_{n_s}T_{\Pi}) \HH^{\lambda},$$
since $T_{\Pi}$ is a unit in $\HH^{\lambda}$.
\begin{lem}\label{lgammabasis} The images of $e_{\chi}$, $e_{\chi}T_{\Pi}$, $e_{\chi}T_{n_s}$ and $e_{\chi}T_{n_s}T_{\Pi}$ in $L^{\lambda}_{\gamma}$
  form an $\Fbar$-basis of $L^{\lambda}_{\gamma}$.
\end{lem}
\begin{proof}This follows from Lemma \ref{basisechi} (ii) and Lemma \ref{relation} (ii).
\end{proof}

\begin{lem} There exists a short exact sequence 
\begin{displaymath}
\xymatrix{0\ar[r]&M^{\lambda}_{\gamma}\ar[r] & L^{\lambda}_{\gamma}\ar[r]& M^{\lambda}_{\gamma}\ar[r]& 0}
\end{displaymath}
of $\HH$-modules.
\end{lem}
\begin{proof}Let $v_1$ be the image of $e_{\chi}T_{n_s}$ in $L^{\lambda}_{\gamma}$ and let $v_2$ be be image of 
$e_{\chi}T_{n_s}T_{\Pi}$ in $L^{\lambda}_{\gamma}$. The 
subspace $\langle v_1, v_2\rangle_{\Fbar}$ is stable under the action of $T_{n_s}$, $T_{\Pi}$ and $e_{\chi'}$, for every 
$\chi':H\rightarrow \Fbar^{\times}$. Hence, by Lemma \ref{generelat} (ii) the subspace is stable under the action of $\HH$. From Lemma \ref{relation} (ii) 
and Definition \ref{list} (ii) it follows that $\langle v_1, v_2\rangle_{\Fbar}\cong M^{\lambda}_{\gamma}$. An easy check shows 
that $L^{\lambda}_{\gamma}/M^{\lambda}_{\gamma}\cong M^{\lambda}_{\gamma}.$  
\end{proof}

\begin{lem}\label{twistbyun} Let $(\pi,\VV)$ be a smooth representation of $G$ and let $\xi\in \Fbar^{\times}$.  Let $\mu_{\xi}$ be an unramified 
quasi-character:
$$\mu_{\xi}: F^{\times}\rightarrow \Fbar^{\times}, \quad x\mapsto \xi^{\val_F(x)} $$ 
where $\val_F$ is the valuation of $F$.  Suppose that $\pi^{I_1}$ contains 
$M_{\gamma}^{\lambda}$, where $\gamma=\{\chi, \chi^{s}\}$ and let 
$V$ be the underlying vector space of $M_{\gamma}^{\lambda}$ in $\VV$. If we consider the representation $(\pi\otimes\mu_{\xi}\circ\det$, $\VV$) of $G$, 
then the action of $\HH$ on $V$ is isomorphic to $M_{\gamma}^{\lambda\xi^{-2}}$.
\end{lem}
\begin{proof}Let 
$$V=\langle v_1, v_2\rangle_{\Fbar}$$
as in Definition \ref{list}. Since $\mu_{\xi}$ is unramified, Corollary \ref{T} implies that the action of $T_{n_s}$ and the idempotents $e_{\chi}$ on 
$V$ does not change. Lemma \ref{generelat} (ii) implies that it is enough to check how $T_{\Pi}$ acts. Since $\det \Pi= -\pif$, twisting by 
$\mu_{\xi}\circ\det$ gives us 
$$v_1 T_{\Pi}= \Pi^{-1}v_1= \xi^{-1}  v_2 \quad \textrm{and} \quad v_2 T_{\Pi} = \Pi^{-1}v_2= \xi^{-1} \lambda v_1.$$
Once  we replace $v_1$ by $\xi v_1$ the isomorphism follows from Definition \ref{list}.  
\end{proof}  

Since, by twisting by an unramified character we may vary $\lambda$ as we wish, we might as well work with $\lambda=1$.
\begin{defi}Let $\gamma=\{\chi,\chi^{s}\}$ then we define  $\HH$-modules
$$M_{\gamma}= M_{\gamma}^{1}\quad \textrm{and}\quad L_{\gamma}=L_{\gamma}^{1}.$$
\end{defi} 

\subsection{Restriction to $\HH_{K}$}
Let $\HH_K= \End_K(\Indu{I_1}{K}{\Eins})$. The natural isomorphism of $K$ representations
$$\Indu{I_1}{K}{\Eins}\cong \{f\in \cIndu{I_1}{G}{\Eins}: \supp f \subseteq K\}$$
gives  an  embedding of algebras
$$\HH_K\hookrightarrow \Hom_K(\Indu{I_1}{K}{\Eins},\cIndu{I_1}{G}{\Eins})\cong
\Hom_G(\cIndu{I_1}{G}{\Eins},\cIndu{I_1}{G}{\Eins})=\HH.$$
As an algebra $\HH_K$ is generated by $T_{n_s}$ and $e_{\chi}$, for all characters $\chi$.
\begin{defi}\label{mchiJK} Let $\chi:H\rightarrow \Fbar^{\times}$ be a character. Let $J_0(\chi)$ be a set, such that $J_0(\chi)=\emptyset$ 
if $\chi\neq \chi^{s}$, and
$J_0(\chi)=\{s\}$, if $\chi=\chi^{s}$. Let $J$ be a subset of $J_0(\chi)$, we define $M_{\chi,J}$ to be a right $\HH_K$-module, whose underlying vector 
space is one dimensional, $M_{\chi,J}=\langle v\rangle_{\Fbar}$ and the action of $\HH_{K}$ is determined by the following:
$$v e_{\chi}= v,$$
$$ v T_{n_s}=0\quad \textrm{if $s\in J$ or $s\not\in J_0(\chi)$},\quad v T_{n_s}=-v,\quad 
\textrm{if $s\not\in J$ and $s\in J_0(\chi)$.}$$
\end{defi}
Given $\chi$ and $J$ as above, we will denote
$$\overline{J}= J_0(\chi)\backslash J.$$ 
\begin{lem}\label{restML} Let $\chi:H\rightarrow \Fbar^{\times}$ be a character and let $\gamma=\{\chi,\chi^{s}\}$, then 
$$M_{\gamma}|_{\HH_K}\cong M_{\chi,J}\oplus M_{\chi^{s},\overline{J}}$$
as $\HH_K$-modules, where $J$ is a subset of $J_0(\chi)$. Moreover, if $\chi\neq \chi^{s}$, then
$$L_{\gamma}|_{\HH_K}\cong (\Indu{I}{K}{\chi}\oplus \Indu{I}{K}{\chi^{s}})^{I_1}$$
as $\HH_K$-modules. 
\end{lem} 
\begin{proof}The first isomorphism follows directly from Definition \ref{list}. Since $J_0(\chi)$ has at most two subsets, it doesn't matter which 
subset we take. For the second isomorphism we observe that the space $(\Indu{I}{K}{\chi})^{I_1}$ is two dimensional, with the basis $\{\varphi_{\chi},
T_{n_s}\varphi_{\chi^{s}}\}$. Moreover, $I$ acts on the basis vectors by characters $\chi$ and $\chi^{s}$ respectively. Now
$$\varphi_{\chi}T_{n_s}= \sum_{u\in I_1/K_1} u n_{s}^{-1}\varphi_{\chi}= e_{\chi} (\sum_{u\in I_1/K_1} u n_{s}^{-1}\varphi)= e_{\chi}T_{n_s}\varphi=
T_{n_s}e_{\chi^{s}}\varphi= T_{n_s} \varphi_{\chi^{s}}$$
and
$$(T_{n_s}\varphi_{\chi^s})T_{n_s}= \sum_{u\in I_1/K_1} u n_{s}^{-1}T_{n_s}e_{\chi^s}\varphi = 
T_{n_s}e_{\chi^{s}} (\sum_{u\in I_1/K_1} u n_{s}^{-1}\varphi)= e_{\chi}T_{n_s}^2\varphi=0$$
and Lemma \ref{lgammabasis} allows us to define the obvious isomorphism on the basis. 
\end{proof}
  
\section{ Irreducible  representations of $\GL_{2}(\Fq)$}\label{GLFq}
\subsection{Carter and Lusztig theory}
In \cite{cl} Carter and Lusztig have constructed all irreducible $\Fbar$-representations of a finite group $\Gamma$, which has a 'split $BN$-pair of 
characteristic $p$'. Since $\GL_2(\Fq)$ is a special case of this, we will  recall their results. Let $\Gamma$ be a finite group with a $BN$-pair 
$(\Gamma, B, N, S)$. Let $H= B\cap N$, then $H$ is normal in $N$, and 
$S$ is the set of Coxeter generators of $W=N/H$. We additionally require that $B=HU$, where $U$ is a normal subgroup of $B$, which is a $p$-group, 
and $H$ is abelian of order prime to $p$. Moreover, we assume that $H=\cap_{n\in N} B^{n}$.   

\begin{thm}\label{calu}\cite{cl} Let $\rho$ be an irreducible representation of $\Gamma$ then 
\begin{itemize}
\item[(i)] the space of $U$ invariants $\rho^{U}$ is one dimensional;
\item[(ii)] suppose that the action of $B$ on $\rho^{U}$ is given by a character $\chi: H\rightarrow \Fbar^{\times}$, via
$B\rightarrow B/U\cong H$
and let 
$ J=\{s\in S:  s\centerdot \rho^U = \rho^{U}\}$
then the pair $(\chi, J)$ determines $\rho$ up to an isomorphism;
\item[(iii)]conversely, given a character $\chi:H\rightarrow \Fbar^{\times}$, let $J_0(\chi)=\{s\in S: \chi^s=\chi\}$ and let  $J$ be a subset of 
$J_0(\chi)$ then there exists an irreducible representation $\rho_{\chi,J}$ of $\Gamma$ with the pair $(\chi, J)$ as above.
\end{itemize}
\end{thm}
\begin{proof} This is \cite{cl} Corollary 7.5, written out in detail, see also \cite{richen} Theorem 3.9 and \cite{curtis} Theorem 4.3.
\end{proof}

Let $\HH_{\Gamma}=\End_{\Gamma}(\Indu{U}{\Gamma}{\Eins})$. We would like to rephrase Theorem \ref{calu} in terms of $\HH_{\Gamma}$-modules.  
For each $s\in S$ we may choose a representative $n_s\in N$. Moreover, according to \cite{cl} Lemma 2.2,  we can  choose $n_s$ in a nice way. 
The obvious equivalent of Definition \ref{defTn} gives an endomorphism $T_{n}\in \HH_{\Gamma}$ for each $n\in N$. Definition \ref{simide} 
for each character $\chi:H\rightarrow \Fbar^{\times}$ gives an idempotent $e_{\chi}\in \HH_{\Gamma}$.

\begin{defi}\label{mchiJ}Let $\chi:H\rightarrow \Fbar^{\times}$ be a character, and let $J$ be a subset of $J_0(\chi)$ we define $M_{\chi,J}$ to be 
a right $\HH_{\Gamma}$-module, whose underlying vector space is one dimensional, $M_{\chi,J}=\langle v\rangle_{\Fbar}$ and the action of 
$\HH_{\Gamma}$ is determined by the following:
$$ve_{\chi}=v$$
and for every $s\in S$ we have 
\begin{displaymath}
v T_{n_s}=\left\{ \begin{array}{ll}
0 & \textrm{if $s\in J$,}\\
-v & \textrm{if $s\in J_0(\chi)$, $s\not\in J$,}\\
0 & \textrm{if $s\not\in J_0(\chi)$.}
\end{array}\right.
\end{displaymath}
\end{defi}

\begin{cor}\label{hgammamod} The functor of $U$ invariants
$$\Rep_{\Gamma}\rightarrow \Mod-\HH_{\Gamma},\quad \rho\mapsto \rho^{U}$$
induces a bijection between the irreducible representations of $\Gamma$ and the irreducible right $\HH_{\Gamma}$-modules. Moreover, if an irreducible 
representation $\rho_{\chi,J}$ corresponds to the pair $(\chi,J)$, in the sense of Theorem \ref{calu} (iii), then
$$\rho_{\chi,J}^{U}\cong M_{\chi,J}$$
as an $\HH_{\Gamma}$-module. 
\end{cor}
\begin{proof} Suppose that $\rho_{\chi,J}$ is an irreducible  representation as above. Then
 $$\rho_{\chi,J}^{U}=\langle v \rangle_{\Fbar}$$
and $H$ acts on $v$, by a character $\chi$.  Moreover, Lemma \ref{basisTinv} and Corollary 
\ref{T} hold (with obvious modifications). Hence, if $\chi:H\rightarrow \Fbar^{\times}$ is a character then
$$ve_{\chi}= \frac{1}{|H|}\sum_{h\in H} \chi(h)h^{-1}v= v.$$
If $\chi\neq \chi'$ then $e_{\chi}e_{\chi'}=0$ and hence $v e_{\chi'}=0$.  The set 
$\{T_{n_s}, e_{\chi}: s\in S, \chi:H\rightarrow \Fbar^{\times}\}$ generates $\HH_{\Gamma}$ as an algebra, so to determine the action, 
it is enough to compute 
$$v T_{n_s}= \sum_{u\in U/U\cap U^{s}}u n_{s}^{-1}v $$
for every $s\in S$. Now the right hand side is given by \cite{cl} Proposition 6.6, which implies that $\rho_{\chi,J}^U\cong M_{\chi,J}$. 

Since, by Theorem \ref{calu} (iii) every irreducible representation of $\Gamma$ corresponds to a pair $(\chi,J)$, it is enough to show that 
every irreducible module $M$ of $\HH_{\Gamma}$ is isomorphic to $M_{\chi,J}$, for some pair $(\chi,J)$. We adopt an argument of Vign\'eras, 
\cite{vig} E.7.1. We consider a representation $\rho(M)=M\otimes_{\HH_{\Gamma}}\Indu{U}{\Gamma}{\Eins}$, where $\Gamma$ acts on the right component of 
the tensor product. We have an injection of right $\HH_{\Gamma}$-modules 
$$M\hookrightarrow \rho(M)^{U},\quad m\mapsto m\otimes \varphi$$
where $\varphi$ takes value $1$ on $U$ and vanishes otherwise. Since, $\varphi$ generates $\Indu{U}{\Gamma}{\Eins}$ as a $\Gamma$-representation, the
image of $M$ in $\rho(M)$ will generate $\rho(M)$ as a $\Gamma$-representation. Hence, if $\rho$ is any non-zero irreducible quotient of $\rho(M)$, 
then since $M$ is irreducible, it 
will be a non-zero submodule of $\rho^{U}$, but by Theorem \ref{calu} (iii), $\rho$ will correspond to some pair $(\chi,J)$ and by above
$M\cong M_{\chi,J}.$
\end{proof}

\begin{remar}Ideally, we would like to have an analogue of the Corollary above for $G$ or more generally for any group of $F$-points of 
a reductive group, split over $F$.
\end{remar}      

Carter and Lusztig, in \cite{cl} construct all the irreducible representations $\rho_{\chi,J}$ in a very elegant way. For each pair $(\chi,J)$ 
they define a $\Gamma$-equivariant homomorphism
$$\Theta^{J}_{w_{0}}: \Indu{B}{\Gamma}{\chi}\rightarrow \Indu{B}{\Gamma}{\chi^{w_{0}}}$$
which depends on the geometry of the Coxeter group $W$, so that 
$$\rho_{\chi,J}\cong \Image \Theta^{J}_{w_{0}}$$
where $w_0$ is the unique element of maximal length in $W$.
 
From now onwards we specialise to our situation, so that $\Gamma=\GL_2(\Fq)$, $B$ is the subgroup of upper-triangular matrices, $U$ is the subgroup 
of unipotent upper-triangular matrices, $H$ is the diagonal matrices, $N$ is the normaliser of $H$ in $\Gamma$, that is the monomial matrices and 
$W= N/H$ is isomorphic to the symmetric group on two letters, $W=\{1, s\}$. Let 
$$n_{s}=\begin{pmatrix} 0 & -1\\ 1 & 0\end{pmatrix}$$
be a fixed representative of $s$ in $N$. In particular, $s$ is the element of the maximal length in $W$ and also the single
Coxeter generator, so that $S=\{s\}$. Hence, if $\chi:H\rightarrow \Fbar^{\times}$, then either $J_0(\chi)=\emptyset$ or $J_0(\chi)=S$.
Since 
$$ K/K_1\cong \Gamma,\quad I/K_1\cong B,\quad I_1/K_1\cong U$$
to ease the notation, we will often identify the spaces 
$$\{f:\Gamma\rightarrow \Fbar: f(ug)=f(g), \quad\forall g\in \Gamma, \quad\forall u \in U\}$$
and 
$$\{f\in \cIndu{I_1}{G}{\Eins}:\supp f\subseteq K\}$$
in the natural way. In particular, we will use the same notation for the elements of $\HH_K$ and $\HH_{\Gamma}$ and we note that the Definitions 
\ref{mchiJK} and \ref{mchiJ} coincide. 

\begin{prop}\label{gltwospecial} For each character $\chi:H\rightarrow \Fbar^{\times}$, such that $\chi=\chi^{s}$,  let 
$$\rho_{\chi,S}=\Image ((1+T_{n_s}): \Indu{B}{\Gamma}{\chi}\rightarrow \Indu{B}{\Gamma}{\chi})$$ 
and let 
$$ \rho_{\chi,\emptyset}=\Image (T_{n_s}: \Indu{B}{\Gamma}{\chi}\rightarrow \Indu{B}{\Gamma}{\chi})$$
then the representations $\rho_{\chi,S}$ and $\rho_{\chi,\emptyset}$ are irreducible. Moreover,
$$\rho_{\chi,S}^{U}=\langle (1+T_{n_s})\varphi_{\chi}\rangle_{\Fbar}\cong M_{\chi,S}\quad \textrm{and}\quad 
\rho_{\chi,\emptyset}^{U}=\langle T_{n_s}\varphi_{\chi}\rangle_{\Fbar}\cong M_{\chi,\emptyset}$$
as $\HH_{\Gamma}$-modules. For each character $\chi:H\rightarrow \Fbar^{\times}$, such that $\chi\neq\chi^{s}$, let 
$$\rho_{\chi,\emptyset}=\Image(T_{n_s}: \Indu{B}{\Gamma}{\chi}\rightarrow \Indu{B}{\Gamma}{\chi^{s}})$$
then the representation $\rho_{\chi,\emptyset}$ is irreducible. Moreover,
$$\rho_{\chi,\emptyset}^{U}=\langle T_{n_s}\varphi_{\chi}\rangle_{\Fbar}\cong M_{\chi,\emptyset}$$
as an $\HH_{\Gamma}$-module. Further, these representations are pairwise non-iso\-mor\-phic, and every irreducible representation of $\Gamma$ is 
isomorphic to $\rho_{\chi,J}$, for some character $\chi$ and a subset $J$ of $J_0(\chi)$.
\end{prop}
\begin{proof} This is a special case of \cite{cl} Theorem 7.1 and Corollary 7.5. The isomorphisms of $\HH_{\Gamma}$-modules are given by the Corollary 
\ref{hgammamod}.
\end{proof}  

\begin{remar} Although we do not use this, we note that Frobenius reciprocity gives us
$$\cIndu{K}{G}{\rho_{\chi,\emptyset}}\cong T_{n_s}(\cIndu{I}{G}{\chi})\le \cIndu{I}{G}{\chi^{s}}$$
and if $\chi=\chi^{s}$ then
$$\cIndu{K}{G}{\rho_{\chi,S}}\cong (1+T_{n_s})(\cIndu{I}{G}{\chi})\le \cIndu{I}{G}{\chi}.$$
Using this, one can relate the central elements of Vign\'eras in \cite{vig} to the `standard ' endomorphisms $T_{\sigma}$ of 
Barthel and Livn\'e in \cite{bl2}.
\end{remar}
\begin{lem}\label{splitiw} Let $\chi:H\rightarrow \Fbar^{\times}$ be a character, such that $\chi=\chi^{s}$. 
Then the homomorphisms $e_{\chi}(1+T_{n_s})e_{\chi}$
and $-e_{\chi}T_{n_s} e_{\chi}$ are orthogonal idempotents. In particular,
$$\Indu{B}{\Gamma}{\chi}\cong \rho_{\chi, \emptyset}\oplus \rho_{\chi, S}.$$
Moreover, let $\chi':\Fq^{\times}\rightarrow \Fbar^{\times}$ be a character such that $\chi=\chi'\circ\det$, then 
$$\rho_{\chi,S}\cong \chi'\circ \det\quad \and \rho_{\chi,\emptyset}\cong St\otimes \chi'\circ\det$$
where $St$ is the Steinberg representation.
\end{lem}
\begin{proof}Since $\chi=\chi^{s}$ we have
$$ e_{\chi}T_{n_s}= T_{n_s} e_{\chi}\quad\textrm{and}\quad e_{\chi} T_{n_s}^{2}=- e_{\chi}T_{n_s}.$$
So the elements above are orthogonal idempotents as claimed. By Proposition \ref{gltwospecial}, the summands they split off are $\rho_{\chi,S}$
and $\rho_{\chi,\emptyset}$.

Since $\chi=\chi^{s}$, the character $\chi$ must factor through the determinant. So $\chi$ extends to a character of $\Gamma$ and hence
$$\Indu{B}{\Gamma}{\chi}\cong \Indu{B}{\Gamma}{\Eins}\otimes \chi'\circ\det.$$
So we may assume that $\chi$ is the trivial character. The Bruhat decomposition says that $\Gamma=BsB \cup B$ and hence by Theorem \ref{calu} (ii) 
$\rho_{\Eins, S}= \Eins$, the trivial representations of $G$. This implies that $\rho_{\Eins, \emptyset}$ is the Steinberg representation. 
\end{proof}

\begin{cor}\label{splitisgood}Let $\chi:H\rightarrow \Fbar^{\times}$ be a character, such that $\chi=\chi^{s}$. Let $\rho$ be any representation
 of $\Gamma$, such that  for some  $v \in \rho^U$ we have
$$\langle v\rangle_{\Fbar}\cong M_{\chi, J}$$
 as an $\HH_{\Gamma}$-module. Then
$$\langle \Gamma v \rangle_{\Fbar}\cong \rho_{\chi, J}$$
as a $\Gamma$-representation.
\end{cor}
\begin{proof} Since $v$ is fixed by $U$   there exists a homomorphism 
$$\psi\in \Hom_{\Gamma}(\Indu{U}{\Gamma}{\Eins}, \rho)$$
such that $ \psi(\varphi)=v.$
The isomorphism of $\HH_{\Gamma}$-modules implies that
$$v=v e_{\chi}=\psi(e_{\chi}\varphi)=\psi(\varphi_{\chi}). $$
Hence, $H$ acts on $v$ by a character $\chi$ and 
$$\psi(\Indu{U}{\Gamma}{\Eins})=\psi(e_{\chi}(\Indu{U}{\Gamma}{\Eins}))=\psi(\Indu{B}{\Gamma}{\chi}).$$
If $J=\emptyset$ then 
$$\psi((1+T_{n_s})\varphi_{\chi})=v(1+T_{n_s})e_{\chi}=0.$$
Hence,  $\rho_{\chi, S}$ is contained in the kernel of $\psi$. By Lemma \ref{splitiw} 
$$\Image \psi\cong \rho_{\chi, \emptyset}.$$
Since, the image is irreducible and contains $v$ we get the result. The proof for $J=S$ is analogous.
\end{proof} 
The Corollary has  a nice application, which complements \cite{vig} E.7.1.
\begin{cor}\label{trivialrep} Let $\pi$ be a smooth representation of $G$ and suppose that there exists a non-zero vector $v\in \pi^{I_1}$ such that 
$$v e_{\Eins}= v,\quad vT_{n_s}=0,\quad vT_{\Pi}=v$$
then $G$ acts trivially on $v$.
\end{cor} 
\begin{proof} As an $\HH_K$ module
$$\langle v \rangle_{\Fbar}\cong M_{\Eins, S}.$$
By Corollary \ref{splitisgood} $K$ acts trivially on $v$. On the other hand 
$$ v= vT_{\Pi}=\Pi^{-1}v.$$
Iwahori decomposition implies that $\Pi$ and $K$ generate $G$ as a group. Hence $G$ acts trivially on $v$.
\end{proof}
\begin{remar} There is a version of this twisted by a character. This example will lead us to better things. See Remark \ref{explicitaction}.
\end{remar} 
\begin{lem}\label{nonsplit} Let $\chi: H\rightarrow \Fbar^{\times}$ be a character, let $J$ be a subset of $J_0(\chi)$, and let 
$\overline{J}=J_0(\chi)\backslash J$.  The  sequence of $\HH_{\Gamma}$-modules 
\begin{displaymath}
\xymatrix{0\ar[r]&M_{\chi, J}\ar[r]& (\Ind_{B}^{\Gamma}{\chi^{s}})^{U}\ar[r] &M_{\chi^{s},\overline{J}}\ar[r] &0}
\end{displaymath}
is exact. Moreover, it splits if and only if $\chi=\chi^{s}$.
\end{lem}
\begin{proof} The space $(\Ind_{B}^{\Gamma}{\chi^{s}})^{U}$ is two dimensional, with the basis $\{T_{n_s}\varphi_{\chi}, \varphi_{\chi^{s}}\}$.

If $\chi=\chi^{s}$ then $e_{\chi}(1+T_{n_s})e_{\chi}$ and $-e_{\chi} T_{n_s} e_{\chi}$ are orthogonal idempotents, which split the sequence.

If $\chi\neq \chi^{s}$ then for any $\lambda$, $\mu\in \Fbar$ we have
$$(\lambda T_{n_s}\varphi_{\chi}+ \mu \varphi_{\chi^{s}})e_{\chi}= \lambda T_{n_s}\varphi_{\chi}, \quad
(\lambda T_{n_s}\varphi_{\chi}+ \mu \varphi_{\chi^{s}})e_{\chi^{s}}= \mu \varphi_{\chi^{s}}$$
and
$$(\lambda T_{n_s}\varphi_{\chi}+ \mu \varphi_{\chi^{s}})T_{n_s}= \mu T_{n_s}\varphi_{\chi}.$$
Hence $M_{\chi,\emptyset}$ is the only proper submodule, so the sequence cannot split. 
\end{proof}

\subsection{Alternative description of irreducible representations}\label{alternative}
Let $\vs{d}{F}$ be an $F$ vector space of homogeneous polynomials in two variables $X$ and $Y$
of the degree $d$. The group $K$ acts on $\vs{d}{F}$ via
$$\begin{pmatrix} a & b \\ c & d\end{pmatrix}(X^{d-i}Y^{i})=(aX+cY)^{d-i}(bX+ dY)^{i}.$$  
For $0\le i \le d$, let
$$m_i= \begin{pmatrix} d \\ i \end{pmatrix} X^{d-i}Y^{i}$$
where $\begin{pmatrix} d \\ i \end{pmatrix}$ denotes the binomial coefficient. Vectors $m_i$, for $0\le i\le d$, form a basis
of $\vs{d}{F}$. Let $\vs{d}{\oF}$ be the $\oF$-lattice in $\vs{d}{F}$ spanned by the $m_i$, for $0\le i\le d$.
 An easy check shows that $\vs{d}{\oF}$ is $K$ invariant. Let 
$$\vs{d}{\Fq}=\vs{d}{\oF}\otimes_{\oF} \oF/\pF.$$ 
 The vectors $m_i\otimes 1$, for $0\le i\le d$, form an $\Fq$-basis
of $\vs{d}{\Fq}$. The subgroup $K_1$ acts trivially on  $\vs{d}{\Fq}$, so we consider $\vs{d}{\Fq}$ as a representation of $\Gamma$.
Let $\Fr$ be the automorphism of $\Gamma$, given by 
$$\Fr: \begin{pmatrix} a & b\\ c & d\end{pmatrix}\mapsto  \begin{pmatrix} a^{p} & b^{p}\\ c^{p} & d^{p}\end{pmatrix}.$$ 
Let $\rho$ be  a representation of $\Gamma$. We will denote by $\rho^{\Fr}$ the representation of $\Gamma$ given by
$$\rho^{\Fr}(g)=\rho(\Fr(g)).$$

\begin{thm}\label{irrgl} Let $\Gamma=\GL_2(\Fq)$ and suppose that $q=p^{n}$. The isomorphism classes of 
irreducible $\Fbar$-representations of $\Gamma$ are parameterised by pairs $(a,\rr)$, where 
\begin{itemize}
\item[-] $a$ is an integer $1\le a \le q-1$ and
\item[-] $\rr$ is an ordered $n$-tuple $\rr=(r_0, r_1, \ldots, r_{n-1})$, where $0\le r_i \le p-1$, for every $i$.
\end{itemize} 
Moreover, the irreducible representations of $\Gamma$ can be realized over $\Fq$ and the irreducible representation corresponding to 
$(a,\rr)$ is given by 
$$V_{\rr,\Fq}\otimes (\det)^{a}\cong \vs{r_0}{\Fq}\otimes \vs{r_1}{\Fq}^{\Fr}\otimes\ldots\otimes\vs{r_i}{\Fq}^{\Fr^{i}}\otimes\ldots\otimes
\vs{r_{n-1}}{\Fq}^{\Fr^{n-1}}\otimes (\det)^{a}.$$
\end{thm}  
\begin{proof} This is shown in \cite{brauer}, see also \cite{bl2} Proposition 1 and \cite{vig} Ap. 6. We remark that since 
$\begin{pmatrix} r\\ i \end{pmatrix}$
is a unit in $\Fq$ if $r\le p-1$, our spaces really coincide with the ones considered in \cite{bl2}.
\end{proof}

We fix some embedding $\iota: \Fq\hookrightarrow \Fbar$ and we will assume that every character $\chi:H\rightarrow \Fbar^{\times}$ factors through
$\iota$.  Once we have done that, we will omit $\iota$ from our notation. We will denote
$$V_{\rr,\Fbar}=V_{\rr,\Fq}\otimes_{\Fq}\Fbar.$$
We need a dictionary between the two descriptions. 

\begin{prop}\label{dictionary} Let $\chi: H\rightarrow \Fbar^{\times}$ and let $a$ be the unique integer, such that $1\le a\le q-1$ and
$$ \chi(\begin{pmatrix} 1 & 0 \\0 & \lambda \end{pmatrix})= \lambda^{a}\quad \forall \lambda\in \Fq^{\times}$$
and let $r$ be the unique integer, such that $1\le r \le q-1$ and
$$ \chi(\begin{pmatrix} \lambda & 0 \\0 & \lambda^{-1} \end{pmatrix})= \lambda^{r}\quad \forall \lambda\in \Fq^{\times}.$$
Suppose that  $r\neq q-1$, and let $\rr=(r_0,\ldots,r_{n-1})$ be the digits of a $p$-adic expansion of $r$
$$r=r_0+ r_1 p+\ldots +r_{n-1} p^{n-1}$$
then $\chi\neq\chi^{s}$ and $\rho_{\chi, \emptyset}$ corresponds to the pair $(a, \rr)$. More precisely
$$\rho_{\chi, \emptyset}\cong \vs{r_0}{\Fbar}\otimes\ldots\otimes\vs{r_{n-1}}{\Fbar}^{\Fr^{n-1}}\otimes (\det)^{a}.$$
Suppose that $r=q-1$, then $\chi=\chi^{s}$, 
$$\rho_{\chi, \emptyset}\cong \vs{p-1}{\Fbar}\otimes\ldots\otimes\vs{p-1}{\Fbar}^{\Fr^{n-1}}\otimes (\det)^{a}\cong \St\otimes(\det)^{a}$$
and
$$\rho_{\chi, S}\cong \vs{0}{\Fbar}\otimes\ldots\otimes\vs{0}{\Fbar}^{\Fr^{n-1}}\otimes (\det)^{a}\cong (\det)^{a}$$
where $\St$ denotes the Steinberg representation.
\end{prop}  
\begin{proof} Every character $\chi:H\rightarrow \Fbar^{\times}$ is of the form
$$\chi: \begin{pmatrix} \lambda & 0 \\ 0 & \mu \end{pmatrix}\mapsto \lambda^{c}\mu^{d}$$
for some integers $c$ and $d$. Moreover, $\chi=\chi^s$ if and only if 
$$c-d\equiv 0 \pmod{q-1}.$$
The integers $a$ and $r$ are uniquely determined by the congruences
$$d\equiv a\pmod{q-1}\quad\mbox{and}\quad c-d\equiv r \pmod{q-1}.$$  
By Theorem \ref{calu} if $\rho$ is an irreducible representation of $\Gamma$, then $\dim \rho^{U}=1$, and by Corollary \ref{hgammamod} 
the irreducible representations of $\Gamma$ correspond to the irreducible modules of the Hecke algebra $\HH_{\Gamma}$. Since we have two complete 
lists of irreducible representations, it is enough to match  up the corresponding irreducible modules. We recall that
$$\rho_{\chi,J}^U\cong M_{\chi, J}$$
as $\HH_{\Gamma}$-modules. 

We observe that the action of $U$ on $\vs{d}{\Fbar}$ fixes the vector $m_0\otimes 1$. Moreover,
$$\begin{pmatrix} \lambda & 0 \\ 0 & \mu \end{pmatrix}m_0 \otimes 1= \lambda^{d}m_0 \otimes 1.$$
Let $(a,\rr)$ be any pair parameterising an irreducible representation of $\Gamma$ and let
$$ r= r_0+ r_1 p+\ldots r_{n-1}p^{n-1}.$$
 By picking such $(m_0\otimes 1)_{r_i}$ in every component of the tensor product we obtain a non-zero vector
$$(m_0\otimes 1)_{\rr}=(m_0\otimes 1)_{r_0}\otimes\ldots \otimes (m_0\otimes 1)_{r_{n-1}}$$
fixed by $U$. The vector $(m_0\otimes 1)_{\rr}$ spans  the space of $U$ invariants, since it  is one dimensional. Moreover,
since the action on the components of the tensor product is twisted by $\Fr$ we obtain
 $$\begin{pmatrix} \lambda & 0 \\ 0 & \mu \end{pmatrix}(m_0 \otimes 1)_{\rr}= (\lambda \mu)^{a}\lambda^{r}(m_0 \otimes 1)_{\rr}.$$

Suppose that we start with an arbitrary character $\chi: H\rightarrow \Fbar^{\times}$ and obtain the integers $a$ and $r$ as in the 
statement of the proposition. 

If $r\neq q-1$, then by above $\chi\neq \chi^{s}$. Let $\rr$ be the $n$-tuple corresponding to $r$. Since, 
$\chi\neq \chi^{s}$, the module  $M_{\chi,\emptyset}$ is the only  irreducible module of $\HH_{\Gamma}$, which is not killed by the
 idempotent $e_{\chi}$. Let $(m_0\otimes 1)_{\rr}$ be the vector constructed above. Since, $H$  acts on $(m_0\otimes 1)_{\rr}$ via 
the character $\chi$, we obtain
$$ M_{\chi, \emptyset}\cong   (\vs{r_0}{\Fbar}\otimes\ldots\otimes\vs{r_{n-1}}{\Fbar}^{\Fr^{n-1}}\otimes (\det)^{a})^{U}$$
as $\HH_{\Gamma}$-modules and that implies the isomorphism between representations.

If $r=q-1$, then $\chi=\chi^{s}$, and the only $\HH_{\Gamma}$-modules, which are  not killed by $e_{\chi}$, are $M_{\chi, S}$ and $M_{\chi,\emptyset}$. 
We observe that $\vs{0}{\Fbar}$ is just the trivial representation. Let $\mathbf{0}=(0,\ldots,0)$, then the representation corresponding to the pair 
$(a,\mathbf{0})$ is just $\Eins\otimes(\det)^a$, which is isomorphic to $\rho_{\chi,S}$, by Proposition \ref{splitiw}. 
The only case left is  $\rr=\mathbf{p-1}=(p-1,\ldots, p-1)$, hence 
$$ M_{\chi, \emptyset}\cong   (\vs{p-1}{\Fbar}\otimes\ldots\otimes\vs{p-1}{\Fbar}^{\Fr^{n-1}}\otimes (\det)^{a})^{U}$$
as $\HH_{\Gamma}$-modules, since the module $M_{\chi, S}$ is already taken. This implies that 
$$\rho_{\chi,\emptyset}\cong V_{\mathbf{p-1}, \Fbar}\otimes(\det)^{a}\cong St\otimes (\det)^{a}$$
where the last isomorphism follows from Proposition \ref{splitiw}.
\end{proof}

\begin{cor}\label{bar} Suppose that $q=p^{n}$ and the representation $\rho_{\chi, J}$ corresponds to the pair $(a,\rr)$. 
Let $r=r_0+r_1p+\ldots+r_{n-1}p^{n-1}$ and let  $\overline{J}= J_0(\chi)\backslash J$, where $J_0(\chi)=\{s\in S: \chi^{s}= \chi\}$. 
Then
$$\rho_{\chi^{s}, \overline{J}}\cong\vs{p-1-r_0}{\Fbar}\otimes\ldots\otimes\vs{p-1-r_{n-1}}{\Fbar}^{\Fr^{n-1}}\otimes (\det)^{a+r}.$$
\end{cor}
\begin{proof}If $r=0$ or $r=q-1$, then $\rr$ is of a special form and the isomorphism follows from Proposition \ref{dictionary}.

If $r\neq 0$ and $r\neq q-1$, we observe that
$$ \chi^s(\begin{pmatrix} 1 & 0 \\0 & \lambda \end{pmatrix})=
\chi(\begin{pmatrix} \lambda & 0 \\0 & \lambda^{-1} \end{pmatrix})\chi(\begin{pmatrix} 1 & 0 \\0 & \lambda \end{pmatrix})= \lambda^{a+r}\quad 
\forall \lambda\in \Fq^{\times}$$
 and
$$ \chi^s(\begin{pmatrix} \lambda & 0 \\0 & \lambda^{-1} \end{pmatrix})=\chi(\begin{pmatrix} \lambda^{-1} & 0 \\0 & \lambda \end{pmatrix})=
 \lambda^{-r}\quad \forall \lambda\in \Fq^{\times}.$$
The claim follows from Proposition \ref{dictionary}.
\end{proof}

\section{Principal indecomposable representations}\label{princindrep}

We will recall some facts from the modular representation theory of finite groups. Let $\Gamma$ be any finite group. We denote by 
$\Rep_{\Gamma}$ the category of $\Fbar$ representations of $\Gamma$ and by 
$\Irr_{\Gamma}$ the set of isomorphism classes of irreducible representations in $\Rep_{\Gamma}$. We note that $\Rep_{\Gamma}$ is equivalent to 
the module category of the ring $\Fbar[\Gamma]$.  

\begin{prop}\label{injrho}  A representation $\inj$ is an injective object in $\Rep_{\Gamma}$ if and only if it is a projective object in
$\Rep_{\Gamma}$.

The isomorphism classes of indecomposable injective (and hence projective)
 objects in $\Rep_{Gamma}$ are parameterised by $\Irr_{\Gamma}$.

More precisely, if $\inj$ is indecomposable and injective, then the maximal semi-simple submodule $\soc(\inj)$ and the maximal 
semi-simple quotient $\inj/ \rad(\inj)$ are both irreducible. Moreover,
$$\soc(\inj)\cong \inj/ \rad(\inj).$$
Conversely, given $\rho\in \Irr_{\Gamma}$, there exists a unique up to isomorphism indecomposable, injective object $\inj \rho$ in $\Rep_{\Gamma}$, such 
that 
$$\rho\cong\soc(\inj \rho).$$
\end{prop}
\begin{proof} See \cite{serre}, Exercises 14.1 and 14.6.
\end{proof}   
We will call indecomposable representations of $\Gamma$, which are injective objects in $\Rep_{\Gamma}$, principal indecomposable representations.
\begin{remar}\label{note} We note that a monomorphism  $\rho\hookrightarrow \inj \rho$ is an injective envelope of $\rho$ in $\Rep_{\Gamma}$.
\end{remar}
\begin{cor}\label{injsumm}We have the following decomposition:
$$\Fbar[\Gamma]\cong\bigoplus_{\rho\in \Irr_{\Gamma}}(\dim \rho)\inj \rho.$$ 
\end{cor}
\begin{proof}Since $\Fbar[\Gamma]$ is an injective and projective object it must decompose into a direct sum of indecomposable injective objects.
Since 
$$\dim \Hom_{\Gamma}(\rho, \Fbar[\Gamma])=\dim \Hom_{\{1\}}(\rho, \Eins)= \dim \rho$$
the representation $\inj \rho$ occurs in the decomposition with the multiplicity $\dim \rho$.
\end{proof}

\begin{prop}\label{restsylow}Let $U$ be a $p$-Sylow subgroup of $\Gamma$. Then a representation $\rho$ is an injective object in 
$\Rep_{\Gamma}$ if and only if $\rho|_U$ is an injective object in $\Rep_U$.
\end{prop}
\begin{proof}This follows easily from \cite{serre}, \S 14.4, Lemma 20.
\end{proof}
\begin{prop}\label{simplep} Suppose that $U$ is a $p$-group, then the only irreducible representation is $\Eins$ and 
hence the only principal indecomposable representation is $\Fbar[U]$.
\end{prop}
\begin{proof} The first part is \cite{serre}, \S 8 , Proposition 26, the last part follows from Corollary \ref{injsumm}.
\end{proof}

\begin{cor}\label{dim}Let $\inj$ be an injective object in $\Rep_{\Gamma}$ and let $U$ be a $p$-Sylow subgroup of $\Gamma$, then
$$\dim \inj= \dim \inj^{U} |U|.$$
\end{cor}
\begin{proof} The restriction $\inj\mid_{U}$ is an injective object in $\Rep_U$. By the above Proposition 
$$\inj\mid_{U}\cong m \Fbar[U].$$
The multiplicity $m$ is given by:
$m=\dim \Hom_{U}(\Eins, \inj)=\dim \inj^{U}.$
\end{proof}

In the rest of the section $\Gamma=\GL_2(\Fq)$ and $U$ is the subgroup of unipotent upper triangular matrices.
 Given $\rho\in\Irr_{\Gamma}$ we are going  to compute $(\inj \rho)^{U}$ as an $\HH_{\Gamma}$-module.  Once we know the modules we are going to show that
if we consider $\inj \rho_{\chi,J}$ and $\inj \rho_{\chi^{s},\overline{J}}$ as representations of $K$, then the action of $\HH_{K}$ on 
$$(\inj \rho_{\chi,J}\oplus \inj \rho_{\chi^{s},\overline{J}})^{I_1}$$
extends to the action of $\HH$, so that if $\chi=\chi^{s}$  then it is isomorphic to a direct sum of supersingular modules and if 
$\chi\neq \chi^{s}$ then
it is isomorphic to a direct sum of $L_{\gamma}$ and supersingular modules. See Propositions \ref{injiwahori} and \ref{injregular} for the precise 
statement.  This calculation,
becomes of importance in Section \ref{II}. Although, the general case includes the case $q=p$, if $q=p$ we give a different, easier way of doing this. 
When $q=p$,  the main result is Proposition \ref{injqp}.

\subsection{The case $q=p$}
We start off with no assumption on $q$. 

\begin{lem} \label{subsp}Suppose that $\chi\neq \chi^s$, then there exists an exact sequence
\begin{displaymath}
\xymatrix{0\ar[r]&\Ind_{B}^{\Gamma}\chi^{s}\ar[r]^{\psi}& \inj \rho_{\chi,\emptyset}}
\end{displaymath} 
of $\Gamma$-representations.
\end{lem}
\begin{proof}Since $\inj \rho_{\chi,\emptyset}$ is an injective module, there exists $\psi$ such that the diagram 
\begin{displaymath}
\xymatrix{ 0\ar[r]& \rho_{\chi,\emptyset}\ar[r]\ar[d]  & \Indu{B}{\Gamma}{\chi^{s}} \ar@{.>}[dl]^{\psi}\\
                  & \inj \rho_{\chi,\emptyset}         &   }
\end{displaymath}

commutes. If $\Ker \psi\neq 0$, then $(\Ker \psi)^{U}$ is a non-zero proper submodule of $(\Indu{B}{\Gamma}{\chi^{s}})^{U}$ not containing 
$M_{\chi, \emptyset}$. By Lemma \ref{nonsplit} this cannot happen.
\end{proof}

\begin{cor}\label{dimge}Suppose that $\chi\neq \chi^{s}$ then 
$$\dim \inj \rho_{\chi,\emptyset}\ge 2q.$$
\end{cor}
\begin{proof}Corollary \ref{dim} implies that
$$\dim \inj \rho_{\chi,\emptyset}=\dim (\inj \rho_{\chi,\emptyset})^{U} |U|.$$
The order of $U$ is $q$ and since by Lemma \ref{subsp} $\Indu{B}{\Gamma}{\chi^{s}}$ is a subspace of $\inj \rho_{\chi,\emptyset}$, we obtain
$$\dim (\inj \rho_{\chi,\emptyset})^{U}\ge 2.$$
\end{proof}

\begin{lem}\label{exseq}Suppose that $q=p$ and $\chi\neq \chi^{s}$ then the sequence of $\Gamma$ representations
\begin{displaymath}
\xymatrix{0\ar[r]&\rho_{\chi, \emptyset}\ar[r]& \Indu{B}{\Gamma}{\chi^{s}}\ar[r] &\rho_{\chi^{s},\emptyset}\ar[r] &0}
\end{displaymath}
is exact.
\end{lem}
\begin{remar}\label{fail} This fails if $q\neq p$.
\end{remar}
\begin{proof}The argument below is taken from \cite{vig} Ap. 6. We know that 
$$\rho_{\chi^{s},\emptyset}\cong T_{n_s}(\Indu{B}{\Gamma}{\chi^{s}})$$
and $\rho_{\chi,\emptyset}$ is isomorphic to the subspace of $\Indu{B}{\Gamma}{\chi^{s}}$ generated by $T_{n_s}\varphi_{\chi}$.
Since, $T_{n_s}^2 \varphi_{\chi}=0$ we always have
$$\rho_{\chi,\emptyset}\le \Ker T_{n_s}.$$
If $q=p$, then by Proposition \ref{dictionary} and Corollary \ref{bar} there exists an integer $r$ such that
$$\dim \rho_{\chi,\emptyset}+\dim \rho_{\chi^{s},\emptyset}= (r+1)+(p-1-r+1)=p+1=\dim \Indu{B}{\Gamma}{\chi^{s}}.$$ 
Hence the sequence is exact.
\end{proof}

\begin{cor}\label{qequalpisgood}Suppose that $q=p$ and let $\chi:H\rightarrow \Fbar^{\times}$ be a character, 
such that $\chi\neq\chi^{s}$. Let $\rho$ be any representation
 of $\Gamma$, such that  for some  $v \in \rho^U$
$$\langle v\rangle_{\Fbar}\cong M_{\chi, \emptyset}$$
 as an $\HH_{\Gamma}$-module. Then
$$\langle \Gamma v \rangle_{\Fbar}\cong \rho_{\chi, \emptyset}$$
as a $\Gamma$-representation.
\end{cor}
\begin{remar}\label{hooray} This fails if $p\neq q$, by Remark \ref{fail}, it is enough to look at $\Indu{B}{\Gamma}{\chi}/ \rho_{\chi^{s},\emptyset}$.
\end{remar}
\begin{proof} Since $v$ is fixed by $U$, there exists a homomorphism 
$$\psi\in \Hom_{\Gamma}(\Indu{U}{\Gamma}{\Eins}, \rho)$$
such that $\psi(\varphi)=v.$ The isomorphism of $\HH_{\Gamma}$-modules implies that
$$v=v e_{\chi}= \psi(e_{\chi}\varphi)=\psi(\varphi_{\chi}).$$
Hence $H$ acts on $v$ by a character $\chi$ and 
$$\psi(\Indu{U}{\Gamma}{\Eins})=\psi(e_{\chi}(\Indu{U}{\Gamma}{\Eins}))=\psi(\Indu{B}{\Gamma}{\chi}).$$
 Now
$$\psi(T_{n_s}\varphi_{\chi^{s}})=v T_{n_s} e_{\chi^{s}}=0.$$
Hence, $\rho_{\chi^{s}, \emptyset}$ is contained in the kernel of $ \psi$. By Lemma \ref{exseq} 
$$\Image \psi\cong \rho_{\chi, \emptyset}.$$
Since, the image is irreducible and contains $v$ we get the result. 
\end{proof}

\begin{lem}\label{qequalpdiminj}Suppose that $q=p$. If $\chi=\chi^{s}$ then
$$\dim \inj \rho_{\chi,J}=p.$$
If $\chi\neq\chi^{s}$ then
$$\dim \inj \rho_{\chi, \emptyset}= 2p.$$
\end{lem}
\begin{proof}
\begin{equation}\notag
\begin{split}
\dim \Fbar[\Gamma]=&\sum_{\rho\in \Irr_{\Gamma}}(\dim \rho)(\dim \inj \rho)\\
                  =&\sum_{\chi=\chi^{s}} (\dim \rho_{\chi,\emptyset})(\dim \inj \rho_{\chi,\emptyset})
+(\dim \rho_{\chi,S})(\dim \inj \rho_{\chi,S})\\
+\frac{1}{2}&\sum_{\chi\neq\chi^{s}}(\dim \rho_{\chi,\emptyset})(\dim \inj \rho_{\chi,\emptyset})+
(\dim \rho_{\chi^{s},\emptyset})(\dim \inj \rho_{\chi^{s},\emptyset}).
\end{split}
\end{equation}
If $\chi=\chi^{s}$ then Corollary \ref{dim} implies that
$$\dim \inj \rho_{\chi, J}\ge p.$$
If $\chi\neq\chi^{s}$ then Corollary \ref{dimge} implies that
$$\dim \inj \rho_{\chi, \emptyset}\ge 2p.$$
Lemma \ref{exseq} and Lemma \ref{splitiw}  imply that 
$$\dim \rho_{\chi,J}+\dim \rho_{\chi^{s}, \overline{J}}=p+1.$$
We put these inequalities together and we obtain
$$\dim \Fbar[\Gamma]\ge\sum_{\chi} (p+1)p= \dim \Fbar[\Gamma]$$
So all the inequalities must be equalities and we obtain the lemma.   
\end{proof}

\begin{cor}\label{qpinjg}Suppose that $q=p$. If $\chi=\chi^{s}$ then 
$$\langle \Gamma  (\inj \rho_{\chi, J})^U \rangle_{\Fbar}\cong \rho_{\chi, J}.$$
In particular,
$$(\inj \rho_{\chi, J})^U\cong M_{\chi,J}$$
as an $\HH_{\Gamma}$-module.

If $\chi\neq\chi^{s}$ then 
$$\langle \Gamma  (\inj \rho_{\chi, \emptyset})^U \rangle_{\Fbar}\cong \Indu{B}{\Gamma}{\chi^{s}}.$$
In particular,
$$(\inj \rho_{\chi, J})^U\cong (\Indu{B}{\Gamma}{\chi^{s}})^{U}$$
as an $\HH_{\Gamma}$-module.
\end{cor}
\begin{proof}If $\chi=\chi^{s}$ then we have an exact sequence 
\begin{displaymath}
\xymatrix{0\ar[r]&\rho_{\chi, J}\ar[r]& \inj \rho_{\chi, J}}
\end{displaymath}
of $\Gamma$-representations. Since, by Lemma \ref{qequalpdiminj} 
$$\dim \rho_{\chi, J}^{U}= \dim (\inj \rho_{\chi,J})^{U}$$
we obtain the Corollary. Similarly, if $\chi\neq \chi^{s}$ then by Lemma \ref{subsp} there exists an exact sequence
\begin{displaymath}
\xymatrix{0\ar[r]&\Ind_{B}^{\Gamma}\chi^{s}\ar[r] & \inj \rho_{\chi,\emptyset}}
\end{displaymath} 
of $\Gamma$-representations. Since, by Lemma \ref{qequalpdiminj}
$$\dim (\Indu{B}{\Gamma}{\chi^{s}})^{U}=\dim (\inj \rho_{\chi,\emptyset})^{U}$$
we obtain the Corollary.
\end{proof}

\begin{prop}\label{injqp} Suppose that $q=p$, let $\chi:H\rightarrow \Fbar^{\times}$ be a character and let $\gamma=\{\chi,\chi^{s}\}$. 
We consider representations $\inj \rho_{\chi,J}$ and  $\inj \rho_{\chi,\overline{J}}$ as representations of $K$, via 
$$K\rightarrow K/K_1\cong \Gamma.$$
If $\chi= \chi^s$ then the action of $\HH_{K}$ on $(\inj \rho_{\chi,\emptyset}\oplus \inj \rho_{\chi,S})^{I_1}$ extends to the action of $\HH$ so that
$$(\inj \rho_{\chi,\emptyset}\oplus \inj \rho_{\chi,S})^{I_1}\cong M_{\gamma}.$$
If $\chi\neq\chi^s$ then the action of $\HH_{K}$ on $(\inj \rho_{\chi,\emptyset}\oplus \inj \rho_{\chi^{s},\emptyset})^{I_1}$ extends to the action of 
$\HH$ so that
$$(\inj \rho_{\chi,\emptyset}\oplus \inj \rho_{\chi^{s},\emptyset})^{I_1}\cong L_{\gamma}.$$
\end{prop}
\begin{proof}Suppose that $\chi= \chi^{s}$ by Corollary \ref{qpinjg} we have 
$$(\inj \rho_{\chi,\emptyset}\oplus \inj \rho_{\chi,S})^{I_1}\cong\langle T_{n_s}\varphi_{\chi}\rangle_{\Fbar}\oplus \langle 
(1+T_{n_s})\varphi_{\chi}\rangle_{\Fbar} \cong M_{\chi,\emptyset}\oplus M_{\chi, S}\cong M_{\gamma}|_{\HH_K}$$
as $\HH_K$-modules, where the last isomorphism follows from Lemma \ref{restML}. It is enough to define the action of $T_{\Pi}$. If we let
$$(T_{n_s}\varphi_{\chi})T_{\Pi}= (1+T_{n_s})\varphi_{\chi}\quad\textrm{and}\quad ((1+T_{n_s})\varphi_{\chi})T_{\Pi}= T_{n_s}\varphi_{\chi}$$
then this gives us the required action. Suppose that $\chi\neq \chi^{s}$, then Corollary \ref{qpinjg} and Lemma \ref{restML} imply that
 $$(\inj \rho_{\chi,\emptyset}\oplus \inj \rho_{\chi^{s},\emptyset})^{I_1}\cong (\Indu{I}{K}{\chi^{s}}\oplus \Indu{I}{K}{\chi})^{I_1}
 \cong L_{\gamma}|_{\HH_K}$$
as $\HH_K$-modules.  The space $(\Indu{I}{K}{\chi^{s}})^{I_1}$ has basis $\{T_{n_s}\varphi_{\chi}, \varphi_{\chi^{s}}\}$ and the space
$(\Indu{I}{K}{\chi})^{I_1}$ has basis
$\{ T_{n_s}\varphi_{\chi^{s}}, \varphi_{\chi} \}$. It is enough to define the action of $T_{\Pi}$ on 
the basis. If we set
$$\varphi_{\chi}T_{\Pi}= \varphi_{\chi^{s}},\quad \varphi_{\chi^{s}}T_{\Pi}= \varphi_{\chi}$$
and 
$$(T_{n_s}\varphi_{\chi})T_{\Pi}= T_{n_s}\varphi_{\chi^{s}},\quad (T_{n_s}\varphi_{\chi^{s}})T_{\Pi}= T_{n_s}\varphi_{\chi}$$
then this gives us the required action.
\end{proof}

\subsection{The general case}\label{generalcase}
Our counting argument breaks down if $p\neq q$. The strategy is to restrict to $\SL_2(\Fq)$, where the principal indecomposable representations have
been worked out by Jeyakumar in \cite{j}. Let
$$\Gamma'=\SL_2(\Fq),\quad B'=B\cap \Gamma', \quad H'= H\cap \Gamma'.$$
We note that $U$ is a subgroup of $\Gamma'$ and $n_{s}\in \Gamma'$. 

\subsubsection{Modular representations of $\SL_2(\Fq)$} 
\begin{thm}\label{irrsl} Suppose that $q=p^{n}$. The isomorphism classes of irreducible $\Fbar$-representations of $\Gamma'$ are parameterised by 
$n$-tuples  $\rr=(r_0,\ldots,r_{n-1})$, where  $0\le r_i \le p-1$, for every $i$. Moreover, every irreducible representation can be
 realized over $\Fq$ and the representation  corresponding to an $n$-tuple $\rr$ is given by 
$$ V_{\rr,\Fq}\cong \vs{r_0}{\Fq}\otimes \vs{r_1}{\Fq}^{\Fr}\otimes\ldots\otimes\vs{r_i}{\Fq}^{\Fr^{i}}\otimes\ldots\otimes
\vs{r_{n-1}}{\Fq}^{\Fr^{n-1}}$$
where $V_{r_i,\Fq}$ are the spaces of Section \ref{alternative}.
\end{thm}  
\begin{proof} This is done by Brauer and Nesbitt, see \cite{brauer}.
\end{proof}

\begin{cor}Let $\rho$ be an irreducible representation of $\Gamma$, then $\rho\mid_{\Gamma'}$ is irreducible. Moreover, given an
irreducible representation $\rho'$ of $\Gamma'$ there exist, precisely $q-1$ isomorphism classes of irreducible representations of $\Gamma$, given
by $\rho\otimes (\det)^{a}$, where $0\le a <q-1$, such that
$$(\rho\otimes(\det)^{a})\mid_{\Gamma'}\cong \rho'.$$
\end{cor}
\begin{proof}This is immediate from Theorem \ref{irrsl} and Theorem \ref{irrgl}.
\end{proof}

\begin{remar}\label{extGL}By counting dimensions, we may show that 
$$(\inj (V_{\rr,\Fbar}\otimes (\det)^{a}))\mid_{\Gamma'}\cong \inj V_{\rr,\Fbar}$$
as $\Gamma'$-representations. However, we will obtain this  directly later on.
\end{remar}

We recall the construction of the indecomposable principal representations for $\SL_{2}(\Fq)$ as it is done in \cite{j}. The idea is to go 
from the Lie algebra to the universal enveloping algebra and then to the group.

Let $\lie$ be the Lie algebra of $\SL_2(\CCC)$. It has a $\CCC$-basis
$$e=\begin{pmatrix}0 & 1\\0 & 0 \end{pmatrix},\quad h= \begin{pmatrix}1 & 0\\0 & -1 \end{pmatrix},\quad f=\begin{pmatrix}0 & 0\\1 & 0 \end{pmatrix}.$$  
Let $\UU$ be the universal enveloping algebra of $\lie$. Let $\UU_{\ZZ}$ be a subring of $\UU$ generated by the elements
$$\frac{e^{k}}{k!}, \quad \frac{f^{k}}{k!}, \quad \forall k \in \ZZ^{+}$$
over $\ZZ$. The ring $\UU_{\ZZ}$ has a $\ZZ$-basis, which is also a $\CCC$-basis for $\UU$. Let $d$ be a non-negative integer and let $V_d$ be the 
irreducible module of $\lie$ of highest weight $d$. The space $V_d$ has a $\CCC$-basis of weight vectors $m_i$, for $0\le i \le d$,
 and the action of $\lie$ is given by
$$e m_0=0,\quad  e m_i=(d-i+1)m_{i-1}, \quad 1\le i \le d,$$ 
$$f m_d=0,\quad  f m_i=(i+1)m_{i+1}, \quad 0\le i \le d-1,$$
$$h m_i=(d-2i)m_i,\quad 0\le i \le d.$$
Let $V_{d,\ZZ}$ be a $\ZZ$-lattice in $V_{d}$ spanned by $m_i$, for $0\le i \le d$. We adopt the convention that $m_i=0$ if $i<0$ or $i>d$. Since, 
$$\frac{e^{k}}{k!}m_i= \begin{pmatrix}d-i+k\\ d-i\end{pmatrix} m_{i-k}$$
and
$$\frac{f^{k}}{k!}m_i= \begin{pmatrix}i+k\\ i\end{pmatrix} m_{i+k}$$
for all $k\in \ZZ^{+}$, the lattice $V_{d, \ZZ}$ is a $\UU_{\ZZ}$-module. Let 
$$\tilde{V}_{d,\Fq}=V_{d,\ZZ}\otimes_{\ZZ}\Fq.$$
For every $\lambda\in \Fq$ we define $x(\lambda)$, $y(\lambda)\in \End(\tilde{V}_{d,\Fq})$, by
$$x(\lambda)(v\otimes 1)=\sum_{k\ge 0} \lambda^{k}(\frac{e^{k}}{k!}v\otimes 1)$$
and    
$$y(\lambda)(v\otimes 1)=\sum_{k\ge 0} \lambda^{k}(\frac{f^{k}}{k!}v\otimes 1).$$
Since $e$ and $f$ act nilpotently on $V_d$ this sum is well defined.  There exists a unique homomorphism
$$\SL_2(\Fq)\rightarrow \End(\tilde{V}_{d,\Fq})$$
such that 
$$\begin{pmatrix} 1 & \lambda\\0 & 1 \end{pmatrix}\mapsto x(\lambda)\quad \mbox{and}\quad \begin{pmatrix} 1 & 0\\ \lambda & 1 \end{pmatrix}
\mapsto y(\lambda).$$
This gives us a representation of $\Gamma'$. To ease the notation, we denote
$$m_{i,\Fq}=m_i\otimes 1.$$
We will refer to $\{m_{i,\Fq}: 0\le i \le d\}$ as the standard basis of $\tilde{V}_{d,\Fq}$.  The action of $\Gamma'$ is determined by
$$\begin{pmatrix}1 & \lambda \\ 0 & 1 \end{pmatrix}m_{i,\Fq}=\sum_{k=0}^{i}\begin{pmatrix}d-k\\ d-i\end{pmatrix} \lambda^{i-k} m_{k,\Fq},$$
$$\begin{pmatrix}1 & 0 \\ \lambda & 1 \end{pmatrix}m_{i,\Fq}=\sum_{k=i}^{d}\begin{pmatrix}k\\ i\end{pmatrix} \lambda^{k-i} m_{k,\Fq}.$$
This gives 
$$\begin{pmatrix}\lambda & 0 \\ 0 & \lambda^{-1} \end{pmatrix}m_{i,\Fq}= \lambda^{d-2i} m_{i, \Fq}.$$

At first we resolve the ambiguities in our notation. 
\begin{lem} Let $V_{d,\Fq}$ be a representation of $\Gamma$ constructed in Section \ref{alternative}. Then
$$V_{d,\Fq}|_{\Gamma'}\cong\tilde{V}_{d,\Fq}.$$
\end{lem}
\begin{proof} The isomorphism is given by 
$$m_i\otimes 1\mapsto m_{i,\Fq}.$$
An easy check shows that the isomorphism respects the action of matrices $\begin{pmatrix}1 & \lambda \\ 0 & 1 \end{pmatrix}$ and 
$\begin{pmatrix}1 & 0 \\ \lambda & 1 \end{pmatrix}$, for all $\lambda\in\Fq$. Since, these matrices generate $\Gamma'$ we are done.
\end{proof}
The Lemma above is the reason, why we wanted to work over $\Fq$. We drop the tilde from our notation and go to $\Fbar$. 

For each $r$, such that $0\le r <p-1$, Jeyakumar finds a $\Gamma'$-invariant subspace $R_r$ of
the representation $V_{p-1-r,\Fbar}\otimes V_{p-1,\Fbar}$, such that $\dim R_r =2p$. Let $R_{p-1}=V_{p-1,\Fbar}$, then 
$\dim R_{p-1}=p$. The main result of \cite{j} can be stated as follows.
\begin{thm}\cite{j}\label{jey} Suppose that $q=p^n$. Let $\rr=(r_0,\ldots, r_{n-1})$ be an $n$-tuple, 
such that $0\le r_i\le p-1$, for every $i$. Let
$$R_{\rr}=R_{r_0}\otimes R_{r_1}^{\Fr}\otimes\ldots\otimes R_{r_{n-1}}^{\Fr^{n-1}}.$$
If $\rr\neq\mathbf{0}$, then 
$$R_{\rr}\cong \inj V_{\rr,\Fbar}.$$
And
$$R_{\mathbf{0}}\cong \inj V_{\mathbf{0},\Fbar}\oplus \inj V_{\mathbf{p-1},\Fbar}$$
where $\mathbf{p-1}=(p-1,\ldots,p-1)$ and $\mathbf{0}=(0,\ldots,0)$.
\end{thm}
\begin{remar} Our indices differ slightly from \cite{j}.
\end{remar}

\subsubsection{Going from $\SL_2(\Fq)$ to $\GL_2(\Fq)$}

We will recall how the subspaces $R_r$ are constructed and show that they are in fact $\Gamma$-invariant. That this should be the case is indicated by 
Remark \ref{extGL}. The twisted tensor product will give us principal indecomposable representations 
of $\Gamma$. Since the spaces $R_r$ have a rather concrete description, this will  enable us to work out the corresponding $\HH_{\Gamma}$-modules.

\begin{lem}\label{GLsubSL}Let $V$ be a representation of $\Gamma$ and let $W$ be a $\Gamma'$-invariant subspace of $V$. If $W$ is invariant under the
 action of $H$, then $W$ is $\Gamma$-invariant.
\end{lem}
\begin{proof}Let $v\in W$ and $g\in \Gamma$. We may write $g=g'g_1$, for some $g'\in\Gamma'$ and $g_1\in H$. Then
$$g v=g'(g_1 v)\in W.$$
Hence $W$ is $\Gamma$-invariant.
\end{proof} 

Let $r$ be an integer such that $0\le r\le p-1$. Let $\{v_i\}$, for $0\le i\le p-1-r$ be the standard basis of $V_{p-1-r,\Fbar}$ and let 
$\{w_j\}$, for $0\le j\le p-1$ be the standard basis of $V_{p-1,\Fbar}$. For $0\le i \le 2p-r-2$, we define vectors $E_i$ in 
$V_{p-1-r,\Fbar}\otimes V_{p-1,\Fbar}$, by
$$E_i=\sum_{k+l=i}v_k\otimes w_l.$$
It is convenient to extend the indexing set to $\ZZ$ by setting  $E_i=0$, if $i<0$ or $i> 2p-2p-r$.
\begin{lem}\label{P} The sequence of $\Gamma$ representations
\begin{displaymath}
\xymatrix{0\ar[r]& V_{2p-r-2,\Fbar}\ar[r]& V_{p-1-r,\Fbar}\otimes V_{p-1,\Fbar}\\
                 & m_{i, \Fbar}\ar[r]    & E_i    }
\end{displaymath}
is exact.
\end{lem}
\begin{proof}If $r=p-1$ then this is true trivially. If $r\neq p-1$ then the map is $\Gamma'$-equivariant by \cite{j} Lemma 4.2.
 So by Lemma \ref{GLsubSL} it is enough to show that it is $H$-equivariant. Since
$$\begin{pmatrix}\lambda & 0\\ 0 & \mu \end{pmatrix} m_{i,\Fbar}= \lambda^{2p-r-2-i}\mu^{i} m_{i,\Fbar}$$
and
$$\begin{pmatrix}\lambda & 0\\ 0 & \mu \end{pmatrix} E_i= \lambda^{2p-r-2-i}\mu^{i} E_i$$
we are done.
\end{proof}

\begin{defi}\label{defRr}\cite{j} Let $r$ be an integer, such that $0\le r <p-1$. For $0\le i \le p-r-1$, let $a_i$ be integers defined by the following 
relation:
$$a_0=0\quad \textrm{and}\quad  a_1= (p-r-2)!$$ 
and
$$a_{i+1}=a_i+\frac{(-1)^{i}(r+1)\ldots (r+i)}{(p-r-2)\ldots (p-r-i-1)}(a_1-a_0).$$
Let $Z$ be a vector in $V_{p-1-r,\Fbar}\otimes V_{p-1,\Fbar}$ given by
$$Z=a_0 (v_0\otimes w_{p-r-1})+a_1(v_1\otimes w_{p-r-2})+\ldots +a_{p-r-1}(v_{p-r-1}\otimes w_{0}),$$
and let $R_{r}$ be a subspace of $V_{p-1-r,\Fbar}\otimes V_{p-1,\Fbar}$ given by
$$R_{r}=\langle E_0,\ldots, E_{2p-r-2}, Z, \frac{f}{1!}Z, \ldots, \frac{f^{r}}{r!}Z\rangle_{\Fbar}.$$
Moreover, for $r=p-1$ we define
$$R_{p-1}= V_{p-1,\Fbar}.$$
\end{defi}  

\begin{prop}\label{Rr} Let $r$ be an integer, such that $0\le r \le p-1$, then $R_r$ is a $\Gamma$-invariant subspace of 
$V_{p-r-1,\Fbar}\otimes V_{p-1,\Fbar}$.
Moreover, if $r\neq p-1$, then 
$$\dim R_r =2p$$
and if $r=p-1$, then
$$ \dim R_{p-1}=p.$$
\end{prop}
\begin{proof}If $r=p-1$ then there is nothing to prove, since $R_{p-1}=V_{p-1,\Fbar}.$
If $r\neq p-1$ then by \cite{j} Theorem 4.7 $R_r$ is $\Gamma'$-invariant and $\dim R_r =2p$. So by Lemma \ref{GLsubSL} it is enough to show that
$R_r$ is $H$-invariant. For $v\in V_{p-r-1,\Fbar}$ and $w\in V_{p-1,\Fbar}$ we have
$$f(v\otimes w)= f v \otimes w+ v\otimes f w.$$ 
Hence, for $0\le k \le r$ we have
$$\frac{f^{k}}{k!}Z\in \langle v_{l+i}\otimes w_{p-r-1-l+j}\mid i+j=k,\quad 0\le l \le  p-r-1\rangle_{\Fbar}$$
with the usual 'vanishing when not defined' convention. Since
\begin{equation}\notag
\begin{split}
\begin{pmatrix}\lambda & 0 \\ 0 & \mu \end{pmatrix} v_{l+i}\otimes w_{p-r-1-l+j}=& \lambda^{p-r-1-l-i}\mu^{l+i}\lambda^{r+l-j}\mu^{p-r-1-l+j}
v_{l+i}\otimes w_{p-r-1-l+j}\\=&\lambda^{p-k-1}\mu^{p-r-1+k} v_{l+i}\otimes w_{p-r-1-l+j}
\end{split}
\end{equation}
the group $H$ acts on each $\frac{f^{k}}{k!}Z$, for $0\le k \le r$  by a character. We combine this with Lemma \ref{P} and obtain that $R_r$ is $H$ 
invariant.
\end{proof}

\begin{lem}\label{fe} We have 
$$\frac{f^{k}}{k!}E_{p-r-1}=0$$
if and only if $k\ge r+1$. For $k\ge 1$ we have
 $$\frac{e^{k}}{k!}E_{p-r-1}=0.$$
In particular, $U$ fixes $E_{p-r-1}$ and the action of $H$ is given by
$$\begin{pmatrix} \lambda & 0 \\ 0 & \mu \end{pmatrix} E_{p-r-1}= \lambda^{r} (\lambda \mu)^{p-r-1} E_{p-r-1}.$$
\end{lem}
\begin{proof}If $r=p-1$, then this is trivial. If $r\neq p-1$ then for $k\ge 0$ we have 
$$\frac{f^{k}}{k!}E_{p-r-1}=\begin{pmatrix}p-r-1+k\\ p-r-1\end{pmatrix} E_{p-r-1+k}.$$
We observe that $E_{p-r-1+k}$ vanishes trivially, if $k\ge p$. If $r+1\le k \le p-1$, then we write $k=r+1+j$, where $0\le j \le p-r-2$. The
 binomial coefficient becomes
$$\begin{pmatrix} j+p \\ p-r-1 \end{pmatrix}.$$
Since $0\le r < p-1$, we have $1\le p-r-1\le p-1$, and since $0\le j < p-r-1$, $p$ divides the binomial coefficient. Hence
$$\frac{f^{k}}{k!}E_{p-r-1}=0$$
for $k\ge r+1$. If $0\le k \le r$, then $p-r-1\le p-r-1+k \le p-1$ and the binomial coefficient does not vanish. Hence
$$\frac{f^{k}}{k!}E_{p-r-1}\neq0$$
for $0\le k \le r$. Let $k\ge 0$, then 
$$\frac{e^{k}}{k!}E_{p-r-1}=\begin{pmatrix}p-1+k\\ p-1\end{pmatrix} E_{p-r-1-k}.$$
We observe that $E_{p-r-1-k}$ vanishes trivially, if $k> p-r-1$. Suppose that $1\le k \le p-r-1$, then  we may write $k=j-1$, where 
$0\le j \le p-r-2 <p-1$. The binomial coefficient becomes
$$\begin{pmatrix}j+p\\ p-1\end{pmatrix}.$$
Since $j<p-1$, $p$ divides the binomial coefficient, and hence 
$$\frac{e^{k}}{k!}E_{p-r-1}=0$$
for all $k\ge 1$. Since the action of $U$ is given in terms of  $\frac{e^{k}}{k!}$ this  implies that $U$ fixes $E_{p-r-1}$. An easy verification 
gives us the action of $H$.
\end{proof}  

\begin{prop}\label{Ep} Let $W_r$ be a subspace of $R_r$ given by
$$W_r= \langle E_{p-r-1}, \ldots , E_{p-1}\rangle_{\Fbar}.$$
Then $W_r$ is $\Gamma$-invariant. Moreover, 
$$W_{r}^{U}= \langle E_{p-r-1}\rangle_{\Fbar}$$
and
$$W_r= \langle \Gamma E_{p-r-1}\rangle_{\Fbar} \cong V_{r,\Fbar}\otimes (\det)^{p-r-1}.$$ 
\end{prop}
\begin{proof}If $r=p-1$ then $W_{p-1}\cong V_{p-1,\Fbar}$
and we are done. Otherwise, since $W_r$ has a basis of eigenvectors for the action of $H$, it is enough to show 
that $W_r$ is $\Gamma'$-invariant. Since the 
action of $\Gamma'$ is given in terms of the action of $\UU_{\ZZ}$ it is enough to show that $W_r$ is invariant under the action of $\UU_{\ZZ}$.
Lemma \ref{fe} implies that $W_r$ has a basis $\frac{f^k}{k!}E_{p-r-1}$, for $0\le k \le r$. We observe that Lemma \ref{fe} also implies that
$$\frac{f^{l}}{l!}(\frac{f^{k}}{k!}E_{p-r-1})=\begin{pmatrix}k+l\\ k\end{pmatrix}\frac{f^{k+l}}{(k+l)!}E_{p-r-1}\in W_r$$
for $0\le k \le r$ and $l\ge 0$. Suppose that $0\le k \le r$ and $l\ge k+1$ then
$$\frac{e^{l}}{l!}(\frac{f^{k}}{k!}E_{p-r-1})=0.$$
This follows from the multiplication in $\UU_{\ZZ}$, see \cite{hump} \S 26.2, and Lemma \ref{fe}. If $0\le  l  \le k \le r$, then
\begin{equation}\notag
\begin{split}
\frac{e^{l}}{l!}(\frac{f^{k}}{k!}E_{p-r-1})=&\begin{pmatrix}p-r-1+k\\ p-r-1\end{pmatrix}\frac{e^{l}}{l!}E_{p-r-1+k}\\
=&\begin{pmatrix}p-r-1+k\\ p-r-1\end{pmatrix}\begin{pmatrix}p-1-k+l\\ p-1-k \end{pmatrix}E_{p-r-1+k-l}\in W_r.
\end{split}
\end{equation}
Hence $W_r$ is invariant under the  action of $\UU_{\ZZ}$ and hence under the action of $\Gamma$. 

We know from Lemma \ref{fe} that $E_{p-r-1}$ is fixed by $U$.  The action of $H$ splits $W_{r}^{U}$ into a direct
 sum of one dimensional subspaces. Suppose that $\dim W_{r}^U\ge 2$. Since $H$ acts on each vector $E_{p-r-1+k}$
by a distinct character for $0\le k \le r$, we must have $E_{p-r-1+j}\in W_r^{U}$, for  some $1\le j \le r$. This implies that 
$$ e E_{p-r-1+j}=(p-j) E_{p-r-2+j}=0.$$
Hence $p$ must divide $j$ and this is impossible. Hence, $\dim W_{r}^{U}=1$.

Since $W_r$ is $\Gamma$-invariant, we have 
$$\langle \Gamma E_{p-r-1}\rangle_{\Fbar}\le W_r.$$
We may choose $r+1$ distinct elements $\lambda_i$ in $\Fq$. Then
$$\begin{pmatrix}1 & 0 \\ \lambda_i & 1\end{pmatrix} E_{p-r-1}= \sum_{k=0}^{r} \lambda_{i}^{k} \frac{f^{k}}{k!} E_{p-r-1}.$$
Let $A$ be an $(r+1)\times (r+1)$ matrix, given by $A_{ k i}= \lambda_{i}^{k}$, for $0\le i, k \le r$, with the convention that $0^0=1$. 
 Then $\det A$ is the 
Vandermonde determinant, which is non-zero, since all the $\lambda_i$ are distinct. Hence, $A$ is invertible and
$$\frac{f^k}{k!}E_{p-r-1}\in \langle \Gamma E_{p-r-1} \rangle_{\Fbar}$$
for all $0\le k \le r$. Hence, $W_r= \langle \Gamma E_{p-r-1}\rangle_{\Fbar}$. 

Since $\dim W_{r}^{U}=1$ and $W_r=\langle \Gamma W_r^{U}\rangle_{\Fbar}$, the representation $W_r$ is irreducible. To decide, which one it is, we may 
proceed as in the proof of Proposition \ref{dictionary}. Since $r<p-1$, the action of 
$B$ on $W_{r}^{U}$  implies that $W_r\cong V_{r,\Fbar}\otimes(\det)^{p-r-1}$.
\end{proof}

\begin{lem}\label{Ez}The vector $E_0$ is fixed by the action of $U$. Moreover, $H$ acts on $E_0$ by
$$\begin{pmatrix}\lambda & 0 \\ 0 & \mu \end{pmatrix}E_0= \lambda^{r} (\lambda \mu)^{p-r-1} (\lambda \mu^{-1})^{p-r-1} E_0.$$
\end{lem}
\begin{proof}Since, $E_0= v_0\otimes w_0$ this is immediate.
\end{proof}    

\begin{defi} Suppose that $q=p^{n}$ and let $\rr=(r_0,\ldots, r_{n-1})$ be the $n$-tuple such that $0\le r_i\le p-1$, then we define
a representation $R_{\rr}$ of $\Gamma$, given by
$$R_{\rr}=R_{r_{0}}\otimes R_{r_{1}}^{\Fr}\otimes\ldots \otimes R_{r_{n-1}}^{\Fr^{n-1}}$$
where $R_{r_{i}}$ are $\Gamma$-representations of  Definition \ref{defRr}.
\end{defi}

\begin{defi}\label{epsilon} Suppose that $q=p^n$ and let $\rr=(r_0,\ldots, r_{n-1})$ be an $n$-tuple, such that $0\le r_i \le p-1$, for every $i$. Let
 $\mathbf{\varepsilon}=(\epsilon_0, \ldots, \epsilon_{n-1})$ be an $n$-tuple, such that $\epsilon_i\in\{0,1\}$ for every $i$. We define a vector
$$b_{\mathbf{\varepsilon}}=E_{(1-\epsilon_0)(p-1-r_0)}\otimes\ldots\otimes E_{(1-\epsilon_{n-1})(p-1-r_{n-1})}$$
in $R_{\rr}$, where $E_{(1-\epsilon_i)(p-1-r_i)}$ is a vector in $R_{r_i}$, for each $0\le i \le n-1$.
\end{defi}     

\begin{defi}\label{Sigmar} Suppose that $q=p^n$ and let $\rr=(r_0,\ldots, r_{n-1})$ be the $n$-tuple, such that $0\le r_i \le p-1$, for $0\le i \le n-1$.
We define $\Sigma_{\rr}$ to be the set of $n$-tuples $(\epsilon_0,\ldots, \epsilon_{n-1})$, such that
$$\epsilon_i=0,\quad \textrm{if $r_i=p-1$ and}\quad \epsilon_i\in\{0,1\}, \quad \textrm{otherwise}.$$
We will write $\mathbf{0}=(0,\ldots, 0)$ and $\mathbf{1}=(1, \ldots, 1)$.
\end{defi}

\begin{remar}We hope to prevent some notational confusion. Since we want Lemma \ref{mrfixit} to hold and since $\dim R_{p-1,\Fbar}^{U}=1$, 
if $r_i=p-1$, we have to make a choice for $\epsilon_i$, between  $0$ and $1$. We choose $0$, since then we can state Lemma \ref{genezero}
in a nice way. However, if $r_i=p-1$, then
$$(1-0)(p-r_i-1)=(1-1)(p-r_i-1)=0$$
so it does not matter, whether $\epsilon_i=0$ or $\epsilon_i=1$, and we will exploit this in our notation. We note that the definition of 
$b_{\mathbf{\varepsilon}}$ is independent of the set $\Sigma_{\rr}$ and we might have $\mathbf{\varepsilon}\in \Sigma_{\rr}$, 
$\mathbf{\varepsilon'}\not\in \Sigma_{\rr}$, but $b_{\mathbf{\varepsilon}}=b_{\mathbf{\varepsilon'}}$.
\end{remar}  

\begin{lem}\label{mrfixit} The set $\{b_{\mathbf{\varepsilon}}: \mathbf{\varepsilon}\in \Sigma_{\rr}\}$ is  a basis of $R_{\rr}^{U}$. 
\end{lem}
\begin{proof}Let $r$ be an integer, such that $0\le r\le p-1$. If $r=p-1$ , then $\dim R_r=p$ and $E_0$ is in $R_{r}^{U}$. If 
$0\le r <p-1$, then $\dim R_r=2p$ and $E_0$ and $E_{p-1-r}$ are two linearly independent vectors in $R_{r}^{U}$. 

Let $\rr$ be an $n$-tuple. Then by above vectors $b_{\mathbf{\varepsilon}}$, for $\mathbf{\varepsilon}\in \Sigma_{\rr}$,
span a linear subspace of $R_{\rr}^{U}$ of dimension $|\Sigma_{\rr}|$. Also by above, $\dim R_{\rr}=|\Sigma_{\rr}|q$. Since, $U$ is a
$p$-Sylow subgroup of $\Gamma'$ of order $q$ and by Theorem \ref{jey} $R_{\rr}$ is an injective object in $\Rep_{\Gamma'}$, Corollary \ref{dim}
implies that 
$$\dim R_{\rr}^{U}= |\Sigma_{\rr}|.$$
Hence, the set  $\{ b_{\mathbf{\varepsilon}}: \varepsilon \in \Sigma_{\rr} \}$ is a basis of $R_{\rr}^{U}$.
\end{proof} 

\begin{lem}\label{mrfixitH} Let $\rr=(r_0,\ldots,r_{n-1})$ be an $n$-tuple, with $0\le r_i \le p-1$,  
let $\mathbf{\varepsilon}=(\epsilon_0,\ldots,\epsilon_{n-1})$ be an $n$-tuple such that $\epsilon_i\in \{0,1\}$, for every $i$, and let
 $b_{\mathbf{\varepsilon}}$ be a vector in $R_{\rr}^{U}$, then the action of 
$H$ is given by 
$$\begin{pmatrix}\lambda & 0 \\ 0 & \mu \end{pmatrix} b_{\mathbf{\varepsilon}}=\lambda^r (\lambda\mu)^{q-1-r}
(\lambda \mu^{-1})^{\mathbf{\varepsilon}\centerdot (\mathbf{p-r-1})} b_{\mathbf{\varepsilon}}$$
where $r=r_0+r_1p+\ldots r_{n-1}p^{n-1}$ and
$$\mathbf{\varepsilon}\centerdot (\mathbf{p-r-1})= \epsilon_0(p-r_0-1)+\epsilon_1(p-r_1-1)p+\ldots+\epsilon_{n-1}(p-r_{n-1}-1)p^{n-1}.$$
\end{lem}
\begin{proof}This follows from Proposition \ref{fe} and Lemma \ref{Ez}. We note that the action on each tensor component is twisted by $\Fr$.
\end{proof} 

\begin{lem}\label{genezero}Suppose that $q=p^{n}$ and let $\rr=(r_0,\ldots, r_{n-1})$ be an $n$-tuple,
 such that $0\le r_i \le p-1$, for each $i$. Let 
$b_{\mathbf{0}}$ be a vector in $R_{\rr}$. Let 
$$r=r_0+r_1p+\ldots r_{n-1} p^{n-1}.$$
 Then
$$\langle \Gamma b_{\mathbf{0}}\rangle_{\Fbar}\cong V_{\rr, \Fbar}\otimes(\det)^{q-1-r}$$
as a $\Gamma$-representation.
\end{lem}
\begin{proof}Let $W_{\rr}$ be the subspace of $R_{\rr}$ given by 
$$W_{\rr}= W_{r_0}\otimes\ldots \otimes W_{r_{n-1}}$$
with the notation of the Proposition \ref{Ep}. We have 
$$0\neq\langle \Gamma b_{\mathbf{0}}\rangle_{\Fbar}\le W_{\rr}.$$
Proposition \ref{Ep} applied to every tensor component implies that
$$W_{\rr}\cong V_{\rr, \Fbar}\otimes(\det)^{q-1-r}$$
which is irreducible. Hence, we must get the whole of $W_{\rr}$. 
\end{proof}

\begin{cor}\label{emptyseteasy} Let $\chi:H\rightarrow \Fbar^{\times}$ and 
 let $a$ and $r$ be  unique integers, such that $1\le a, r\le q-1$ and
$$ \chi(\begin{pmatrix} 1 & 0 \\0 & \lambda \end{pmatrix})= \lambda^{a}\quad \forall \lambda\in \Fq^{\times},\quad
 \chi(\begin{pmatrix} \lambda & 0 \\0 & \lambda^{-1} \end{pmatrix})= \lambda^{r}\quad \forall \lambda\in \Fq^{\times}.$$
Let $r=r_0+r_1p+\ldots + r_{n-1}p^{n-1}$, where $0\le r_i \le p-1$ for each $i$, and let
$\rr=(r_0,\ldots,r_{n-1})$. If $\chi\neq \chi^s$ then
$$\inj \rho_{\chi, \emptyset}\cong R_{\rr}\otimes(\det)^{a+r}.$$
If $\chi=\chi^{s}$ then
$$\inj \rho_{\chi, \emptyset}\cong R_{\mathbf{p-1}}\otimes(\det)^{a}\cong V_{\mathbf{p-1},\Fbar}\otimes (\det)^{a}.$$
\end{cor}
\begin{proof} Lemma \ref{genezero} implies the existence of an exact sequence
\begin{displaymath}
\xymatrix{0\ar[r]& V_{\rr,\Fbar}\ar[r]& R_{\rr}\otimes (\det)^r}
\end{displaymath}
 of $\Gamma$-representations. It is enough to show that $R_{\rr}$ is  an indecomposable injective object in $\Rep_{\Gamma}$. The rest follows 
from Propositions \ref{dictionary} and  \ref{injrho}.

 Theorem \ref{jey} says that the restriction of $R_{\rr}$ to $\Gamma'$ is indecomposable. In particular, $R_{\rr}$  must be indecomposable as 
a $\Gamma$-representation. Moreover, Theorem \ref{jey} says that the restrictions of $R_{\rr}$ to $\Gamma'$ is an injective object
in $\Rep_{\Gamma'}$. Since $U$ is a $p$-Sylow subgroup of both $\Gamma$ and $\Gamma'$, Proposition \ref{restsylow} implies that
$R_{\rr}$  is an injective object in $\Rep_{\Gamma}$. Finally, the last isomorphism follows directly from the definition
of $R_{\mathbf{p-1}}$. 
\end{proof}

\subsubsection{Computation of $\HH_{\Gamma}$-modules}
We will compute the action of $T_{n_s}$ on $R_{\rr}^{U}$.
\begin{prop}\label{hardcore}Let $q=p^{n}$ and let $\rr=(r_0, \ldots, r_{n-1})$ be the $n$-tuple, such that $0\le r_i \le p-1$, for every $i$.
Let $\mathbf{\varepsilon}\in \Sigma_{\rr}$ and let $b_{\mathbf{\varepsilon}}$ be a vector in $R_{\rr}$.
\begin{itemize}
\item[(i)] Suppose that for some index $j$,  $\epsilon_j=0$ and $r_j\neq p-1$ then
$$\sum_{u\in U} u n_{s}^{-1}b_{\mathbf{\varepsilon}}=0.$$
\item[(ii)] Suppose that $\rr\neq \mathbf{0}$. Moreover, suppose that for every $i$, if $\epsilon_i=0$, then 
$r_i=p-1$  then $b_{\mathbf{\varepsilon}}=b_{\mathbf{1}}$ and
$$\sum_{u\in U} u n_{s}^{-1}b_{\mathbf{1}}=(-1)^{1+|\rr|}b_{\mathbf{0}}$$
where $|\rr|=  r_0+  r_1 p +\ldots +  r_{n-1} p^{n-1}.$  
\item[(iii)] Suppose that $\rr=\mathbf{0}$ and $\mathbf{\varepsilon}=\mathbf{1}$, then
$$\sum_{u\in U}u n_{s}^{-1}b_{\mathbf{1}}=-(b_{\mathbf{0}}+b_{\mathbf{1}}).$$
\end{itemize}
This covers all the possible pairs $(\rr, \varepsilon)$, such that $\varepsilon\in \Sigma_{\rr}$.
\end{prop}
\begin{remar} We note that $b_{\mathbf{1}}$ is well defined even if $\mathbf{1}\not \in \Sigma_{\rr}$. See Definitions \ref{epsilon} and \ref{Sigmar}.
\end{remar}
\begin{proof}  Let $r$ be an integer such that $0\le r \le p-1$ and let $\epsilon \in \{0,1\}$ such that  $\epsilon=0$, if $r=p-1$. Let 
$E_{(1-\epsilon)(p-r-1)}$ be a vector in $R_r$. We observe that
$$n_{s}^{-1}E_{(1-\epsilon)(p-r-1)}=(-1)^{p-1+\epsilon(p-r-1)}E_{p-1+\epsilon(p-r-1)}.$$
If $r\neq p-1$ this follows from Lemma \ref{P},   and 
if $r=p-1$, this follows from the isomorphism $R_{p-1}\cong V_{p-1,\Fbar}$. Moreover, 
$$\frac{e^k}{k!}E_{p-1}=0,\quad\textrm{if $k> r$}\quad\textrm{and}\quad \frac{e^{k}}{k!}E_{2p-2-r}=0,\quad\textrm{if $k> 2p-2-r$. }$$

Let $\rr$ be an $n$-tuple, $\rr=(r_0, \ldots, r_{n-1})$ and let $\mathbf{\varepsilon}\in \Sigma_{\rr}$. Then 
\begin{equation}\notag
\begin{split}
\sum_{u\in U}u n_{s}^{-1} b_{\mathbf{\varepsilon}}=\sum_{k_0,\ldots, k_{n-1}\ge 0} (-1)^{q-1-\mathbf{\varepsilon}\centerdot (\mathbf{p-1}-\rr)}
\sum_{\lambda\in \Fq} \lambda^{k_0+\ldots+k_{n-1}p^{n-1}}\\
 \frac{e^{k_0}}{k_0!}E_{p-1+\epsilon_0(p-r_0-1)}\otimes\ldots\otimes \frac{e^{k_{n-1}}}{k_{n-1}!}E_{p-1+\epsilon_{n-1}(p-r_{n-1}-1)}
\end{split}
\end{equation}
where $\mathbf{\varepsilon}\centerdot (\mathbf{p-1}-\rr)=\sum_{i=0}^{n-1} \epsilon_i(p-1-r_i)p^{i}$.
We have acted by $n_{s}^{-1}$ on each tensor 
component and then expanded the action of $u\in U$ on each tensor component and rearranged the summation. We will show that the 
the terms in the sum vanish, unless 
$$(k_0,\ldots,k_{n-1})=(p-1,\ldots,p-1)\quad \textrm{or} \quad (k_0,\ldots,k_{n-1})=(2(p-1),\ldots, 2(p-1))$$
and $\rr$ and $\mathbf{\varepsilon}$ are of a special form.

\textit{Step 1.} We claim that  if $\epsilon_i=1$ then it is enough to consider $k_i=r_i$ and if $\epsilon_i=1$, then it is enough to consider
 $k_i=p-1$ and $k_i=2p-r_i-2$, since all the other terms in the sum vanish. 
 
We observe that for each $i$, if $\epsilon_i=0$ then it is enough to consider $0\le k_i\le r_i$ and if $\epsilon_i=1$ then
it is enough to consider $0\le k_i \le 2p-2-r_i$. This follows by  looking at a single tensor component as above. 
Moreover, we observe that
$$\frac{e^{k_i}}{k_i !}E_{p-1+\epsilon_i(p-r_i-1)}\in \langle E_{p-1+\epsilon_i(p-r_i-1)-k_i}\rangle_{\Fbar}.$$
The vector $\sum_{u\in U} u n_{s}^{-1}b_{\mathbf{\varepsilon}}$ is fixed by $U$. By Lemma \ref{mrfixit} vectors $b_{\mathbf{\varepsilon'}}$, for 
$\mathbf{\varepsilon'}\in \Sigma_{\rr}$, form a basis of $R_{\rr}^{U}$. Hence, for each $i$, it is enough to consider $k_i$ of the form
$$k_i=p-1+(\epsilon_i-1)(p-r_i-1)\quad\textrm{and}\quad k_i=p-1+\epsilon_i(p-r_i-1)$$
 since all the other terms must vanish. If $\epsilon_i=0$ and  $k_i$ is of the form as
above then the inequality $k_i\le r_i$ can be fulfilled if and only if $k_i=r_i$. If $\epsilon_i=1$, then $k_i\le 2p-r_i-2$ implies
that $k_i=p-1$ or $k_i=2p-r_i-2$. 

\textit{Step 2.} Let $k=k_0+k_1 p+\ldots+k_{n-1} p^{n-1}$. We claim that if $\varepsilon=\mathbf{1}$ and $\rr=\mathbf{0}$, then it is enough to consider 
two cases $k=q-1$ and 
$$(k_0,\ldots,k_{n-1})=(2(p-1),\ldots, 2(p-1))$$
 and otherwise it is enough to consider the case $k=q-1$, since all the other terms in the sum vanish.
 
Step 1 implies that it is enough to consider $n$-tuples $(k_0,\ldots,k_{n-1})$, such that $0\le k \le 2(q-1)$.
Moreover, the upper bound is obtained if and only if $\rr=\mathbf{0}$, $\mathbf{\varepsilon}=\mathbf{1}$ and 
$(k_0, \ldots, k_{n-1})=(2(p-1), \ldots, 2(p-1)).$
If $k=0$ or $k> 0$ and $q-1$ does not divide $k$ then  
$$\sum_{\lambda\in \Fq}\lambda^{k}=0.$$
We note that $0^{0}=1$ comes from the action by the identity matrix. If $k> 0$ and $q-1$ divides $k$, then
$$\sum_{\lambda\in \Fq}\lambda^{k}=-1.$$
This establishes the claim.

\textit{Step 3.} We claim that if  $k=q-1$, then it is enough to consider $k_i=p-1$, for every $i$, since all the other terms in the sum vanish. 

We use Step 1 to define integers $a_i$ and $a'_i$, such that for each $i$ 
$$a_i+a'_i=k_i$$
and $0\le a_i, a'_i\le p-1$, as follows. If $\epsilon_i=0$, then $a_i=r_i$ and $a'_i=0$. If 
$\epsilon_i=1$ and $k_i=p-1$, then $a_i=p-1$ and $a'_i=0$. If $\epsilon_i=1$ and $k_i=2p-r_i-2$, then $a_i=p-1$ and $a'_i=p-1-r_i$. Then
$q-1=k$ implies that
$$ a_0+a_1 p+\ldots+ a_{n-1}p^{n-1}= (p-1-a'_0)+ (p-1-a'_1) p+\ldots+ (p-1-a'_{n-1} )p^{n-1}.$$
Since $0\le a_i, a'_i \le p-1$, for every $i$, this implies 
$$a_i=p-1-a'_i,\quad \forall i.$$  
 If $\epsilon_i=1$ and $k_i=p-1$, then we are done. Otherwise if $\epsilon_i=0$ or $\epsilon_i=1$ and $k_i=2p-2-r_i$ then above implies that
$r_i=p-1$. This establishes the claim.

\textit{Step 4.} Suppose that for some index $j$, $\epsilon_j=0$ and $r_j\neq p-1$, then Steps 1, 2 and 3 imply that all the terms vanish. So we 
obtain part (i) of the Proposition. We note that this case includes $\rr=\mathbf{0}$ and $\mathbf{\varepsilon}\neq \mathbf{1}$.
 
\textit{Step 5.} Suppose that $\rr\neq \mathbf{0}$. Moreover, suppose that for every $i$, if $\epsilon_i=0$ then $r_i=p-1$. 
 We will compute what happens on each tensor component if $k_i=p-1$. If  $\epsilon_i=0$, then  $r_i=p-1$ and
$$\frac{e^{p-1}}{(p-1)!}E_{p-1}=\begin{pmatrix}p-1  \\ 0 \end{pmatrix} E_{0}= E_{p-r_i-1}.$$
If $\epsilon_i=1$ then
$$\frac{e^{p-1}}{(p-1)!}E_{2p-2-r_i}=\begin{pmatrix}p-1  \\ 0 \end{pmatrix} E_{p-r_i-1}=E_{p-r_i-1}.$$
The above calculation gives us 
$$\sum_{u\in U} u n_{s}^{-1} b_{\mathbf{\varepsilon}}= (-1)^{\mid\rr\mid+ 1}b_{\mathbf{0}}.$$
Since, by Steps 2 and 3, it is enough to consider a single term in the sum 
$$(k_0,\ldots,k_{n-1})=(p-1,\ldots,p-1)$$
and by Definition \ref{epsilon}, $b_{\mathbf{0}}= E_{p-r_0-1}\otimes\ldots\otimes E_{p-r_{n-1}-1}$. Moreover, if $p=2$, then $1=-1$ and if 
$p\neq 2$ then
$$(-1)^{p-1+\epsilon_i(p-1-r_i)}= (-1)^{r_i}$$
trivially, if $\epsilon_i=1$ and since $r_i=p-1$ if $\epsilon_i=0$. We get an extra $-1$ from summing over $\lambda\in \Fq$. This accounts for the sign.
We claim that in this case $b_{\mathbf{\varepsilon}}=b_{\mathbf{1}}$. Indeed, if $r_i\neq p-1$ then $\epsilon_i=1$ and  if $r_i=p-1$, then
$$(1-\epsilon_i)(p-1-r_i)= (1-1)(p-1-r_i)=0.$$
Hence, $b_{\varepsilon}=b_{\mathbf{1}}$, see \ref{epsilon}. This establishes part (ii) of the Proposition.

\textit{Step 6.} The only case left is  $\rr=\mathbf{0}$ and $\mathbf{\varepsilon}=\mathbf{1}$. The only difference to Step 5 
is that we get a contribution from $(k_0,\ldots,k_{n-1})=(2(p-1),\ldots, 2(p-1))$. More, precisely
$$\frac{e^{2p-2}}{(2p-2)!}E_{2p-2}=\begin{pmatrix}2p-2  \\ 0 \end{pmatrix} E_{0}=E_{0}.$$
And by Definition \ref{epsilon}, $b_{\mathbf{1}}=E_0\otimes\ldots \otimes E_0.$ Hence,
$$\sum_{u\in U} u n_{s}^{-1} b_{\mathbf{1}}= -(b_{\mathbf{1}}+b_{\mathbf{0}}).$$
The minus sign comes from summing over $\lambda\in \Fq$. This establishes part (iii) of the Proposition.
\end{proof} 

\begin{remar} We think of $\otimes (\det)^{a}$ as a twist, that is, it changes the action, but does not change the underlying vector space. Moreover,
since $U\le \Gamma'$ and $n_s\in \Gamma'$, Proposition \ref{hardcore} does not change if we twist the action by $(\det)^{a}$.
\end{remar} 
\begin{remar} We know that something like
$$\sum_{u\in U} u n_{s}^{-1}b_{\mathbf{1}}= (-1)^{1+|\rr|}b_{\mathbf{0}}$$
has to happen by Lemma \ref{subsp}.
\end{remar}
\begin{lem}\label{gensteinberg} Let $b_{\mathbf{1}}$ and $b_{\mathbf{0}}$ be  vectors in $R_{\mathbf{0}}$. Then
$$\langle \Gamma (b_{\mathbf{1}}+b_{\mathbf{0}})\rangle_{\Fbar}\cong V_{\mathbf{p-1},\Fbar}.$$
\end{lem}
\begin{proof}The vector $b_{\mathbf{1}}+ b_{\mathbf{0}}$ is fixed by $U$. Moreover, by Lemma \ref{mrfixitH} $H$ acts trivially on it.
By Proposition \ref{hardcore}
$$(b_{\mathbf{1}}+b_{\mathbf{0}})T_{n_s}=\sum_{u\in U} u n_{s}^{-1}(b_{\mathbf{1}}+b_{\mathbf{0}})=-(b_{\mathbf{1}}+b_{\mathbf{0}}).$$
Hence 
$$\langle b_{\mathbf{1}}+b_{\mathbf{0}} \rangle_{\Fbar}\cong M_{\Eins, \emptyset}$$
as $\HH_{\Gamma}$-module and  Lemma \ref{splitisgood} gives us the result.
\end{proof}

\begin{cor}\label{sleepy}Let $\chi:H\rightarrow \Fbar^{\times}$ be a character, such that $\chi=\chi^{s}$  and 
 let $a$ be the unique integer, such that $1\le a\le q-1$ and
$$ \chi(\begin{pmatrix} 1 & 0 \\0 & \lambda \end{pmatrix})= \lambda^{a}\quad \forall \lambda\in \Fq^{\times}$$
then
$$\inj \rho_{\chi, S}\oplus \inj\rho_{\chi,\emptyset}\cong R_{\mathbf{0}}\otimes (\det)^{a}.$$
\end{cor}
\begin{proof} This is a rerun of the proof of Corollary \ref{emptyseteasy}. Lemma \ref{genezero} and Lemma \ref{gensteinberg} 
imply the existence of an exact sequence
\begin{displaymath}
\xymatrix{0\ar[r]& V_{\mathbf{0},\Fbar}\oplus V_{\mathbf{p-1},\Fbar}\ar[r]& R_{\mathbf{0}}}
\end{displaymath}
 of $\Gamma$-representations.
So it is enough to show that $R_{\mathbf{0}}$ is an injective object in $\Rep_{\Gamma}$ and that it has at most $2$ direct summands. The rest 
follows from Proposition \ref{injrho} and Proposition \ref{dictionary}. Theorem \ref{jey} says that the restriction of 
$R_{\mathbf{0}}$ to $\Gamma'$ has exactly $2$ direct
summands, hence $R_{\mathbf{0}}$ may have at most $2$ direct summands. Moreover, Theorem \ref{jey} says that the restriction of $R_{\mathbf{0}}$
to $\Gamma'$ is an injective object in $\Rep_{\Gamma'}$. Since $U$ is a $p$-Sylow subgroup of $\Gamma$ and $\Gamma'$ contains $U$, Proposition
\ref{restsylow} implies that $R_{\mathbf{0}}$ is an injective object in $\Rep_{\Gamma}$.
\end{proof}

\begin{defi} Let $\alpha: H\rightarrow \Fbar^{\times}$ be a character, given by 
$$\alpha: \begin{pmatrix} \lambda & 0 \\ 0 & \mu\end{pmatrix} \mapsto \lambda \mu^{-1}.$$
\end{defi}

\begin{lem}\label{twistbys} Suppose that $q=p^{n}$ and let $\chi:H\rightarrow \Fbar^{\times}$ be a character. Let $r$ be the unique integer, 
such that $0\le r < q-1$ and
$$ \chi(\begin{pmatrix} \lambda & 0 \\0 & \lambda^{-1} \end{pmatrix})= \lambda^{r}\quad \forall \lambda\in \Fq^{\times}.$$
Let $r=r_0+r_1p+\ldots + r_{n-1}p^{n-1}$, where $0\le r_i \le p-1$ for each $i$, and let
$$\rr=(r_0,\ldots,r_{n-1}).$$   
Let $\mathbf{\varepsilon}=(\epsilon_0,\ldots,\epsilon_{n-1})$ be an $n$-tuple, such that $\epsilon_i\in\{0,1\}$ for every $i$, then 
$$(\chi \alpha^{\mathbf{\varepsilon}\centerdot(\mathbf{p-1}-\rr)})^{s}=
\chi\alpha^{(\mathbf{1-\varepsilon})\centerdot(\mathbf{p-1-\rr})}.$$
 Moreover, if $\rr=\mathbf{0}$, then we suppose that $\mathbf{\varepsilon}\neq \mathbf{0}$ and
$\mathbf{\varepsilon}\neq \mathbf{1}$, then
$$(\chi\alpha^{\mathbf{\varepsilon}\centerdot(\mathbf{p-1}-\rr)})^{s}\neq
\chi\alpha^{\mathbf{\varepsilon}\centerdot(\mathbf{p-1}-\rr)}$$
where $\mathbf{\varepsilon}\centerdot(\mathbf{p-1}-\rr)=\sum_{i=0}^{n-1} \epsilon_i(p-r_i-1)p^{i}$.
\end{lem}
\begin{proof}Since twisting by $s$ does not affect $\det$ we may assume that
$$\chi(\begin{pmatrix}\lambda & 0 \\ 0 & \mu \end{pmatrix})=\lambda^r \quad \forall \lambda, \mu\in \Fq^{\times}.$$
Then the first part of the lemma amounts to 
$$\mu^{r}(\mu \lambda^{-1})^{\mathbf{\varepsilon}\centerdot(\mathbf{p-1}-\rr)}=\lambda^{r}
 ( \lambda \mu^{-1})^{q-1-r-\mathbf{\varepsilon}\centerdot(\mathbf{p-1}-\rr)}=
\lambda^{r}( \lambda \mu^{-1})^{\mathbf{(1-\varepsilon)}\centerdot(\mathbf{p-1}-\rr)}.$$
 For the second part, we observe that  the equality holds if and only if 
$$\mu^{r+2\mathbf{\varepsilon}\centerdot(\mathbf{p-1}-\rr) }=\lambda^{r+2\mathbf{\varepsilon}\centerdot(\mathbf{p-1}-\rr) }$$
for every $\lambda, \mu \in \Fq^{\times}$. Hence, equality holds if and only if 
$$ \sum_{i=0}^{n-1}(r_i+2(p-1-r_i)\epsilon_i) p^{i}\equiv 0 \pmod{q-1}.$$
Since, $\epsilon_i\in\{0,1\}$  we have 
$$0\le r_i + 2(p-1-r_i)\epsilon_i\le 2(p-1).$$
The congruence implies that $r+2\mathbf{\varepsilon}\centerdot (\mathbf{p-1}-\rr)$ must  take values $0, q-1$ or $2(q-1)$. The extreme values are obtained 
if and only if $\rr=\mathbf{0}$ and $\varepsilon=\mathbf{0}$ or $\rr=\mathbf{0}$ and $\varepsilon=\mathbf{1}$. By our assumptions, both cases are 
excluded. If
 $$r+2\mathbf{\varepsilon}\centerdot (\mathbf{p-1}-\rr)=q-1$$
then we rewrite this as
$$ \sum_{i=0}^{n-1}(p-1-r_i)\epsilon_i p^{i} =\sum_{i=0}^{n-1}(p-1-r_i)(1-\epsilon_i)p^{i}.$$
Hence, for every $i$ we must have 
$$(p-1-r_i)\epsilon_i=(p-1-r_i)(1-\epsilon_i).$$
Since $2\epsilon_i \neq 1$, for every $i$, this forces $r_i=p-1$, for every $i$, but $r<q-1$, hence this case is also excluded.
\end{proof}

\begin{defi}\label{deltasigmar} Suppose that $q=p^{n}$ and let $\rr=(r_0,\ldots, r_{n-1})$ be an $n$-tuple, such that $0\le r_i\le p-1$ for every $i$. 
We define 
$$\delta\in \Sigma_{\rr}$$ 
given by $\delta_i=1$ if $r_i\neq p-1$ and $\delta_i=0$ if $r_i=p-1$. 

We further define $\Sigma_{\rr}'$ to be a subset of $\Sigma_{\rr}$ given by
$$\Sigma_{\rr}'=\Sigma_{\rr}\backslash \{\mathbf{0},\delta\}.$$
\end{defi}
\begin{remar}\label{deltaone}We note that if $p=q$ or $\rr=(p-1,\ldots,p-1)$, then $\Sigma_{\rr}'=\emptyset$ and we always have $b_{\delta}=b_{\mathbf{1}}$.
\end{remar}

\begin{lem}\label{nondeg} Suppose that $q=p^{n}$ and let $\chi:H\rightarrow \Fbar^{\times}$ be a character. Let $r$ be the unique integer, 
such that $0\le r < q-1$ and
$$ \chi(\begin{pmatrix} \lambda & 0 \\0 & \lambda^{-1} \end{pmatrix})= \lambda^{r}\quad \forall \lambda\in \Fq^{\times}.$$
Let $r=r_0+r_1p+\ldots + r_{n-1}p^{n-1}$, where $0\le r_i \le p-1$ for each $i$, and let
$$\rr=(r_0,\ldots,r_{n-1}).$$ 
If $r=0$ then we consider $\inj \rho_{\chi,S}$ and if $r\neq 0$ we consider $\inj \rho_{\chi,\emptyset}$ as representations of $K$ on which $K_1$
acts trivially. 

Suppose that $\mathbf{\varepsilon}\in \Sigma_{\rr}'$. If $r=0$ then we regard $b_{\mathbf{\varepsilon}}$ and 
$b_{\mathbf{1-\varepsilon}}$ as vectors in $(\inj \rho_{\chi, S})^{I_1}$ via the isomorphism of Corollary \ref{sleepy}. If $r\neq 0$ then
we regard $b_{\mathbf{\varepsilon}}$ and $b_{\mathbf{1-\varepsilon}}$ as vectors in $(\inj \rho_{\chi,\emptyset})^{I_1}$ via the isomorphism 
of Corollary \ref{emptyseteasy}.

Then the action of $\HH_K$ on 
$\langle b_{\mathbf{\varepsilon}}, b_{\mathbf{1-\varepsilon}}\rangle_{\Fbar}$ extends to the action of $\HH$, so that
$$\langle b_{\mathbf{\varepsilon}}, b_{\mathbf{1-\varepsilon}}\rangle_{\Fbar}\cong M_{\gamma_{\mathbf{\varepsilon}}}$$
as an $\HH$-module, where 
$$\gamma_{\mathbf{\varepsilon}}=\gamma_{\mathbf{1-\varepsilon}}=\{\chi \alpha^{\mathbf{\varepsilon}\centerdot(\mathbf{p-1}-\rr)},
(\chi \alpha^{\mathbf{\varepsilon}\centerdot(\mathbf{p-1}-\rr)})^{s}\}.$$
\end{lem}

\begin{proof}To ease the notation, let 
$$\psi=\chi \alpha^{\mathbf{\varepsilon}\centerdot(\mathbf{p-1}-\rr)}.$$
We observe that if $b_{\mathbf{1-\varepsilon}}=b_{\mathbf{0}}$, then $\mathbf{\varepsilon}=\delta$ and if $b_{\mathbf{1-\varepsilon}}=b_{\delta}$
then $\mathbf{\varepsilon}=\mathbf{0}$. Since $\mathbf{\varepsilon}\in \Sigma_{\rr}'$ neither of the above can occur.

By Lemma \ref{mrfixitH} and taking into account the twist by a power of $\det$, $I$ acts on $b_{\mathbf{\varepsilon}}$ via the character $\psi$.
By the same argument and Lemma \ref{twistbys} $I$ acts on $b_{\mathbf{1-\varepsilon}}$ via the character $\psi^{s}$. Hence, 
$$b_{\mathbf{\varepsilon}}e_{\psi}=b_{\mathbf{\varepsilon}}\quad\textrm{and}\quad b_{\mathbf{1-\varepsilon}}e_{\psi^{s}}=b_{\mathbf{1-\varepsilon}}.$$
Moreover,  Lemma
\ref{twistbys} says that $\psi\neq \psi^{s}$. The case $\rr=\mathbf{0}$ is not a problem, since
$\mathbf{\varepsilon}\in \Sigma_{\mathbf{0}}'$ implies that $\mathbf{1-\varepsilon}\in \Sigma_{\mathbf{0}}'$. Since $H$ acts 
on  $b_{\mathbf{\varepsilon}}$ and $b_{\mathbf{1-\varepsilon}}$ by different characters, they  are 
linearly independent. Proposition \ref{hardcore} implies that
$$b_{\mathbf{\varepsilon}}T_{n_s}=\sum_{u\in I_1/K_1}u n_{s}^{-1} b_{\mathbf{\varepsilon}}=0\quad
\textrm{and}
\quad b_{\mathbf{1-\varepsilon}}T_{n_s}=\sum_{u\in I_1/K_1}u n_{s}^{-1} b_{\mathbf{1-\varepsilon}}=0.$$
Hence, by Lemma \ref{restML}
$$\langle b_{\mathbf{\varepsilon}}, b_{\mathbf{1-\varepsilon}}\rangle_{\Fbar}\cong
\langle b_{\mathbf{\varepsilon}}\rangle_{\Fbar}\oplus\langle b_{\mathbf{1-\varepsilon}}\rangle_{\Fbar}\cong M_{\psi, \emptyset}\oplus 
M_{\psi^{s},\emptyset}\cong M_{\gamma_{\mathbf{\varepsilon}}}|_{\HH_K}$$
as $\HH_K$-modules. So we define
$$b_{\mathbf{\varepsilon}}T_{\Pi}= b_{\mathbf{1-\varepsilon}}\quad\textrm{and}\quad b_{\mathbf{1-\varepsilon}}T_{\Pi}= b_{\mathbf{\varepsilon}}$$
which gives us the required isomorphism of $\HH$-modules.
\end{proof}  

\begin{prop}\label{injiwahori} Suppose that $q=p^n$ and let $\chi:H\rightarrow \Fbar^{\times}$ be a character, such that $\chi=\chi^{s}$. 
We consider the representation
$$\inj \rho_{\chi,\emptyset}\oplus \inj \rho_{\chi, S}$$
 as a representation of $K$, such that $K_1$ acts trivially. We may extend the action of $\HH_{K}$ on 
$$(\inj \rho_{\chi,\emptyset}\oplus \inj \rho_{\chi, S})^{I_1}$$
to the action of $\HH$, such that $(\inj \rho_{\chi,\emptyset}\oplus \inj \rho_{\chi, S})^{I_1}$ as an $\HH$-module is isomorphic to a 
direct sum of $2^{n-1}$ supersingular modules of $\HH$. 

More precisely, for every $\mathbf{\varepsilon}\in \Sigma_{\mathbf{0}}$ we consider $b_{\mathbf{\varepsilon}}$ as vectors in 
$$(\inj \rho_{\chi,\emptyset}\oplus \inj \rho_{\chi, S})^{I_1}$$ via the isomorphism of Corollary \ref{sleepy}. Then the action of $\HH_{K}$
can be extended to the action of $\HH$ so that
$$\langle b_{\mathbf{0}}, b_{\mathbf{0}}+b_{\mathbf{1}} \rangle_{\Fbar}\cong M_{\gamma}$$
where $\gamma=\{\chi\}$. If $\mathbf{\varepsilon}\in \Sigma_{\mathbf{0}}'$, then
$$\langle b_{\mathbf{\varepsilon}},  b_{\mathbf{1-\varepsilon}}\rangle_{\Fbar}\cong M_{\gamma_{\mathbf{\varepsilon}}}$$
where $\gamma_{\mathbf{\varepsilon}}=\gamma_{\mathbf{1-\varepsilon}}=\{\chi\alpha^{\mathbf{\varepsilon}\centerdot(\mathbf{p-1})},
 \chi(\alpha^{\mathbf{\varepsilon}\centerdot(\mathbf{p-1})})^{s}\}$.
 \end{prop}
\begin{proof} Since, by 
Lemma \ref{mrfixit}  $b_{\mathbf{\varepsilon}}$ for $\mathbf{\varepsilon}\in \Sigma_{\mathbf{0}}$ form a basis of $R_{\mathbf{0}}^{U}$, the second 
part implies the first. Since $\Sigma_{\mathbf{0}}'=\Sigma_{\mathbf{0}}\backslash \{\mathbf{0}, \mathbf{1}\}$, the last part of the Proposition
is given by Lemma \ref{nondeg}.  

Lemmas \ref{genezero} and \ref{gensteinberg} imply that
$$\langle b_{\mathbf{0}} \rangle_{\Fbar}\cong M_{\chi, S},\quad \langle b_{\mathbf{1}}+b_{\mathbf{0}}\rangle_{\Fbar}\cong M_{\chi,\emptyset}$$
as an $\HH_K$-module.  Hence, by Lemma \ref{restML}
$$\langle b_{\mathbf{0}}, b_{\mathbf{1}}+b_{\mathbf{0}}\rangle_{\Fbar}\cong M_{\chi,S}\oplus M_{\chi, \emptyset}\cong M_{\gamma}|_{\HH_{K}}$$
as $\HH_{K}$-modules. Hence, if we define 
$$b_{\mathbf{0}}T_{\Pi}= b_{\mathbf{0}}+b_{\mathbf{1}}\quad \textrm{and}\quad (b_{\mathbf{0}}+b_{\mathbf{1}})T_{\Pi}=b_{\mathbf{0}}$$
we get the required isomorphism. 
\end{proof}

\begin{prop}\label{injregular}Suppose that $q=p^n$, let $\chi:H\rightarrow \Fbar^{\times}$ be a character, such that $\chi\neq\chi^{s}$, and 
 let $a$ and $r$ be  unique integers, such that $1\le a, r\le q-1$ and
$$ \chi(\begin{pmatrix} 1 & 0 \\0 & \lambda \end{pmatrix})= \lambda^{a}\quad \forall \lambda\in \Fq^{\times},\quad 
 \chi(\begin{pmatrix} \lambda & 0 \\0 & \lambda^{-1} \end{pmatrix})= \lambda^{r}\quad \forall \lambda\in \Fq^{\times}.$$
Let $r=r_0+r_1p+\ldots + r_{n-1}p^{n-1}$, where $0\le r_i \le p-1$ for each $i$, and let
$$\rr=(r_0,\ldots,r_{n-1}).$$  
Then
$$\inj \rho_{\chi, \emptyset}\oplus \inj \rho_{\chi^{s}, \emptyset}\cong R_{\rr}\otimes(\det)^{a+r} \oplus R_{\mathbf{p-1}-\rr}\otimes (\det)^{a}$$
where $\mathbf{p-1}-\rr=(p-1-r_0, \ldots, p-1-r_{n-1})$.

We regard the representation $\inj \rho_{\chi, \emptyset}\oplus \inj \rho_{\chi^{s}, \emptyset}$ as a representation of $K$, on which $K_1$ acts 
trivially. Let $c$ and $d$ be the cardinality of the  sets:
$$c=|\{r_i: r_i\neq p-1\}|\quad\textrm{and}\quad d=|\{r_i: r_i \neq 0\}|$$
then we may extend the action of $\HH_{K}$ on 
$$(\inj \rho_{\chi,\emptyset}\oplus \inj \rho_{\chi^{s}, \emptyset})^{I_1}$$
to the action of $\HH$, such that $(\inj \rho_{\chi,\emptyset}\oplus \inj \rho_{\chi^{s}, \emptyset})^{I_1}$ as an $\HH$-module is isomorphic to a 
direct sum of $L_{\gamma}$ and $2^{c-1}+2^{d-1}-2$ supersingular modules of $\HH$. 

More precisely, let $b_{\mathbf{\varepsilon}}$, for $\mathbf{\varepsilon}\in \Sigma_{\rr}$, be a basis of $(\inj \rho_{\chi,\emptyset})^{I_1}$
and let $\bar{b}_{\mathbf{\varepsilon}}$, for $\mathbf{\varepsilon}\in \Sigma_{\mathbf{p-1}-\rr}$, be a basis of 
$(\inj \rho_{\chi^{s},\emptyset})^{I_1}$ via the isomorphism above. Then the action of $\HH_{K}$
can be extended to the action of $\HH$ so that
$$\langle b_{\mathbf{0}}, b_{\mathbf{1}}, \bar{b}_{\mathbf{0}}, \bar{b}_{\mathbf{1}} \rangle_{\Fbar}\cong L_{\gamma}$$
and 
$$\langle b_{\mathbf{0}}, \bar{b}_{\mathbf{0}}\rangle_{\Fbar}\cong M_{\gamma}$$
where $\gamma=\{\chi, \chi^{s}\}$. If $\mathbf{\varepsilon}\in \Sigma_{\rr}'$, then
$$\langle b_{\mathbf{\varepsilon}},  b_{\mathbf{1-\varepsilon}}\rangle_{\Fbar}\cong M_{\gamma_{\mathbf{\varepsilon}}}$$
where $\gamma_{\mathbf{\varepsilon}}=\gamma_{\mathbf{1-\varepsilon}}=\{\chi\alpha^{\mathbf{\varepsilon}\centerdot(\mathbf{p-1}-\rr)},
 (\chi\alpha^{\mathbf{\varepsilon}\centerdot(\mathbf{p-1}-\rr)})^{s}\}$. If $\mathbf{\varepsilon}\in \Sigma_{\mathbf{p-1}-\rr}'$ then
$$\langle \bar{b}_{\mathbf{\varepsilon}}, \bar{b}_{\mathbf{1-\varepsilon}}\rangle_{\Fbar}\cong M_{\bar{\gamma}_{\mathbf{\varepsilon}}}$$
where $\bar{\gamma}_{\mathbf{\varepsilon}}=\bar{\gamma}_{\mathbf{1-\varepsilon}}=\{\chi^{s}\alpha^{\mathbf{\varepsilon}\centerdot \rr},
 (\chi^{s}\alpha^{\mathbf{\varepsilon}\centerdot \rr })^{s}\}$.
\end{prop}
\begin{proof} The first part of the Proposition follows from Corollary \ref{emptyseteasy} and Corollary \ref{bar}. For the second part we observe 
that since $\chi\neq \chi^{s}$, we have $r\neq q-1$ and hence vectors $b_{\mathbf{0}}$, $b_{\mathbf{1}}$, $\bar{b}_{\mathbf{0}}$ and 
$\bar{b}_{\mathbf{1}}$ are linearly independent. Lemma \ref{genezero} implies that 
$$\langle b_{\mathbf{0}}\rangle_{\Fbar}\cong M_{\chi, \emptyset}\quad \textrm{and}\quad \langle \overline{b}_{\mathbf{0}}\rangle_{\Fbar}\cong 
M_{\chi^{s},\emptyset} $$
as $\HH_{K}$-modules. Lemma \ref{mrfixitH} with the appropriate twist by a power of $\det$ says that $H$ acts on $b_{\mathbf{1}}$ by a character
$\chi \alpha^{\mathbf{1}\centerdot(\mathbf{p-1}-\rr)}$ and $H$ acts on $\overline{b}_{\mathbf{1}}$ by a character 
$\chi^{s} \alpha^{\mathbf{1}\centerdot \rr }$. Lemma \ref{twistbys} implies that
$$\chi \alpha^{\mathbf{1}\centerdot(\mathbf{p-1}-\rr)}=\chi^{s}\quad
\textrm{and}\quad \chi^{s} \alpha^{\mathbf{1}\centerdot \rr }=\chi.$$ 
Hence,
$$b_{\mathbf{1}}e_{\chi^{s}}=b_{\mathbf{1}}\quad\textrm{and}\quad \bar{b}_{\mathbf{1}}e_{\chi}=\bar{b}_{\mathbf{1}}.$$
Proposition \ref{hardcore} implies that 
$$(-1)^{r+1}b_{\mathbf{1}}T_{n_{s}}=(-1)^{r+1}\sum_{u\in I_1/ K_1} u n_{s}^{-1}b_{\mathbf{1}}=b_{\mathbf{0}}$$
and 
$$(-1)^{q-r}\bar{b}_{\mathbf{1}}T_{n_{s}}=(-1)^{q-r}\sum_{u\in I_1/ K_1} u n_{s}^{-1}\bar{b}_{\mathbf{1}}=\bar{b}_{\mathbf{0}}.$$
Hence, by Lemma \ref{restML} 
$$\langle b_{\mathbf{0}}, b_{\mathbf{1}}, \bar{b}_{\mathbf{0}}, \bar{b}_{\mathbf{1}} \rangle_{\Fbar}\cong L_{\gamma}\mid_{\HH_K}$$
as $\HH_K$-modules. We note that if $p=2$ then $1=-1$ and if $p\neq 2$ then $(-1)^{q-r}=(-1)^{r+1}$. So if we define
$$b_{\mathbf{1}}T_{\Pi}= \bar{b}_{\mathbf{1}},\quad \bar{b}_{\mathbf{1}}T_{\Pi}= b_{\mathbf{1}}, 
\quad b_{\mathbf{0}}T_{\Pi}=\bar{b}_{\mathbf{0}},\quad \bar{b}_{\mathbf{0}}T_{\Pi}= b_{\mathbf{0}}$$
we get the required isomorphism of $\HH$-modules. Moreover,
$$\langle b_{\mathbf{0}}, \bar{b}_{\mathbf{0}}\rangle_{\Fbar}\cong M_{\gamma}$$
as $\HH$-module. The last part of the Proposition follows from Lemma \ref{nondeg}. Since $\dim (\inj \rho_{\chi,\emptyset})^{I_1}=2^{c}$ and
 $\dim (\inj \rho_{\chi^{s},\emptyset})^{I_1}=2^{d}$ an easy calculation gives us the number of indecomposable summands.
\end{proof}
\begin{remar}If $p=q$, then $\Sigma'_{\rr}=\emptyset$ and Propositions \ref{injiwahori} and \ref{injregular} specialise  to Proposition \ref{injqp}.
\end{remar}
The following Proposition can be seen as a consolation for the Remark \ref{hooray}.
\begin{prop}\label{consol} Suppose that $q=p^{n}$, $\chi\neq \chi^{s}$ and let $\rho$ be a representation of $\Gamma$, such that
$\rho^{U}\cong M_{\chi,\emptyset}\oplus M_{\chi^{s},\emptyset}$ as an $\HH_{\Gamma}$-module, and $\rho=\langle\Gamma \rho^{U}\rangle_{\Fbar}$, then
$$\rho\cong \rho_{\chi,\emptyset}\oplus \rho_{\chi^{s},\emptyset}.$$
\end{prop}
\begin{proof} If $\rho$ is a semi-simple representation of $\Gamma$, then Corollary \ref{hgammamod} implies the Lemma. Suppose that $\rho$ is not
semi-simple. Let $\soc(\rho)$ be the maximal semi-simple subrepresentation of $\rho$. Since $\rho$ is generated by $\rho^{U}$ as a 
$\Gamma$-representation, the space $(\soc(\rho))^{U}$ is one dimensional, and hence $\soc(\rho)$ is an irreducible representation of $\Gamma$. By Corollary
\ref{hgammamod} and symmetry we may assume that 
$$\soc(\rho)\cong \rho_{\chi,\emptyset}.$$
Since, $\soc(\rho)$ is irreducible, $\rho$ is an essential extension of $\rho_{\chi,\emptyset}$. By this we mean that every non-zero $\Gamma$ invariant 
subspace of $\rho$ intersects $\rho_{\chi,\emptyset}$ non-trivially. This implies that there exists an exact sequence 
\begin{displaymath}
\xymatrix{0\ar[r]&\rho\ar[r]& \inj \rho_{\chi, \emptyset}}
\end{displaymath}
of $\Gamma$-representations. After twisting by a power of determinant we may assume that $\chi$ is given by
$\chi(\begin{pmatrix}\lambda & 0 \\ 0 & \mu \end{pmatrix})=\lambda^{r}$,
where $0< r <q-1$. The inequalities are strict, since $\chi\neq \chi^{s}$. Let $\rr$ be the corresponding $n$-tuple. Let $\varepsilon\in \Sigma_{\rr}$
and $b_{\varepsilon}\in (\inj \rho_{\chi,\emptyset})^U$, then $H$ acts on $b_{\varepsilon}$ by the character 
$\chi \alpha^{\varepsilon \centerdot(\mathbf{p-1}-\rr)}$. In particular, if $\varepsilon'\in \Sigma_{\rr}$, such that $\varepsilon'\neq \varepsilon$, 
then $H$ acts on $b_{\varepsilon}$ and $b_{\varepsilon'}$ by different characters. As a consequence of this, the submodule 
$M_{\chi^{s}, \emptyset}$ of $\rho^U$ must be mapped to 
some subspace $\langle b_{\mathbf{\varepsilon}}\rangle_{\Fbar}$ of $(\inj \rho_{\chi,\emptyset})^{U}$, where $\varepsilon\in \Sigma_{\rr}$. By examining 
the action of $H$, we get that 
$\chi^{s}=\chi \alpha^{\mathbf{\varepsilon}\centerdot(\mathbf{p-1}-\rr)}.$
This implies that 
$$\varepsilon\centerdot (\mathbf{p-1}-\rr)+r\equiv 0 \pmod{q-1}.$$
Since $0<r<q-1$ and $\varepsilon\in \Sigma_{\rr}$, we have 
$$0<\varepsilon\centerdot (\mathbf{p-1}-\rr)+r\le q-1.$$
Hence, we get an equality on the right hand side, which implies that, for each $i$, $(1-\epsilon_i)(p-1-r_i)=0$. 
So $\varepsilon=\delta$, and $b_{\epsilon}=b_{\mathbf{1}}$, see 
\ref{deltasigmar} and \ref{deltaone}. However, by Proposition \ref{hardcore} (ii)
$$b_{\mathbf{1}}T_{n_s}=(-1)^{r+1}b_{\mathbf{0}}\neq 0.$$  
We obtain a contradiction, since $T_{n_s}$ kills $M_{\chi^{s},\emptyset}$. 
\end{proof}         
\section{Coefficient systems}\label{coeffsys}
We closely follow \cite{ss1} and \cite{ss2}\S V, where the $G$-equivariant coefficient systems of $\CCC$-vector spaces are treated. In fact, the results of
this Section do not depend on the underlying field. Our motivation to use coefficient systems stems from \cite{rs}, where the equivariant coefficient
systems of $\Fbar$-vector spaces of finite Chevalley groups are considered. 
\subsection{Definitions}
 The Bruhat - Tits tree $X$ of $G$ is the simplicial complex, whose vertices are the similarity classes $[L]$ of
 $\oF$-lattices  in a $2$-dimensional $F$-vector space $V$ and whose edges are $1$-simplices, given by families
 $\{[L_{0}], [L_{1}]\}$ of 
similarity classes such that
$$ \pif L_{0} \subset L_{1} \subset L_{0}.$$ 
We denote by $X_0$ the set of all vertices and by $X_1$ the set of all edges. 
\begin{defi} Let $\sigma$ be a simplex in $X$, then  we define
$$\KK(\sigma)=\{g\in \Aut_{F}(V): g\sigma=\sigma\}.$$
\end{defi}
By fixing a basis $\{v_1, v_2\}$ of $V$ we identify $G$ with $\Aut_F(V)$. Let 
$$\sigma_{0}=[\oF v_1+\oF v_2]\quad \textrm{and} \quad \sigma_{1}=\{[\oF v_1+\oF v_2], [\oF v_1+\pF v_2]\}.$$
Then $\sigma_{0}$ is a vertex  and $\KK(\sigma_{0})=F^{\times}K$, and $\sigma_{1}$ is an edge containing
 a vertex $\sigma_{0}$. Moreover,  $\KK(\sigma_{1})$ is the group generated by $I$ and $\Pi$.  

\begin{defi} A coefficient system $\VV$ (of $\Fbar$-vector spaces) on $X$ consists of 
\begin{itemize} 
\item[-] $\Fbar$ vector spaces $V_{\sigma}$ for each simplex $\sigma$ of $X$, and
\item[-] linear maps $r^{\sigma'}_{\sigma}:V_{\sigma'}\rightarrow V_{\sigma}$ for each pair $\sigma\subseteq\sigma'$
 of simplices of $X$ such that for every simplex $\sigma$, $r^{\sigma}_{\sigma}=\id_{V_{\sigma}}$.
\end{itemize}
\end{defi} 

\begin{defi} We say the group $G$ acts on the coefficient system $\VV$, if for every $g\in G$ and for every simplex $\sigma$ there is given
a linear map
$$g_{\sigma}: V_{\sigma}\rightarrow V_{g\sigma},$$
such that
\begin{itemize} 
\item[-] $g_{h\sigma}\circ h_{\sigma}=(gh)_{\sigma}$, for every $g$, $h\in G$ and for every simplex $\sigma$,
\item[-] $1_{\sigma}=\id_{V_{\sigma}}$ for every simplex $\sigma$,
\item[-] the following diagram commutes for every $g\in G$ and every pair of simplices $\sigma\subseteq \sigma'$:
\begin{displaymath}
\xymatrix {V_{\sigma}\ar[r]^{g_{\sigma}} & V_{g\sigma} \\
          V_{\sigma'} \ar[u]^{r^{\sigma'}_{\sigma}} \ar[r]^{g_{\sigma'}} & V_{g\sigma'}\ar[u]_{r^{g\sigma'}_{g\sigma}}}.
\end{displaymath}  
 \end{itemize}
\end{defi}
In particular, the stabiliser $\KK(\sigma)$ acts linearly on $V_{\sigma}$ for any simplex $\sigma$.
 \begin{defi} A $G$-equivariant coefficient system $(V_{\tau})_{\tau}$ on $X$ is a coefficient system on $X$ together with a $G$-action , such that 
the action of the stabiliser $\KK(\sigma)$ of a simplex $\sigma$ on $V_{\sigma}$ is smooth. 
\end{defi}
\begin{remar} The definition given in \cite{ss2} \S V, requires the action to factor through a discrete quotient.
\end{remar} 
Let $\coeff$ denote the category of all equivariant coefficient systems on $X$, equipped with the obvious morphisms.

The following observation will turn out to be very useful. Suppose that $G$ acts on a coefficient system $\VV=(V_{\sigma})_{\sigma}$.
Let $\tau'$ be an edge containing a vertex $\tau$. There exists $g\in G$, such that $\tau'=g \sigma_1$ and $\tau=g \sigma_0$. Then
$$V_{\tau}=g_{\sigma_0} V_{\sigma_0},\quad V_{\tau'}= g_{\sigma_1} V_{\sigma_1}$$
and 
$$ r^{\tau'}_{\tau}= g_{\sigma_0}\circ r^{\sigma_1}_{\sigma_0} \circ (g^{-1})_{\tau'}.$$ 
\subsection{Homology}
Let $X_{(0)}$ be the set of vertices on the tree and let $X_{(1)}$ be the set of oriented edges on the tree. We will say that two vertices
 $\sigma$ and $\sigma'$ are neighbours if $\{\sigma, \sigma'\}$ is an edge. And we will write
$$(\sigma, \sigma')$$
to mean a directed edge going from $\sigma$ to $\sigma'$. Let $\VV=(V_{\tau})_{\tau}$ be an equivariant coefficient system. We define a
space of oriented $0$-chains to be 
$$C^{or}_c(X_{(0)}, \VV)=\Fbar\mbox{-vector space of all maps }\omega: X_{(0)}\rightarrow \mathop{\/\dot{\bigcup}}\limits_{\sigma\in X_{0}}V_{\sigma}$$
such that 
\begin{itemize}
\item[-] $\omega$ has finite support and
\item[-] $\omega(\sigma)\in V_{\sigma}$ for every vertex $\sigma$.
\end{itemize}
Similarly, the space of oriented $1$-chains is
$$C^{or}_c(X_{(1)}, \VV)=\Fbar\mbox{-vector space of all maps }\omega: X_{(1)}\rightarrow
 \mathop{\/\dot{\bigcup}}\limits_{\{\sigma, \sigma'\}\in X_{1}}V_{\{\sigma,\sigma'\}}$$
such that
\begin{itemize}
\item[-] $\omega$ has finite support,
\item[-] $\omega((\sigma,\sigma'))\in V_{\{\sigma, \sigma'\}}$ ,
\item[-] $\omega((\sigma', \sigma))=-\omega((\sigma, \sigma'))$ for every oriented edge $(\sigma,\sigma')$.
\end{itemize}
The group $G$ acts on $C^{or}_c(X_{(0)}, \VV)$ via
$$(g\omega)(\sigma)=g_{g^{-1}\sigma}(\omega(g^{-1}\sigma))$$
and on $C^{or}_c(X_{(1)}, \VV)$ via
$$(g\omega)((\sigma, \sigma'))= g_{\{g^{-1}\sigma, g^{-1}\sigma'\}}(\omega((g^{-1}\sigma, g^{-1}\sigma'))).$$
The action on both spaces is smooth.

The boundary map is given by
$$\partial:C^{or}_c(X_{(1)}, \VV)\rightarrow C^{or}_c(X_{(0)}, \VV)$$
$$\omega\mapsto (\sigma\mapsto \sum_{\sigma'}r^{\{\sigma, \sigma'\}}_{\sigma}(\omega((\sigma, \sigma'))))$$
where the sum is taken over all the neighbours of $\sigma$. The map $\partial$ is $G$-equivariant. 

We define $H_0(X,\VV)$ to be the cokernel of $\partial$. It is naturally a smooth representation of $G$.
\subsection{Some computations of $H_0(X,\VV)$}
Throughout this section we fix an equivariant coefficient system $\VV=(V_{\tau})_{\tau}$, with the restriction maps given by $r^{\tau'}_{\tau}$.
 Our first lemma follows immediately from the definition of $\partial$.
\begin{lem}\label{delta} Let $\omega$ be an oriented $1$-chain supported on a single edge $\tau=\{\sigma,\sigma'\}$. Let
 $$v=\omega((\sigma, \sigma')).$$
Then $$\partial(\omega)=\omega_{\sigma}-\omega_{\sigma'},$$
where  $\omega_{\sigma}$ and $\omega_{\sigma'}$ are $0$-chains supported only on  $\sigma$ and $\sigma'$ respectively. Moreover,
 $$\omega_{\sigma}(\sigma)=r^{\tau}_{\sigma}(v)$$ 
and
$$\omega_{\sigma'}(\sigma')=r^{\tau}_{\sigma'}(v).$$

\end{lem} 
\begin{lem}\label{inj} Let $\omega$ be a $0$-chain supported on a single vertex $\sigma$. Suppose that the restriction map $r^{\sigma_1}_{\sigma_0}$
is an injection, then the image of $\omega$ in $H_0(X,\VV)$ is non-zero.
\end{lem}
\begin{proof} Since every restriction map is conjugate to $r^{\sigma_1}_{\sigma_0}$ by some element of $G$, it follows that every restriction
map is injective. 

Let $\omega'$ be a non-zero oriented $1$-chain. We may think of the support of $\omega'$ as the union of edges of a
 finite subgraph $\mathcal T$ of $X$. Since 
all the restriction maps are injective, Lemma \ref{delta} implies that $\partial(\omega')$ will not vanish on the vertices of $\mathcal T$, which have only
 one neighbour in $\mathcal T$. In particular, $\partial(\omega')$ will be supported on at least $2$ vertices. Hence, 
$\omega\not\in \partial C^{or}_c(X_{(1)},\VV).$
\end{proof}
 
\begin{lem}\label{surj}Let $\omega$ be $0$-chain. Suppose that the restriction map $r^{\sigma_1}_{\sigma_0}$ is
 surjective, then there exists a $0$-chain
$\omega_0$, supported on a single vertex $\sigma_0$, such that
$$\omega+ \partial C^{or}_c(X_{(1)},\VV)=\omega_0+\partial C^{or}_c(X_{(1)},\VV).$$
\end{lem}
\begin{proof} Since every restriction map is conjugate to $r^{\sigma_1}_{\sigma_0}$ by some element of $G$, it follows that every restriction map
is surjective. 

It is enough to prove the statement when $\omega$ is supported on a single vertex $\tau$, since an arbitrary $0$-chain is a finite sum of such. If 
$\tau=\sigma_0$ then we are done. Otherwise, there exists a directed path going from $\sigma_0$ to $\tau$, consisting of finitely many directed edges
$(\sigma_0, \tau_1), \ldots, (\tau_m, \tau).$

We argue by induction on $m$. Let $v=\omega(\tau)$. Since $r^{\{\tau_m, \tau\}}_{\tau}$ is surjective there exists $v'\in V_{\{\tau_m, \tau\}}$,
such that 
$$r^{\{\tau_m, \tau\}}_{\tau}(v')=v.$$
Let $\omega'$ be an oriented $1$-chain supported on the single edge $\{\tau_m, \tau\}$ with $\omega'((\tau_m, \tau))=v'$. By Lemma \ref{delta}
$\omega+\partial(\omega')$ is supported on a single vertex $\tau_m$. Since, the number of edges in the directed path has decreased by one, the claim
 follows from induction. 
\end{proof}   

The following special case will be used in the calculations of modules of  the Hecke algebra.
\begin{lem}\label{act} Let $\omega_0$ be a $0$-chain supported on a single vertex $\sigma_0$. Let
 $$v_0=\omega_0(\sigma_0)$$
 and suppose that there exists $v_1\in V_{\sigma_1}$, such that
 $$r^{\sigma_1}_{\sigma_0}(v_1)=v_0.$$
Let $\omega'$ be a $0$-chain supported on a single vertex $\sigma_0$ with 
$$\omega'(\sigma_0)= r^{\sigma_1}_{\sigma_0}((\Pi^{-1})_{\sigma_1}(v_1)),$$ then
$$\Pi^{-1} \omega_0 + \partial C^{or}_c(X_{(1)},\VV)=\omega'+\partial C^{or}_c(X_{(1)},\VV).$$
\end{lem}
\begin{proof} We observe that $\Pi \sigma_0=\Pi^{-1}\sigma_0$ and  $\sigma_1=\{\sigma_0,\Pi\sigma_0\}$. The $0$-chain $\Pi^{-1}\omega_0$ is supported on
 a single vertex $\Pi\sigma_0$ with
 $$(\Pi^{-1}\omega_0)(\Pi\sigma_0)=(\Pi^{-1})_{\sigma_0}(v_0).$$   
Let $\omega_1$ be an oriented $1$-chain supported on a single edge $\sigma_1$ with
 $$\omega_1((\sigma_0, \Pi \sigma_0))= (\Pi^{-1})_{\sigma_1}(v_1).$$
 From Lemma \ref{delta} we know that $\partial(\omega_1)$ is supported only on $\sigma_0$ and $\Pi \sigma_0$. Moreover,
\begin{equation}
\begin{split}
\partial(\omega_1)(\Pi \sigma_0)&=r^{\sigma_1}_{\Pi\sigma_0}(\omega_1((\Pi \sigma_0, \sigma_0))=r^{\sigma_1}_{\Pi\sigma_0}(-(\Pi^{-1})_{\sigma_1}(v_1))\\
                        &=-(r^{\sigma_1}_{\Pi\sigma_0}\circ (\Pi^{-1})_{\sigma_1})(v_1)=-((\Pi^{-1})_{\sigma_0}\circ r^{\sigma_1}_{\sigma_0})(v_1)=
 -(\Pi^{-1})_{\sigma_0}(v_0),
\end{split}\notag
\end{equation}
and 
$$\partial(\omega_1)(\sigma_0)=r^{\sigma_1}_{\sigma_0}(\omega_1((\sigma_0, \Pi\sigma_0)))= r^{\sigma_1}_{\sigma_0}((\Pi^{-1})_{\sigma_1}(v_1)).$$
Hence
$$\partial(\omega_1)= \omega'- \Pi^{-1}\omega_0$$ 
and that establishes the claim.
\end{proof}

\begin{prop}\label{isores}Suppose that the restriction map $r^{\sigma_1}_{\sigma_0}$ is an isomorphism of vector spaces. Then
$$H_0(X,\VV)|_{\kzr}\cong V_{\sigma_0}$$
and
$$H_0(X,\VV)|_{\kon}\cong V_{\sigma_1}.$$
Moreover, the diagram 
\begin {displaymath}
\xymatrix {V_{\sigma_0} \ar[r]^{\cong} & H_0(X,\VV) \\
          V_{\sigma_1} \ar[u]^{r^{\sigma_1}_{\sigma_0}} \ar[r]^{\cong} & H_0(X,\VV) \ar[u]_{\id}}
\end {displaymath}
of $F^{\times}I$-representations commutes.
\end{prop}
\begin{proof}Let $C^{or}_{c}(\sigma_0,\VV)$ be a subspace of $C^{or}_c(X_{(0)},\VV)$ consisting of the $0$-chains whose support lies
 in the simplex $\sigma_0$, with the understanding that the $0$-chain which vanishes on every simplex is supported on the empty simplex.
Let $\jmath$ be the composition
$$\jmath:C^{or}_{c}(\sigma_0,\VV)\hookrightarrow C^{or}_c(X_{(0)}, \VV)\rightarrow H_0(X, \VV).$$
Then $\jmath$ is $\kzr$ equivariant. Moreover, Lemma \ref{inj} says that $\jmath$ is an injection and Lemma \ref{surj} says that it is a surjection.
Hence 
$$\jmath: C^{or}_c(\sigma_0,\VV)\cong H_0(X, \VV)|_{\kzr}.$$
Let $\ev_0$ be the map
$$\ev_0:C^{or}_c(\sigma_0,\VV)\rightarrow V_{\sigma_0}$$
$$\omega\mapsto \omega(\sigma_0)$$
then $\ev_0$ is an isomorphism of $\kzr$-representations. Hence
$$ \jmath\circ(\ev_{0})^{-1}: V_{\sigma_0}\cong\ H_0(X,\VV)|_{\kzr}.$$
Since $\VV$ is $G$-equivariant, the map $r^{\sigma_1}_{\sigma_0}$ is $F^{\times}I=\kon\cap\kzr$-equivariant and  since it is isomorphism of vector spaces,
 we obtain that
$$\jmath\circ(\ev_0)^{-1}\circ r^{\sigma_1}_{\sigma_0}:V_{\sigma_1}|_{F^{\times}I}\cong H_0(X,\VV)|_{F^{\times}I}.$$
We claim that this isomorphism is in fact $\kon$-equivariant. Let $v_1\in V_{\sigma_1}$, let $v_0=r^{\sigma_1}_{\sigma_0}(v_1)$ and let 
$\omega_0\in C^{or}_{c}(\sigma_0,\VV)$, such that $\omega_{0}(\sigma_0)=v_0$. Then
$$(\jmath\circ(\ev_0)^{-1}\circ r^{\sigma_1}_{\sigma_0})(v_1)=\omega_0+\partial C^{or}_{c}(X_{(1)}, \VV).$$
By Lemma \ref{act}
$$ \Pi^{-1}\omega_0 +\partial C^{or}_{c}(X_{(1)},\VV)= \omega'+ \partial C^{or}_{c}(X_{(1)}, \VV),$$
where $\omega'\in C^{or}_{c}(\sigma_0,\VV)$ with $\omega'(\sigma_0)= r^{\sigma_1}_{\sigma_0}((\Pi^{-1})_{\sigma_1}(v_1))$. This implies that
$$\Pi^{-1}(\jmath\circ(\ev_0)^{-1}\circ r^{\sigma_1}_{\sigma_0})(v_1)=
(\jmath\circ(\ev_0)^{-1}\circ r^{\sigma_1}_{\sigma_0})((\Pi^{-1})_{\sigma_1}(v_1)).$$
Since $\Pi^{-1}$ and $F^{\times}I$ generate $\kon$ this proves the claim.

The commutativity of the diagram follows from the way the isomorphisms are constructed.
\end{proof}  

\subsection{Constant functor}
The content of this Section is essentially \cite{rs} Lemma 1.1 and  Theorem 1.2. Let $\Rep_G$ be the category of smooth $\Fbar$-representations of $G$.
Let $\pi$ be a smooth representation of $G$ with the underlying vector space $\WW$. Let $\sigma$ be a simplex on the tree $X$, we set
$$(\const_{\pi})_{\sigma}=\WW.$$
If $\sigma$ and $\sigma'$ are two simplices, such that $\sigma\subseteq \sigma'$ then we define the restriction map 
$$r^{\sigma'}_{\sigma}=\id_{\WW}.$$
For every $g\in G$ and every simplex $\sigma$ in $X$ we define a linear map $g_{\sigma}$ by
$$g_{\sigma}:(\const_{\pi})_{\sigma}\rightarrow (\const_{\pi})_{g\sigma},\quad v\mapsto \pi(g)v.$$
This gives a $G$-equivariant coefficient system on $X$, which we denote by $\const_{\pi}$.
\begin{defi} We define the constant functor
$$\const: \Rep_{G}\rightarrow\coeff,\quad \pi\mapsto \const_{\pi}.$$
\end{defi}

\begin{lem}\label{Hconst} Let $\pi$ be a smooth representation of $G$, then
$$H_0(X,\const_{\pi})\cong \pi$$
as a $G$-representation. 
\end{lem}
\begin{proof}We have an evaluation map
$$\ev: C^{or}_{c}(X_{(0)},\const_{\pi})\rightarrow \pi, \quad \omega \mapsto \sum_{\sigma\in X_{(0)}}\omega(\sigma).$$  
Since the restriction maps are just $\id_{\WW}$, Lemma \ref{delta} implies that the image of the boundary map
$\partial C^{or}_{c}(X_{(1)},\const_{\pi})$ is contained in the kernel of $\ev$. Hence, we get a $G$-equivariant map
$$H_0(X,\const_{\pi})\rightarrow \pi.$$
It is enough to show that this is an isomorphism of vector spaces, and this is implied by Proposition \ref{isores}.
\end{proof}

\begin{prop}Let $\VV=(V_{\sigma})_{\sigma}$ be a $G$-equivariant coefficient system with the restriction maps $r^{\sigma'}_{\sigma}$
and let $(\pi,\WW)$ be a smooth representation of $G$, then
$$\Hom_{\coeff}(\VV,\const_{\pi})\cong \Hom_{G}(H_0(X,\VV), \pi).$$
\end{prop}
\begin{proof}Any morphism of $G$-equivariant coefficient systems will induce a $G$-equivariant homomorphism in the $0$-th homology. Hence by
Lemma \ref{Hconst}  we have a map 
$$\Hom_{\coeff}(\VV,\const_{\pi})\rightarrow \Hom_{G}(H_0(X,\VV), \pi).$$
We will construct an inverse of this. Let $\phi\in \Hom_{G}(H_0(X,\VV), \pi)$, let $\sigma$ be a vertex on the tree $X$, let
$v$ be a vector in $V_{\sigma}$, and let $\omega_{\sigma,v}$ be a $0$-chain, such that 
$$ \supp \omega_{\sigma,v}\subseteq \sigma,\quad \omega_{\sigma,v}(\sigma)=v,$$
 then we define
$$\phi_{\sigma}:V_{\sigma}\rightarrow \WW,\quad v\mapsto \phi(\omega_{\sigma,v}+\partial C^{or}_{c}(X_{(1)},\VV)).$$ 
Let $\tau$ be an edge in $X$ with vertices $\sigma$ and $\sigma'$, we define 
$$\phi_{\tau}: V_{\tau}\rightarrow \WW,\quad v\mapsto \phi_{\sigma}(r^{\tau}_{\sigma}(v)).$$
Lemma \ref{delta} implies that the definition of $\phi_{\tau}$ does not depend on the choice of vertex. Hence, the collection of linear maps
$(\phi_{\sigma})_{\sigma}$ is a morphism of coefficient systems, which induces $\phi$ on the $0$-th homology. An easy check shows that 
$(\phi_{\sigma})_{\sigma}$ respect the $G$-action on $\VV$ and $\const_{\pi}$.
\end{proof}   
 
\subsection{Diagrams}\label{diag}

\begin{defi}\label{diagrams}
Let  $\diag$ be the category, whose objects are diagrams 
\begin {displaymath}
\xymatrix {D_{0} \\
          D_{1} \ar[u]^{r}}
\end {displaymath}
where $(\rho_0, D_{0})$ is a a smooth  $\Fbar$-representation of $\kzr$, $(\rho_1, D_{1})$ is a smooth $\Fbar$-representation of $\kon$,
 and
 $r\in \Hom_{F^{\times}I}(D_{1}, D_{0})$.

The morphisms between two objects $(D_0, D_1,r)$ and $(D'_{0}, D'_{1}, r')$ are pairs $(\psi_{0}, \psi_{1})$, such that
$\psi_0\in \Hom_{\kzr}(D_0, D'_{0})$,  $\psi_1\in \Hom_{\kon}(D_1, D'_{1})$ and the diagram:
\begin {displaymath}
\xymatrix {D_{0} \ar[r]^{\psi_{0}} & D'_{0} \\
          D_{1} \ar[u]^{r} \ar[r]^{\psi_{1}} & D'_{1} \ar[u]^{r'}}
\end {displaymath}
of $F^{\times}I$ representations commutes.
\end{defi}
The main result of this section is Theorem \ref{equiv}, which says that the categories $\diag$ and $\coeff$ are equivalent. It is easier to work with 
objects of $\diag$ than the coefficient systems.
  
\begin{defi} Let $\VV=(V_{\sigma})_{\sigma}$ be an object in $\coeff$. Let $ \DD: \coeff \rightarrow \diag$ be a functor, given by
 \begin {displaymath}
\VV \mapsto\xymatrix {V_{\sigma_{0}} \\
          V_{\sigma_{1}} \ar[u]^{r^{\sigma_{1}}_{\sigma_{0}}}}.
\end {displaymath}
\end{defi}

We will construct a functor $\CC:\diag\rightarrow \coeff$ and show that the functors $\CC$ and $\DD$ induce an equivalence of categories.

\subsubsection{Underlying vector spaces}
 Let $D=(D_0, D_1, r)$ be an object in $\diag$. Let $i\in \{0,1\}$, we define 
$\cIndu{\KK(\sigma_i)}{G}{\rho_i}$, to be a  representation of $G$  whose underlying vector space consists of functions
$$f:G\rightarrow D_i$$
such that
$$f(kg)= \rho_i(k)f(g) \quad \forall g\in G, \quad  \forall k\in \KK(\sigma_i)$$
and $\supp f $ is compact modulo the centre. The group $G$ acts by the right translations, that is 
$$(gf)(g_1)= f(g_1 g).$$ 
Let $\tau$ be a vertex on the tree $X$, then there exists $g\in G$, such that $\tau= g \sigma_0$. Let 
$$\FF_{\tau}=\{f\in \cIndu{\kzr}{G}{\rho_{0}}: \supp f \subseteq \kzr g^{-1}\}.$$
The space $\FF_{\tau}$ is independent of the choice of $g$. Let $\tau'$ be an edge on the tree $X$, then there exists $g\in G$ such that 
$\tau= g \sigma_1$. We define  
 $$\FF_{\tau'}=\{f\in \cIndu{\kon}{G}{\rho_{1}}: \supp f \subseteq \kon g^{-1}\}.$$
We observe that $\FF_{\tau'}$ is also independent of the choice of $g$.
\subsubsection{Restriction maps}\label{restmap}
Let $i\in\{0,1\}$, then $\FF_{\sigma_i}$ is naturally isomorphic to $D_i$ as a $\KK(\sigma_i)$ representation. The isomorphism is given by
$$\ev_{i}:\FF_{\sigma_i}\rightarrow D_i,\quad f \mapsto f(1).$$
The inverse is given by
$$\ev_{i}^{-1}: D_i \rightarrow \FF_{\sigma_i},\quad v\mapsto f_v$$
where $f_v(k)=\rho_i(k)v$, if $k\in \KK(\sigma_i)$, and $0$ otherwise. Let
$$r^{\sigma_1}_{\sigma_0}=\ev_{0}^{-1}\circ r \circ \ev_{1}.$$
Then $r^{\sigma_1}_{\sigma_0}$ is an $F^{\times}I$-equivariant map from $\FF_{\sigma_{1}}$ to $\FF_{\sigma_0}$. If $v\in D_1$ then it sends
$$r^{\sigma_1}_{\sigma_0}:f_v\mapsto f_{r(v)}.$$
We observe, for the purposes of Theorem \ref{equiv}, that 
$$\tilde{D}=(\FF_{\sigma_0}, \FF_{\sigma_1}, r^{\sigma_1}_{\sigma_0})$$
 is an object of $\diag$. Moreover, 
$\ev=(\ev_0, \ev_1)$ is an isomorphism of diagrams between $D$ and $\tilde{D}$. We will show later on that $\ev$ induces a natural transformation
 between certain functors.

 Let $\tau'$ be an edge containing a vertex $\tau$, then there exists $g\in G$,
such that $\tau=g \sigma_0$ and $\tau'=g \sigma_1$. Moreover, $g$ can only be replaced by $gk$, where $k\in \kzr \cap \kon= F^{\times} I$. We define 
$$r^{\tau'}_{\tau} : \FF_{\tau'} \rightarrow \FF_{\tau},\quad f\mapsto g r^{\sigma_1}_{\sigma_0}( g^{-1}f)$$
where  $g$ acts on the space $\cIndu{\kzr}{G}{D_{0}}$ and $g^{-1}$  on the space $\cIndu{\kon}{G}{D_{1}}$. Since, $r$ is 
$F^{\times}I$-equivariant we have 
 $$\rho_0(k)\circ r^{\sigma_1}_{\sigma_0} \circ\rho_1(k^{-1}) =r^{\sigma_1}_{\sigma_0} $$
 for all $k\in F^{\times} I$. Hence, the map $r^{\tau'}_{\tau}$ is independent of the choice of $g$. Explicitly, let $v=f(g^{-1})$, then
$$r^{\tau'}_{\tau} :f\mapsto g f_{r(v)}.$$
 Let $\tau$ be any simplex  then we define the map $r^{\tau}_{\tau}=\id_{\FF_{\tau}}$.   
\subsubsection{$G$-action}
So far from a diagram we have constructed a coefficient system. We need to show that $G$ acts on it. Let $i\in\{0,1\}$ and let
 $f\in\cIndu{\KK(\sigma_i)}{G}{D_{i}}$. For any $g\in G$ we have
$$\supp (gf)=(\supp f)g^{-1}.$$
Hence for any simplex $\tau$ we obtain  a linear map
$$g_{\tau}:\FF_{\tau}\rightarrow \FF_{g \tau},\quad f\mapsto gf.$$
Moreover, $1_{\tau}=\id_{\FF_{\tau}}$ and $g_{h\tau}\circ h_{\tau}=(gh)_{\tau}$, for any $g$, $h\in G$.
Let $\tau'$ be an edge containing a vertex  $\tau$. We need to show that the  diagram:
\begin {displaymath}
\xymatrix {\FF_{\tau} \ar[r]^{g_{\tau}} & \FF_{g\tau} \\
          \FF_{\tau'} \ar[u]^{r^{\tau'}_{\tau}} \ar[r]^{g_{\tau'}} & \FF_{g\tau'} \ar[u]_{r^{g\tau'}_{g\tau}}}
\end {displaymath}
commutes. There exists $g_1\in G$ such that $\tau=g_1\sigma_0$ and $\tau'=g_1 \sigma_1$. Moreover, such $g_1$ is determined up to a multiple 
$g_1 k$, where $k\in F^{\times}I$. Let $f\in \FF_{\tau'}$ and let $v=f(g_{1}^{-1})$, then
$$r^{\tau'}_{\tau}(f)=g_1 f_{r(v)}.$$
Hence
$$(g_{\tau}\circ r^{\tau'}_{\tau})(f)=g g_1 f_{r(v)}.$$
Since $g\tau'=gg_1 \sigma_1$, $g\tau=gg_1\sigma_0$ and $(g f)( (g g_1)^{-1})= f(g_{1}^{-1})=v$ we obtain
$$(r^{g\tau'}_{g\tau}\circ g_{\tau'})(f)= r^{g\tau'}_{g\tau}(g g_1 f_v)= g g_1 f_{r(v)}.$$
Hence the diagram commutes.
\subsubsection{Morphisms}\label{morph} Let $D'=(D'_{0}, D'_{1}, r')$ be another diagram, let $\psi=(\psi_0, \psi_1)$ be a morphism of diagrams
$$\psi: D\rightarrow D'$$
and let $\FF'=(\FF'_{\tau})_{\tau}$ be a coefficient system associated to $D'$ via the construction above.
 Let $\tau$ be any simplex on the tree. If 
$\tau$ is a vertex let $i=0$ and  if $\tau$ is an edge, let $i=1$. There exists some $g\in G$ such that $\tau=g \sigma_i$. Let $f\in V_{\tau}$
and let $v=f(g^{-1})$ we define a map 
$$\psi_{\tau}:\FF_{\tau}\rightarrow \FF'_{\tau},\quad f\mapsto g f_{\psi_{i}(v)}$$
where $f_{\psi_i(v)}$ is the unique function in $\FF'_{\sigma_i}$, such that  $f_{\psi_i(v)}(1)=\psi_i(v)$. 
Since the map $\psi_i$ is $\KK(\sigma_i)$-equivariant, $\psi_{\tau}$ is independent of the choice of $g$.

We will show that the maps $(\psi_{\tau})_{\tau}$ are compatible with the restriction maps. Let $\tau'$ be an edge containing a vertex 
$\tau$. We claim that  the diagram
\begin {displaymath}
\xymatrix {\FF_{\tau} \ar[r]^{\psi_{\tau}} & \FF'_{\tau} \\
          \FF_{\tau'} \ar[u]^{r^{\tau'}_{\tau}} \ar[r]^{\psi_{\tau'}} & \FF'_{\tau'} \ar[u]_{(r')^{\tau'}_{\tau}}}
\end {displaymath}
commutes. There exists $g\in G$ such that $\tau=g\sigma_0$ and $\tau'=g\sigma_1$. Let $f\in \FF_{\tau'}$ and let $v=f(g^{-1})$. Then
$$(\psi_{\tau}\circ r^{\tau'}_{\tau})(f)=\psi_{\tau}(g f_{r(v)})=g f_{\psi_0(r(v))}$$
and 
$$((r')^{\tau'}_{\tau}\circ\psi_{\tau'})(f)=(r')^{\tau'}_{\tau}(g f_{\psi_1(v)})= g f_{r'(\psi_1(v))}.$$
Since $(\psi_0,\psi_1)$ is a morphism of diagrams
$$\psi_0(r(v))=r'(\psi_1(v)).$$
Hence the diagram commutes as claimed and $(\psi_{\tau})_{\tau}$ are compatible with the restriction maps.  

Finally, we will show that the maps $(\psi_{\tau})_{\tau}$ are compatible with the $G$-action. Let $\tau$ be any  simplex on the tree.
To ease the notation, for every $h\in G$ we denote by $h_{\tau}$ the action of $h$ on both $(\FF_{\tau})_{\tau}$ and  $(\FF'_{\tau})_{\tau}$.
Let $\tau$ be a simplex on the tree $X$ and let $h\in G$. We claim that the diagram
\begin {displaymath}
\xymatrix {\FF_{h\tau} \ar[r]^{\psi_{h\tau}} & \FF'_{h\tau} \\
          \FF_{\tau} \ar[u]^{h_{\tau}} \ar[r]^{\psi_{\tau}} & \FF'_{\tau} \ar[u]_{h_{\tau}}}
\end {displaymath}
commutes. If $\tau$ is an edge let $i=1$, if $\tau$ is a vertex let $i=0$. There exists $g\in G$, such that  $\tau=g\sigma_i$. Let $f\in \FF_{\tau}$ and
 let $v=f(g^{-1})$, then
$$ \psi_{h\tau}(h_{\tau}(f))=\psi_{h\tau}(hg f_{v})= hg f_{\psi_i(v)}$$
and 
$$ h_{\tau}(\psi_{\tau}(f))=h_{\tau}(g f_{\psi_i(v)})=hg f_{\psi_i(v)}.$$ 
Hence, the diagram commutes as claimed and  the collection $(\psi_{\tau})_{\tau}$ defines a morphism of equivariant coefficient systems. 
\subsubsection{Equivalence} 
\begin{defi} Let $\CC$ be a functor
$$\CC:\diag\rightarrow \coeff$$
which sends a diagram $D$ to the coefficient system $(\FF_{\tau})_{\tau}$ as above.
\end{defi}
One needs to check that given three diagrams and two morphisms between them
\begin{displaymath}
\xymatrix{D\ar[r]^{\psi}& D'\ar[r]^{\psi'} &D''}
\end{displaymath}
we have 
$$\CC(\psi'\circ\psi)=\CC(\psi')\circ\CC(\psi).$$
However, that is immediate from the construction of $\CC(\psi)$ in Section \ref{morph}.
\begin{thm}\label{equiv} The functors $\CC$ and $\DD$ induce an equivalence of categories between $\diag$ and $\coeff$.
\end{thm}
\begin{proof} Let $D=(D_0,D_1,r)$ be an object in $\diag$. Then
 $$(\DD\circ\CC)(D)=\tilde{D}=(\FF_{\sigma_0},\FF_{\sigma_1}, r^{\sigma_1}_{\sigma_0})$$
 with the notation of Section \ref{restmap}. The isomorphism 
$$\ev:\tilde{D}\cong D$$
 of Section \ref{restmap} is given by the evaluation at $1$. We claim that it induces an isomorphism of
functors between $\DD\circ\CC$ and $\id_{\diag}$. We only need to check what happens to morphisms. Let $D'=(D'_0, D'_1, r')$ be another object in 
the category of diagrams and let $\psi=(\psi_0, \psi_1)$ be a morphism 
$$\psi: D\rightarrow D'.$$
Let $(\DD\circ\CC)(D')=\tilde{D}'=(\FF'_{\sigma_0}, \FF'_{\sigma_1}, (r')^{\sigma_1}_{\sigma_0})$ and
 let
 $$(\DD\circ\CC)(\psi)=\tilde{\psi}=(\tilde{\psi}_0,\tilde{\psi}_1)$$
 be a morphism induced by a functor $\DD\circ\CC$. We need to show that the diagram:
\begin {displaymath}
\xymatrix {\tilde{D}'\ar[r]^{\ev} & D' \\
          \tilde{D} \ar[u]^{\tilde{\psi}} \ar[r]^{\ev} & D \ar[u]_{\psi}}
\end {displaymath}
commutes. Let $i\in\{0,1\}$, let $f\in \FF_{\sigma_i}$ and let $v=f(1)$ then
$$(\psi_i\circ\ev_i)(f)=\psi_i(v).$$
From Section \ref{morph} $\tilde{\psi}_i(f)$ is the unique function in $\FF'_{\sigma_i}$, taking value $\psi_i(v)$ at $1$. Hence
$$(\ev_i\circ\tilde{\psi}_i)(f)= \psi_i(v).$$
This implies that the diagram commutes. 

Conversely, we need to show that the functor $\CC\circ \DD$ is isomorphic to $\id_{\coeff}$. Let $\VV=(V_{\tau})_{\tau}$ be a $G$-equivariant
coefficient system with the restriction maps  $t^{\tau'}_{\tau}$. 
Then $\DD(\VV)$ is a diagram given by:
\begin {displaymath}
\xymatrix {V_{\sigma_0} \\
          V_{\sigma_1} \ar[u]^{t^{\sigma_1}_{\sigma_0}}}.
\end {displaymath}  
Let  $k\in \kzr$ then it acts on $V_{\sigma_0}$ by a linear map $k_{\sigma_0}$. Similarly,
 if $k\in \kon$ then it acts on  on $V_{\sigma_1}$ by a 
linear map $k_{\sigma_1}$. Let
 $$(\CC\circ\DD)(\VV)=\FF=(\FF_{\tau})_{\tau}$$
 with the restriction maps $r^{\tau'}_{\tau}$. We will construct a canonical isomorphism $\ev=(\ev_{\tau})_{\tau}$ 
$$\ev: \FF\cong \VV$$
of $G$ equivariant coefficient systems. Let $\tau$ be a simplex on the tree. If $\tau$ is a vertex
let $i=0$ and if $\tau$ is an edge let $i=1$. There exists $g\in G$ such that $\tau=g\sigma_i$. 
For  $f\in \FF_{\tau}$ we let $v=f(g^{-1})$. Then $v$ is a vector in $V_{\sigma_i}$. We define a map $\ev_{\tau}$, by
$$\ev_{\tau}:\FF_{\tau}\rightarrow V_{\tau},\quad f\mapsto g_{\sigma_i} v$$
where $g_{\sigma_i}$ is the linear map coming from the $G$ action on $\VV$. If we replace $g$ by $gk$, for some $k\in \KK(\sigma_i)$, then 
$$ (gk)_{\sigma_i} (f((gk)^{-1}))= (g_{\sigma_{i}}\circ k_{\sigma_i}\circ k^{-1}_{\sigma_{i}})(f(g^{-1}))= g_{\sigma_i}(f(g^{-1})).$$
Hence, the  map $\ev_{\tau}$ is independent of the choice of $g$. Moreover, $\ev_{\tau}$ is an isomorphism of vector spaces with the inverse 
given as follows. Let $w\in V_{\tau}$, let $v=(g^{-1})_{\tau} w$, then 
$v$ is a vector in $W_{\sigma_{i}}$. Let $f_v$ be the unique  function in $ \FF_{\tau}$ such that $f_{v}(1)=v$. Then $(\ev_{\tau})^{-1}$ is given by
$$(\ev_{\tau})^{-1}:V_{\tau}\rightarrow \FF_{\tau},\quad w\mapsto g f_v$$
where the action by $g$ is on the space $ \cIndu{\KK(\sigma_i)}{G}{V_{\sigma_i}}$.

 The collection of maps $(\ev_{\tau})_{\tau}$ is $G$-equivariant. Let $h\in G$, then $h f$ belongs to the space $\FF_{h\tau}$ and
$$ \ev_{h\tau}(hf)= (hg)_{\sigma_i}((hf)((hg)^{-1}))=(h_{\tau}\circ g_{\sigma_i})(f(g^{-1}))
=h_{\tau}( \ev_{\tau}(f)).$$

We need to show that the maps $\ev_{\tau}$ are compatible with the restriction maps. Let $\tau'$ be an edge containing a vertex $\tau$.
 We need to show that the diagram 
\begin {displaymath}
\xymatrix {\FF_{\tau}\ar[r]^{\ev_{\tau}} & V_{\tau} \\
          \FF_{\tau'} \ar[u]^{r^{\tau'}_{\tau}} \ar[r]^{\ev_{\tau'}} & V_{\tau'}\ar[u]_{t^{\tau'}_{\tau}}}
\end {displaymath}  
commutes. There exists $g\in G$ such that $\tau=g \sigma_0$ and $\tau'=g \sigma_1$. Let $f$ be a
 function in $\FF_{\tau'}$. Let $v_1=f(g^{-1})$ , then $v_1$ is a 
vector in $V_{\sigma_1}$. Let $v_0= t^{\sigma_1}_{\sigma_0}(v_1)$. Then $r^{\tau'}_{\tau}(f)$ is the unique function of $\FF_{\tau}$
 taking value  $v_0$ at $g^{-1}$. Hence
$$(\ev_{\tau}\circ r^{\tau'}_{\tau})(f)= g_{\sigma_{0}}v_0.$$
On the other hand
$$(t^{\tau'}_{\tau}\circ\ev_{\tau'})(f)= t^{\tau'}_{\tau}(g_{\sigma_1} v_1).$$
The action of $G$ on $\VV$ respects the restriction maps, in the sense that  the diagram: 
\begin{displaymath}
\xymatrix {V_{\sigma_0}\ar[r]^{g_{\sigma_0}} & V_{\tau} \\
          V_{\sigma_1} \ar[u]^{t^{\sigma_1}_{\sigma_0}} \ar[r]^{g_{\sigma_1}} & V_{\tau'}\ar[u]_{t^{\tau'}_{\tau}}}.
\end{displaymath}  
commutes. Hence,
$$ t^{\tau'}_{\tau}(g_{\sigma_1}v_1)=g_{\sigma_0}v_0.$$
Hence our original diagram commutes and $\ev=(\ev_{\tau})_{\tau}$ defines an isomorphism of $G$-equivariant coefficient systems. 

In order to show that the morphism $\ev$ induces an isomorphism of functors between $\CC\circ \DD$ and $\id_{\coeff}$ we need to check 
what happens to the morphisms. However the proof is almost identical to the one given for $\diag$ so we omit it.  
\end{proof}

\begin{cor}\label{thepoint}Let $(\rho_0, V_0)$ be a smooth representation of $\kzr$ and $(\rho_1, V_1)$ a smooth representation of $\kon$. Suppose
that there exists an $F^{\times}I$-equi\-va\-riant isomorphism 
$$r:V_1\cong V_0,$$
then there exists a unique (up to isomorphism) smooth representation $\pi$ of $G$, such that 
$$\pi |_{\kzr}\cong \rho_0,\quad \pi |_{\kon}\cong \rho_1$$
and the diagram 
\begin {displaymath}
\xymatrix {V_{0} \ar[r]^{\cong} & \pi \\
          V_{1} \ar[u]^{r} \ar[r]^{\cong} & \pi \ar[u]_{\id}}
\end {displaymath}
of $F^{\times}I$-representations commutes.
\end{cor}
\begin{proof} Let $D$ be the object in $\diag$, given by $D=(V_0,V_1, r)$. 
Let $\CC(D)$ be a coefficient system corresponding to $D$, with the restriction maps $r^{\tau'}_{\tau}$. 
Since $(\DD\circ\CC)(D)\cong D$ and $r$ is an isomorphism, the map $r^{\sigma_1}_{\sigma_0}$ is an isomorphism and Proposition \ref{act} implies that 
$H_0(X,\CC(D))$ satisfies the conditions of the Corollary.

The statement of the Corollary can be rephrased as follows: there exists a unique up to isomorphism smooth representation $\pi$ of $G$, such that
$$D\cong \DD(\const_{\pi}).$$
If $\pi'$ was another such, then 
$$\DD(\const_{\pi'})\cong D\cong \DD(\const_{\pi}).$$
Hence, by Theorem \ref{equiv}
$$\const_{\pi'}\cong \const_{\pi}.$$
Lemma \ref{Hconst} implies that
$$\pi'\cong H_0(X,\const_{\pi'})\cong H_0(X,\const_{\pi})\cong \pi$$
and we obtain uniqueness. 
\end{proof}

\begin{remar}\label{explicitaction} Let $\tilde{W}$ be a subgroup of $G$ generated by $s$ and $\Pi$. The Iwahori decomposition says that
$G=I\tilde{W}I$. Let $\pi$ be a representation constructed as above, $v\in \pi$ and $g\in G$. Then $gv$ may be determined by decomposing 
$g=u_1 w u_2$, where $u_1, u_2\in I$, $w\in \tilde{W}$, and then chasing around the diagram. 

The simplest example illustrating \ref{thepoint} is the trivial diagram $\tilde{\Eins}=(\Eins,\Eins,\id)$. The proof of Corollary \ref{trivialrep}
can be reinterpreted as a construction of a morphism $\tilde{\Eins}\hookrightarrow \DD(\const_{\pi})$. This gives us an injection of $G$ representations
$$\Eins\cong H_0(X,\CC(\tilde{\Eins}))\hookrightarrow H_0(X,\const_{\pi})\cong \pi.$$
\end{remar}

\section{Supersingular representations}
\subsection{Coefficient systems $\VV_{\gamma}$}\label{VV}
Let $\chi:H\rightarrow \Fbar^{\times}$ be a character, and let $\rho_{\chi, J}$ be an irreducible representation of 
$\Gamma$, with the notations of Section \ref{GLFq}. We consider $\chi$ as a character of $I$ and $\rho_{\chi,J}$ as a representation of $K$, via 
$$K\rightarrow K/K_1\cong \Gamma \quad \textrm{and}\quad I\rightarrow I/I_1\cong H.$$
Let $\tilde{\rho}_{\chi, J}$ be the extension of $\rho_{\chi, J}$ to $F^{\times}K$ such that our fixed uniformiser $\pif$ acts trivially, and let 
$\tilde{\chi}$ be the extension of $\chi$ to $F^{\times}I$, such that $\pif$ acts trivially. The space of $I_1$-invariants of $\tilde{\rho}_{\chi,J}$ 
is one dimensional and $F^{\times}I$ acts on it via the character $\tilde{\chi}$. We fix a vector $v_{\chi, J}$ such that
$$\rho_{\chi,J}^{I_1}=\langle v_{\chi, J}\rangle_{\Fbar}.$$
\begin{lem}\label{ac} There exists a unique action of $\kon$  on $(\tilde{\rho}_{\chi,J}\oplus \tilde{\rho}_{\chi^{s}, \overline{J}})^{I_1}$, 
extending the action of  $F^{\times}I$, such that
$$\Pi^{-1}v_{\chi,J}=v_{\chi^{s},\overline{J}}\quad \textrm{and}\quad \Pi^{-1}v_{\chi^{s},\overline{J}}=v_{\chi, J}.$$
Moreover, with this action
$$(\tilde{\rho}_{\chi,J}\oplus \tilde{\rho}_{\chi^{s}, \overline{J}})^{I_1}\cong \Indu{F^{\times}I}{\kon}{\tilde{\chi}}$$
as $\kon$-representations.
\end{lem}
\begin{proof}We note that if $t\in T$ is a diagonal matrix then $\Pi t \Pi^{-1}= sts$, hence $(\tilde{\chi})^{\Pi}\cong \tilde{\chi^{s}}$ as 
representations of $F^{\times}I$ and  Mackey's decomposition gives us
$$(\Indu{F^{\times}I}{\kon}{\tilde{\chi}})|_{F^{\times}I}\cong \tilde{\chi}\oplus\tilde{\chi^{s}}.$$
Since 
$$(\tilde{\rho}_{\chi,J}\oplus \tilde{\rho}_{\chi^{s}, \overline{J}})^{I_1}\cong \tilde{\chi}\oplus\tilde{\chi^{s}}$$
as $F^{\times}I$-representation, we can extend the action. Explicitly, we consider  $f\in \Indu{F^{\times}I}{\kon}{\tilde{\chi}}$, such that
$\supp f =F^{\times}I$ and $f(g)=\tilde{\chi}(g)$, for all $ g\in F^{\times}I.$
 Then the map 
$$f \mapsto v_{\chi, J}, \quad \Pi^{-1}f\mapsto v_{\chi^{s},\overline{J}}$$
induces the required isomorphism. Since, $\Pi$ and $F^{\times}I$ generate $\kon$ the action is unique. 
\end{proof}

\begin{defi}\label{dgamma} Let $\chi:H\rightarrow \Fbar^{\times}$ be a character, and let $\gamma=\{\chi, \chi^{s}\}$ we define 
$D_{\gamma}$ to be an object in $\diag$,
 given by  
\begin {displaymath}
\xymatrix {\tilde{\rho}_{\chi, J} \oplus \tilde{\rho}_{\chi^{s}, \overline{J}} \\
          (\tilde{\rho}_{\chi, J} \oplus \tilde{\rho}_{\chi^{s}, \overline{J}})^{I_1} \ar[u]}
\end {displaymath}
where the action of $\kon$ on $(\tilde{\rho}_{\chi, J} \oplus \tilde{\rho}_{\chi^{s}, \overline{J}})^{I_1}$ is given by Lemma \ref{ac}.
 Moreover, we define
$\VV_{\gamma}$ to be a coefficient system, given by 
$$\VV_{\gamma}=\CC(D_{\gamma}).$$
\end{defi}
\begin{lem} The diagram $D_{\gamma}$ is independent up to isomorphism of the choices made for $v_{\chi,J}$ and  $v_{\chi^{s},\overline{J}}$.
\end{lem}
\begin{proof} Suppose that instead we choose vectors $v'_{\chi,J}$ and $v'_{\chi^{s},\overline{J}}$ and let $D'_{\gamma}$ be the corresponding diagram.
 Since, the spaces $\rho_{\chi,J}^{I_1}$ and $\rho_{\chi^{s},\overline{J}}^{I_1}$ are one dimensional there exist $\lambda, \mu \in \Fbar^{\times}$, 
such that 
$$\lambda v_{\chi,J}= v'_{\chi,J},\quad \mu v_{\chi^{s},\overline{J}}=v'_{\chi^{s},\overline{J}}.$$
The isomorphism
$$\lambda\id\oplus \mu\id: \tilde{\rho}_{\chi,J}\oplus\tilde{\rho}_{\chi^{s},\overline{J}}\rightarrow 
\tilde{\rho}_{\chi,J}\oplus\tilde{\rho}_{\chi^{s},\overline{J}}$$
induces an isomorphism of diagrams $D_{\gamma}\cong D'_{\gamma}$.
\end{proof}

Since $D_{\gamma}$ and $\DD(\VV_{\gamma})$ are canonically isomorphic, to ease the notation, we identify them. Let 
$\omega_{\chi, J}$, $\omega_{\chi^{s}, \overline{J}}\in C^{or}_c(X_{(0)},\VV_{\gamma})$ supported on a single vertex $\sigma_0$, such that
$$\omega_{\chi,J}(\sigma_0)= v_{\chi,J}\quad \textrm{and} \quad \omega_{\chi^{s},\overline{J}}(\sigma_0)=v_{\chi^{s}, \overline{J}}.$$
Let 
$$\bar{\omega}_{\chi, J}= \omega_{\chi, J} +\partial C^{or}_c(X_{(1)},\VV_{\gamma})\quad \textrm{and}\quad 
\bar{\omega}_{\chi^{s}, \overline{J}}= \omega_{\chi^{s}, \overline{J}} +\partial C^{or}_c(X_{(1)},\VV_{\gamma})$$
be their images in $H_0(X,\VV_{\gamma})$.
\begin{lem}\label{supermod} We have
$$\langle \bar{\omega}_{\chi, J}, \bar{\omega}_{\chi^{s},\overline{J}}\rangle_{\Fbar}\cong M_{\gamma}$$
as right $\HH$-modules. 
\end{lem}
\begin{proof} Since the restriction maps in $\VV_{\gamma}$ are injective, Lemma \ref{inj} says that $\bar{\omega}_{\chi,J}$ and 
$\bar{\omega}_{\chi^{s},\overline{J}}$ are non-zero. We have 
$$\langle v_{\chi, J}\rangle_{\Fbar}=(\tilde{\rho}_{\chi,J})^{I_1}\cong M_{\chi, J}\quad\textrm{and}\quad 
\langle v_{\chi^{s},\overline{ J}}\rangle_{\Fbar}=(\tilde{\rho}_{\chi^{s},\overline{J}})^{I_1}\cong M_{\chi^{s},\overline{J}}$$ 
as $\HH_K$-modules. Hence $\bar{\omega}_{\chi,J}$ and $\bar{\omega}_{\chi^{s},\overline{J}}$ are fixed by $I_1$ and 
$$\langle \bar{\omega}_{\chi,J}\rangle_{\Fbar}\oplus \langle \bar{\omega}_{\chi^{s},\overline{J}}\rangle_{\Fbar}\cong M_{\chi, J} \oplus 
M_{\chi^{s}, \overline{J}}$$
as $\HH_K$-modules. Corollary \ref{T} and Lemma \ref{act} imply that
$$ \bar{\omega}_{\chi, J}T_{\Pi}= \Pi^{-1}\bar{\omega}_{\chi, J}= \bar{\omega}_{\chi^{s},\overline{J}}\quad \textrm{and}\quad  
\bar{\omega}_{\chi^{s}, \overline{J}}T_{\Pi}= \Pi^{-1}\bar{\omega}_{\chi^{s},\overline{ J}}= \bar{\omega}_{\chi,J}.$$
Hence 
$$\langle \bar{\omega}_{\chi, J}, \bar{\omega}_{\chi^{s},\overline{J}}\rangle_{\Fbar}\cong M_{\gamma}$$
as $\HH$-modules. 
\end{proof}

\begin{lem}\label{omegagener} The vector $\bar{\omega}_{\chi,J}$ (resp. $\bar{\omega}_{\chi^{s},\overline{J}}$) generates $H_0(X,\VV_{\gamma})$ as
a $G$-representation. 
\end{lem}
\begin{proof} Lemma \ref{act} implies that  $ \Pi^{-1} \bar{\omega}_{\chi,J}=\bar{\omega}_{\chi^{s},\overline{J}}.$
Hence, it is enough to show that $\omega_{\chi,J}$ and $\omega_{\chi^{s},\overline{J}}$ generate $C^{or}_c(X_{(0)},\VV_{\gamma})$ as a 
$G$-representation. Since, $\rho_{\chi,J}$ and $\rho_{\chi^{s},\overline{J}}$ are irreducible $K$-representations, $\omega_{\chi,J}$ and 
$\omega_{\chi^{s},\overline{J}}$ will generate the space
$$C_{c}^{or}(\sigma_0,\VV_{\gamma})= \{\omega\in C^{or}_{c}(X_{(0)},\VV_{\gamma}):\supp \omega \subseteq \sigma_0\}$$
as a $K$-representation. Since the action of $G$ on the vertices of $X$ is transitive, the space $C^{or}_{c}(\sigma_0,\VV_{\gamma})$ will generate 
$C^{or}_{c}(X_{(0)},\VV_{\gamma})$ as a $G$-representation.
\end{proof}

\begin{cor}\label{superquot} Let $\pi$ be a non-zero irreducible quotient of $H_0(X,\VV_{\gamma})$, then $\pi$ is a supersingular representation. 
\end{cor}
\begin{proof} Lemma \ref{omegagener} implies that the images of $\bar{\omega}_{\chi,J}$ and $\bar{\omega}_{\chi^{s},\overline{J}}$ in $\pi$ are
 non-zero. Hence, by Lemma \ref{supermod}, $\pi^{I_1}$ will contain a supersingular module $M_{\gamma}$, 
then Corollary \ref{supersingularity} implies that $\pi$ is supersingular. 
\end{proof}

\begin{prop}\label{homHz} Let $\pi$ be a smooth representation of $G$ and suppose that there exists $v_1, v_2\in \pi^{I_1}$ such that 
$$ \langle K v_1\rangle_{\Fbar}\cong \rho_{\chi, J},\quad \langle K v_2 \rangle_{\Fbar}\cong \rho_{\chi^{s}, \overline{J}},
\quad \Pi^{-1} v_1= v_2, \quad \Pi^{-1} v_2 = v_1,$$
then there exists a $G$-equivariant map $\phi: H_0(X,\VV_{\gamma})\rightarrow \pi $
such that
$$\phi(\bar{\omega}_{\chi,J})=v_1\quad\textrm{and}\quad \phi(\bar{\omega}_{\chi^{s},\overline{J}})=v_2$$
where $\gamma=\{\chi, \chi^{s}\}$.
\end{prop} 
\begin{proof} By Lemma \ref{Hconst} and Theorem \ref{equiv}, it is enough to construct a morphism of diagrams $D_{\gamma}\rightarrow \DD(\const_{\pi})$.
However, such morphism is immediate.
\end{proof}

\begin{cor}\label{allnice}Let $\pi$ be a smooth representation of $G$ and suppose that  one of the following holds: $\chi=\chi^{s}$, or $p=q$, then
$$\Hom_G(H_0(X,\VV_{\gamma}),\pi)\cong \Hom_{\HH}(M_{\gamma}, \pi^{I_1}).$$
\end{cor}
\begin{remar}This fails if $q\neq p$ and $\chi\neq\chi^{s}$. Proposition \ref{Injmod} gives an example.
\end{remar}
\begin{proof} Lemmas \ref{supermod} and \ref{omegagener} imply that we always have an injection
$$\Hom_G(H_0(X,\VV_{\gamma}),\pi)\hookrightarrow \Hom_{\HH}(M_{\gamma}, \pi^{I_1}).$$
By Lemma \ref{restML} $M_{\gamma}|_{\HH_K}\cong M_{\chi,J}\oplus M_{\chi^{s},\overline{J}}$. Under the assumptions made, Corollaries \ref{T}, 
\ref{splitisgood} and respectively \ref{qequalpisgood} give us vectors $v_1, v_2\in \pi^{I_1}$ as in
 Proposition \ref{homHz}, hence  the injection is an isomorphism.
\end{proof}

\begin{cor}\label{shouldbeall}Let $\pi$ be a smooth representation, and suppose that $\pi^{I_1}\cong M_{\gamma}$, then
$$\dim \Hom_{G}(H_0(X,\VV_{\gamma}),\pi)=1.$$
\end{cor}
\begin{proof} It is enough to consider the case $p\neq q$ and $\chi\neq\chi^{s}$. Since Corollary \ref{allnice} implies the statement in the other cases.
Let $\rho=\langle K \pi^{I_1}\rangle_{\Fbar}$, then $\rho^{I_1}=\pi^{I_1}$. Hence
$$\rho^{I_1}\cong M_{\gamma}|_{\HH_K}\cong M_{\chi,\emptyset}\oplus M_{\chi^{s},\emptyset}$$ as 
an $\HH_K$-module. Proposition \ref{consol} implies that $\rho\cong\rho_{\chi,\emptyset}\oplus \rho_{\chi^{s},\emptyset}$. The action of $\Pi$ on 
$\pi^{I_1}$ is given by Corollary \ref{T}. Now we may apply Proposition
\ref{homHz} to get a non-zero homomorphism. So the dimension is at least one. The  module $M_{\gamma}$ is irreducible, and Lemmas \ref{supermod} and 
\ref{omegagener} imply that the dimension is at most one. 
\end{proof}

\subsection{Injective envelopes}\label{injectenve}
For the convenience of the reader we recall some general facts about injective envelopes. Let $\kK$ be a pro-finite group and let $\Rep_{\kK}$ 
be the category of smooth $\Fbar$-representations of $\kK$. We assume that
$\kK$ has an open normal pro-$p$ subgroup $\PP$. 
\begin{defi} Let $\pi\in \Rep_{\kK}$ and let $\rho$ be a $\kK$-invariant subspace of $\pi$. We say that $\pi$ is an essential extension of $\rho$
if for every non-zero $\kK$-invariant subspace $\pi'$ of $\pi$, we have $\pi'\cap \rho\neq 0$.

Let $\rho\in \Rep_{\kK}$ and let $\Inj$ be an injective object in $\Rep_{\kK}$. A monomorphism $\iota: \rho \hookrightarrow \Inj$ is called an injective 
envelope of $\rho$, if $\Inj$ is an essential extension of $\iota(\rho)$.
\end{defi}

\begin{prop} Every representation $\rho\in \Rep_{\kK}$ has an injective envelope $\iota: \rho \hookrightarrow \Inj \rho$. Moreover, 
 injective envelopes are unique up to isomorphism.
\end{prop}
\begin{proof} \cite{sw}, \S 3.1.
\end{proof} 

\begin{lem}\label{mininjobj} Let $\Inj$ be an injective object in $\Rep_{\kK}$ and let $\iota: \rho \rightarrow \Inj \rho$ 
be an injective envelope of $\rho$ in
$\Rep_{\kK}$. Let  $\phi$ be a monomorphism $\phi: \rho \hookrightarrow \Inj$, then there exists a monomorphism 
$\psi: \Inj \rho \hookrightarrow  \Inj$ such that $\phi= \psi \circ \iota$.
\end{lem}
\begin{proof} Since $\Inj$ is an injective object there exists $\psi$ such that the diagram 
\begin{displaymath}
\xymatrix{ 0\ar[r]& \rho\ar[r]^{\iota}\ar[d]_{\phi}  & \Inj \rho \ar@{.>}[dl]^{\psi}\\
                  & \Inj          &   }
\end{displaymath}
of $\kK$-representations commutes. Since $\phi$ is an injection  $\Ker \psi\cap \iota(\rho)= 0$. This implies that 
$\Ker \psi =0$, as $\Inj \rho$ is an essential extension of $\iota(\rho)$.
\end{proof}

\begin{lem}\label{injinv} Let $\rho\in Rep_{\kK}$ be an irreducible representation and let $\iota:\rho\hookrightarrow \Inj \rho$ be an injective 
envelope of $\rho$ in $\Rep_{\kK}$, then $\rho\hookrightarrow (\Inj \rho)^{\PP}$ is an injective envelope of $\rho$ in $\Rep_{\kK/\PP}$.
\end{lem}
\begin{proof} 
We note that since $\PP$ is an open normal pro-$p$ subgroup of $\kK$ and $\rho$ is irreducible, 
Lemma \ref{propinvariants}  implies that $\PP$ acts trivially on $\rho$. Hence, $\iota(\rho)$ is a subspace of $(\Inj \rho)^{\PP}$. Moreover,
$(\Inj \rho)^{\PP}$ is an essential extension of $\iota(\rho)$, since 
$\Inj \rho$ is an essential extension of $\iota(\rho)$. 

Let $\LLL: \Rep_{\kK/\PP}\rightarrow \Rep_{\kK}$ be a functor sending a representation $\xi$ to its inflation $\LLL(\xi)$ 
to a representation of $\kK$, via $\kK\rightarrow \kK/\PP$. Then 
$$\Hom_{\kK/\PP}(\xi,(\Inj \rho)^{\PP})\cong \Hom_{\kK}(\LLL(\xi),\Inj \rho)$$
where the isomorphism is canonical. Since, the functor $\LLL$ is exact and $\Inj \rho$ is an injective object in $\Rep_K$, 
the functor $\Hom_{\kK/\PP}( * ,(\Inj \rho)^{\PP})$ is 
exact. Hence,  $(\Inj \rho)^{\PP}$ is an injective object in $\Rep_{\kK/\PP}$, which establishes the Lemma.
\end{proof}

\begin{defi} Let $\pi\in \Rep_{\kK}$, we denote by $\soc \pi$ the subspace of $\pi$, generated by all irreducible subrepresentations of $\pi$.
\end{defi}

\begin{lem}\label{socle} Let $\rho\in \Rep_{\kK}$ be irreducible, and let $\iota: \rho\hookrightarrow \Inj \rho$ be an injective envelope of $\rho$, then
$\soc(\Inj \rho)\cong \rho$.
\end{lem}
\begin{proof} Let $\tau$ be any non-zero $\kK$ invariant subspace of $\Inj \rho$, which is irreducible as a representation of $\kK$. 
Since $\Inj \rho$ is an essential extension of $\iota(\rho)$ and $\rho$ is irreducible, we have $\tau=\iota(\rho)$. Hence, $\soc(\Inj \rho)=\iota(\rho)$.
\end{proof}
\subsubsection{Admissibility}\label{admissibility}
Let $\GG$ be a locally pro-finite group and let $\Rep_{\GG}$ be the category of smooth $\Fbar$-representations of $\GG$.
\begin{defi} A representation $\pi\in \Rep_{\GG}$ is called admissible, if for every open subgroup $\kK$ of $\GG$, the space $\pi^{\kK}$ of 
$\kK$-invariants is finite dimensional.
\end{defi}

\begin{lem}\label{ohdear} Suppose that $\GG$ has an open pro-$p$ subgroup $\PP$. A representation  $\pi\in \Rep_{\GG}$ is admissible 
if and only if $\pi^{\PP}$ is finite dimensional.
\end{lem}
\begin{proof} If $\pi$ is admissible, then $\pi^{\PP}$ is finite dimensional. Suppose that $\pi^{\PP}$ is finite dimensional and let 
$\Eins \hookrightarrow \Inj \Eins$ be an injective envelope of the trivial representation in $\Rep_{\PP}$, then there exists $\psi$, such that the diagram 
\begin{displaymath}
\xymatrix{ 0\ar[r]& \pi^{\PP}\ar[r]\ar[d]  & \pi|_{\PP} \ar@{.>}[dl]^{\psi}\\
                  & (\dim \pi^{\PP})\Inj \Eins         &   }
\end{displaymath}
of $\PP$-representations commutes. This implies that  $(\Ker \psi)^{\PP}=0$, and hence by Lemma \ref{propinvariants}, $\psi$ is injective.

Let $\kK$ be any open subgroup of $\GG$. Since  $\PP$ is an open compact subgroup of $\GG$,  we may choose an open  subgroup $\PP'$ of $\GG$ such that
$\PP'$ is a subgroup of $\PP\cap \kK$ and $\PP'$ is normal in $\PP$. It is enough to show that $\pi^{\PP'}$ is finite dimensional. Since $\psi$ is 
an injection, it is enough to show that $(\Inj \Eins)^{\PP'}$ is finite dimensional. Since $\PP$ is pro-$p$ and $\PP'$ is a normal open subgroup
of $\PP$, Lemma \ref{injinv} and Proposition \ref{simplep} imply that 
$$(\Inj \Eins)^{\PP'}\cong \Fbar[\PP/\PP']$$
which is finite dimensional.
\end{proof} 
\subsection{Coefficient systems $\II_{\gamma}$}\label{II}
Let $\chi:H\rightarrow \Fbar^{\times}$ be a character, and let  
\begin{displaymath}
\rho_{\chi,J}\hookrightarrow \Inj \rho_{\chi,J},\quad \rho_{\chi^{s},\overline{J}}\hookrightarrow \Inj \rho_{\chi^{s},\overline{J}}
\end{displaymath} 
be injective envelopes of $\rho_{\chi,J}$ and $\rho_{\chi^{s},\overline{J}}$ in $\Rep_K$, respectively.
 We may extend the action of $K$ to the action of $F^{\times}K$, so that our fixed uniformiser $\pif$ acts trivially. We get an exact sequence
\begin{displaymath}
\xymatrix{0\ar[r]&\tilde{\rho}_{\chi,J} \oplus \tilde{\rho}_{\chi^{s},\overline{J}}\ar[r]& \widetilde{\Inj} \rho_{\chi,J}\oplus 
\widetilde{\Inj} \rho_{\chi^{s},\overline{J}}}
\end{displaymath} 
of $F^{\times}K$-representations. This gives a commutative diagram  
\begin{displaymath}
\xymatrix{0\ar[r]&\tilde{\rho}_{\chi,J} \oplus \tilde{\rho}_{\chi^{s},\overline{J}}\ar[r]& \widetilde{\Inj} \rho_{\chi,J}\oplus 
\widetilde{\Inj} \rho_{\chi^{s},\overline{J}}\\
0\ar[r]&(\tilde{\rho}_{\chi,J} \oplus \tilde{\rho}_{\chi^{s},\overline{J}})^{I_1}\ar[r]\ar[u]& \widetilde{\Inj} \rho_{\chi,J}\oplus 
\widetilde{\Inj} \rho_{\chi^{s},\overline{J}}\ar[u]}
\end{displaymath} 
of $F^{\times}I$-representations. We will show that we may extend the action of $F^{\times}I$ on $(\widetilde{\Inj} \rho_{\chi,J}\oplus 
\widetilde{\Inj} \rho_{\chi^{s},\overline{J}})|_{F^{\times}I}$ to the action of $\kon$, so that we get an object $Y_{\gamma}$ in $\diag$, together
with an embedding $D_{\gamma}\hookrightarrow Y_{\gamma}$. Since the categories $\diag$ and $\coeff$ are equivalent, this will gives us an embedding of 
coefficient systems $\VV_{\gamma}\hookrightarrow \II_{\gamma}$. We will show that the image
$$\pi_{\gamma}=\Image(H_0(X,\VV_{\gamma})\rightarrow H_0(X,\II_{\gamma}))$$
is an irreducible supersingular representation of $G$. All the hard work was done in Propositions \ref{injiwahori} and \ref{injregular}, the 
construction of $Y_{\gamma}$ and the proof of irreducibility follow from the 'general non-sense' of Section \ref{injectenve}. This gives hope 
that similar construction might work for other groups. 

\begin{lem}\label{restI} Let $\rho$ be an irreducible representation of $K$ and let
\begin{displaymath}
\rho\hookrightarrow \Inj \rho
\end{displaymath}
 be an injective envelope of $\rho$ in $\Rep_K$, then
$$(\Inj \rho)|_{I}\cong \bigoplus_{\chi}\dim \Hom_{H}(\chi, (\inj \rho)^{U})\Inj \chi$$
where the sum is taken over all  irreducible representations of $H$, which we  identify with the irreducible representations of $I$ and
\begin{displaymath}
\rho\hookrightarrow \inj \rho, \quad \chi\hookrightarrow \Inj \chi
\end{displaymath} 
are the injective envelopes of $\rho$ in $\Rep_{\Gamma}$ and of $\chi$ in $\Rep_{I}$, respectively.
\end{lem}
\begin{proof}If $\chi$ is an irreducible representation of $I$, then Lemma \ref{propinvariants} implies that $I_1$ acts trivially on $\chi$. Since
$I/I_1\cong H$, the irreducible representations of $I$ and $H$ coincide. Moreover, since $H$ is abelian, all the 
irreducible representations of $H$ are one dimensional. Since, the order of $H$ is prime to $p$, all the representations of $H$ are semi-simple. 
Therefore
$$(\Inj \rho)^{I_1}\cong \bigoplus_{\chi} m_{\chi} \chi$$
as a representation of $I$, where the multiplicity $m_{\chi}$ of $\chi$ is given by 
$$m_{\chi}=\dim \Hom_{I}(\chi, \Inj \rho).$$ 
Lemma \ref{injinv} implies that $(\Inj \rho)^{K_1}\cong \inj \rho$ as representations of $K/K_1\cong \Gamma$. Corollary \ref{injsumm} implies that
$\inj \rho$ is finite dimensional. In particular, $m_{\chi}$ is finite for every $\chi$. Moreover,
$$\Hom_{I}(\chi, \Inj \rho)\cong \Hom_I(\chi, (\Inj \rho)^{K_1})\cong\Hom_B(\chi, \inj \rho) \cong \Hom_H(\chi, (\inj \rho)^U).$$
Hence, $m_{\chi}=\dim \Hom_H(\chi, (\inj \rho)^U)$. We consider an exact sequence
\begin{displaymath}
\xymatrix{0\ar[r]&(\Inj \rho)^{I_1}\ar[r]& (\Inj \rho)|_{I}}
\end{displaymath}
of $I$-representations. The restriction $(\Inj \rho)|_{I}$ is an injective object in $\Rep_{I}$. Lemma \ref{mininjobj} implies that
$$(\Inj \rho)|_I\cong \mathcal N \oplus \bigoplus_{\chi} m_{\chi} \Inj \chi$$
for some representation $\mathcal N$. Since $\Rep_H$ is semi-simple and $\Inj \chi$ is an essential extension of $\chi$, 
Lemma \ref{injinv} implies that $(\Inj \chi)^{I_1}\cong \chi$. By comparing the dimensions of $I_1$-invariants of both sides 
we get that $\dim \mathcal N^{I_1}=0$ and Lemma \ref{propinvariants} implies that $\mathcal N=0$.
\end{proof}

\begin{lem}\label{extact} Let $\chi:H\rightarrow \Fbar^{\times}$ be a character. We consider $\chi$ and $\chi^{s}$ as 
one dimensional representations of $I$, 
via $I/I_1\cong H$. Let 
\begin{displaymath}
\chi\hookrightarrow \Inj \chi,\quad \chi^{s}\hookrightarrow \Inj \chi^{s}
\end{displaymath}
be  injective envelopes of $\chi$ and $\chi^{s}$ in $\Rep_I$, respectively. Let $V_1$ be the underlying
vector space of $\Inj \chi$ and let $V_2$ be the underlying vector space of $\Inj \chi^{s}$. Further, let $v_1$  and $v_2$ be vectors in $V_1$ and $V_2$ 
respectively, such  that
$$\langle v_1 \rangle_{\Fbar}=(\Inj \chi)^{I_1},\quad \langle v_2 \rangle_{\Fbar}=(\Inj \chi^{s})^{I_1}.$$
Then there exists an action of $\kon$ on $V_1\oplus V_2$, extending the action of $I$, so that our fixed uniformiser $\pif$ acts trivially and 
$$\Pi^{-1}v_1 = v_2,\quad \Pi^{-1}v_2= v_1.$$
\end{lem}
\begin{proof}  Let $t\in T$ be any diagonal matrix, then $ s t s= \Pi t\Pi^{-1}$. Hence 
$$\chi^{s} \cong \chi^{\Pi}$$
as $I$-representations, where $\chi^{\Pi}$ denotes the action of $I$, on the underlying vector space of $\chi$, twisted by $\Pi$. 
So we get an exact sequence
\begin{displaymath}
\xymatrix{0\ar[r]&\chi^{s}\ar[r]& (\Inj \chi)^{\Pi}}
\end{displaymath} 
of $I$-representations. Twisting by $\Pi$ is an exact functor in $\Rep_I$ and 
$$\Hom_I(\xi, (\Inj \chi)^{\Pi})\cong \Hom_I(\xi^{\Pi}, \Inj \chi).$$
Since $\Inj \chi$ is an injective object in $\Rep_I$, this implies that
$(\Inj \chi)^{\Pi}$ is an injective object in $\Rep_I$. Since $\Inj \chi$ is an essential extension of $\chi$, $(\Inj \chi)^{\Pi}$ is an essential
extension of $\chi^{s}$. Since injective envelopes are unique up to isomorphism, there exists an isomorphism $\phi$ of $I$-representations 
$$ \phi: (\Inj \chi)^{\Pi}\cong \Inj \chi^{s}.$$
The proof of Lemma \ref{restI} shows that the space $(\Inj \chi)^{I_1}$ is one dimensional. Hence, after replacing $\phi$ by a scalar multiple 
we may assume that $\phi(v_1) =v_2$. We may extend the action of $I$ on $V_1$ and $V_2$ to the action of $F^{\times}I$ by making $\pif$ act trivially. 
 We denote the corresponding representations by $\widetilde{\Inj} \chi$ and $\widetilde{\Inj}\chi^{s}$. For trivial reasons  
$$\phi: (\widetilde{\Inj}\chi)^{\Pi}\cong \widetilde{\Inj}\chi^{s}.$$ 
We consider the induced representation
$\Indu{F^{\times}I}{\kon}{\widetilde{\Inj}\chi}$. Let $\ev_1$ and $\ev_{\Pi}$ be the evaluation maps at $1$ and $\Pi$ respectively, then
we get an $F^{\times}I$-equivariant isomorphism: 
$$\Indu{F^{\times}I}{\kon}{\widetilde{\Inj}\chi}\cong V_1\oplus V_2, \quad f\mapsto \ev_1(f)+ \phi(\ev_{\Pi}(f)).$$
The action of $\kon$ on the left hand side gives us the action of $\kon$ on $V_1\oplus V_2$. Let $v\in V_1$ and $w\in V_2$, then the action of $\Pi^{-1}$
is given by
$$\Pi^{-1}(v+w)=\phi^{-1}(w)+ \phi(v)$$ 
and hence $\Pi^{-1}v_1=v_2$ and $\Pi^{-1}v_2=v_1$.
\end{proof} 
We will construct a diagram $Y_{\gamma}$. This will involve making some choices. Suppose that $q=p^n$, let $\chi:H\rightarrow \Fbar^{\times}$
be a character and let $\gamma=\{\chi,\chi^{s}\}$.  We consider an irreducible representation $\rho_{\chi,J}$ of $K$. Lemma \ref{dictionary} 
gives us a pair $(\rr, a)$, where $\rr$ is 
the usual $n$-tuple and $a$ is an integer modulo $q-1$. Let $\rho_{\chi,J}\hookrightarrow \Inj \rho_{\chi,J}$ be an injective envelope of 
$\rho_{\chi,J}$ in $\Rep_{K}$. Let $\WW_{\rr}$ be the underlying vector space of $\Inj \rho_{\chi,J}$. We may assume that 
$\WW_{\rr}$ depends only on the $n$-tuple $\rr$. Since, if $\chi'=\chi\otimes (\det)^{c}$, then
$\rho_{\chi',J}\cong \rho_{\chi,J}\otimes (\det)^{c}$ and a simple argument shows that $\rho_{\chi',J}\hookrightarrow (\Inj \rho_{\chi,J})\otimes(\det)^c$
is an injective envelope of $\rho_{\chi',J}$ in $\Rep_K$. Let 
$$Y_{\gamma,0}=(\widetilde{\Inj}\rho_{\chi,J}\oplus \widetilde{\Inj}\rho_{\chi^{s},\overline{J}}, \WW_{\rr}\oplus \WW_{\mathbf{p-1}-\rr})$$
where tilde denotes the extension of the action of $K$ to the action of $F^{\times}K$, so that $\pif$ acts trivially.  We are going to construct 
an action of $\kon$ on $Y_{\gamma,0}|_{F^{\times}I}$, which extends the action of $F^{\times}I$, and this will give us $Y_{\gamma}$.
 However, this can be done in a lot of ways, and not all of them suit our purposes. Lemma \ref{injinv} and Remark \ref{note} imply that
$$(Y_{\gamma,0})^{K_1}\cong \inj \rho_{\chi,J}\oplus \inj \rho_{\chi^{s},\overline{J}}$$
as $K$-representations, where on the right hand  side we adopt the notation of Propositions \ref{injiwahori} and \ref{injregular}. In particular, 
$$(Y_{\gamma, 0})^{I_1}\cong ( \inj \rho_{\chi,J}\oplus \inj \rho_{\chi^{s},\overline{J}})^{I_1}$$
as $\HH_{K}$-modules. In Lemma \ref{mrfixit} we have worked out a basis consisting of eigenvectors for the action of $I$ of (a model of) 
$( \inj \rho_{\chi,J}\oplus \inj \rho_{\chi^{s},\overline{J}})^{I_1}$. The above isomorphism gives us a basis $\BB$ of $(Y_{\gamma,0})^{I_1}$.
Lemma \ref{restI} gives an $F^{\times}I$-equivariant decomposition: 
$$\zeta: \WW_{\rr}\oplus \WW_{\mathbf{p-1}-\rr}\cong \bigoplus_{b\in \BB} \WW(b)$$
such that $\zeta(b)\in \WW(b)$, for every $b\in \BB$, and the representation, given by the action of  $I$ on $\WW(b)$, is an injective object in $\Rep_I$,
which is also an essential extension of $\langle \zeta(b) \rangle_{\Fbar}$. To simplify things we view $\zeta$ as identification and omit it
from our notation. 
 
If $\chi=\chi^{s}$ then we pair up the basis vectors as in Proposition \ref{injiwahori}:
$$\BB= \{b_{\mathbf{0}}, b_{\mathbf{0}}+b_{\mathbf{1}}\}\bigcup_{\{\varepsilon,\mathbf{1}-\varepsilon\}\subseteq \Sigma'_{\mathbf{0}}}
\{b_{\varepsilon}, b_{\mathbf{1}-\varepsilon}\}.$$
If $\chi\neq \chi^{s}$ then we pair up the basis vectors as in Proposition \ref{injregular}:
$$\BB= \{b_{\mathbf{0}},\bar{b}_{\mathbf{0}}\}\cup\{b_{\mathbf{1}},\bar{b}_{\mathbf{1}}\}
\bigcup_{\{\varepsilon,\mathbf{1}-\varepsilon\}\subseteq \Sigma'_{\rr}}\{b_{\varepsilon}, b_{\mathbf{1}-\varepsilon}\}
\bigcup_{\{\varepsilon,\mathbf{1}-\varepsilon\}\subseteq \Sigma'_{\mathbf{p-1}-\rr}}\{\bar{b}_{\varepsilon}, \bar{b}_{\mathbf{1}-\varepsilon}\}.$$
Let $\{b,b'\}$ be any such pair and suppose that $I$ acts on $b$ via a character $\psi$, then $I$ will act on $b'$, via a character $\psi^{s}$. 
 We denote 
 $$\WW(b,b')=\WW(b)\oplus\WW(b').$$
 Lemma \ref{extact} implies that there exists an action of $\kon$ on $\WW(b,b')$, extending the action of $F^{\times}I$, such that  
$$\Pi^{-1}b=b', \quad \Pi^{-1}b'= b.$$
This amounts to fixing an isomorphism of vector spaces $\phi: \WW(b)\cong\WW(b')$, such that $\phi(b)=b'$ and which induces an isomorphism of $I$ 
representations $\phi: (\Inj \psi)^{\Pi}\cong \Inj \psi^{s}$.

If $\chi=\chi^{s}$ then $Y_{\gamma,0}$ decomposes into $F^{\times}I$-invariant subspaces:
$$ \WW(b_{\mathbf{0}}, b_{\mathbf{0}}+b_{\mathbf{1}})
\bigoplus_{\{\varepsilon,\mathbf{1}-\varepsilon\}\subseteq \Sigma'_{\mathbf{0}}} \WW(b_{\varepsilon}, b_{\mathbf{1}-\varepsilon}).$$
If $\chi\neq \chi^{s}$ then $Y_{\gamma,0}$ decomposes into $F^{\times}I$-invariant subspaces: 
$$ \WW(b_{\mathbf{0}}, \bar{b}_{\mathbf{0}})\oplus \WW(b_{\mathbf{1}}, \bar{b}_{\mathbf{1}})
\bigoplus_{\{\varepsilon,\mathbf{1}-\varepsilon\}\subseteq \Sigma'_{\rr}} \WW(b_{\varepsilon}, b_{\mathbf{1}-\varepsilon})
\bigoplus_{\{\varepsilon,\mathbf{1}-\varepsilon\}\subseteq \Sigma'_{\mathbf{p-1}-\rr}} \WW(\bar{b}_{\varepsilon}, \bar{b}_{\mathbf{1}-\varepsilon}).$$
Let $Y_{\gamma,1}$ be a representation of $\kon$, whose underlying vector space is $\WW_{\rr}\oplus \WW_{\mathbf{p-1}-\rr}$, and the action of 
$\kon$ extends the action of $F^{\times}I$ on each direct summand, as it was done for $\WW(b,b')$. 

\begin{defi}\label{ygamma} Let $Y_{\gamma}$ be an object in $\diag$, given by 
$$Y_{\gamma}=(Y_{\gamma,0},Y_{\gamma,1}, \id)$$
and let $\II_{\gamma}$ be the corresponding coefficient system 
$$\II_{\gamma}=\CC(Y_{\gamma}).$$
\end{defi}

\begin{remar} The definition of $Y_{\gamma}$ depends on all the choices we have made.
\end{remar}

\begin{prop}\label{Injmod} Let $\chi:H\rightarrow \Fbar^{\times}$ be a character and let $\gamma=\{\chi,\chi^{s}\}$. Suppose that $\chi=\chi^s$, then
$$H_0(X,\II_{\gamma})^{I_1}\cong M_{\gamma} \bigoplus_{\{\varepsilon,\mathbf{1}-\varepsilon\}\subseteq \Sigma'_{\mathbf{0}}} M_{\gamma_{\varepsilon}}$$
as $\HH$-modules, where 
$\gamma_{\mathbf{\varepsilon}}=\gamma_{\mathbf{1-\varepsilon}}=\{\chi\alpha^{\mathbf{\varepsilon}\centerdot(\mathbf{p-1})},
 \chi(\alpha^{\mathbf{\varepsilon}\centerdot(\mathbf{p-1})})^{s}\}$. Suppose that $\chi\neq \chi^{s}$, then
$$H_0(X,\II_{\gamma})^{I_1}\cong L_{\gamma}\bigoplus_{\{\varepsilon,\mathbf{1}-\varepsilon\}\subseteq \Sigma'_{\rr}}M_{\gamma_{\mathbf{\varepsilon}}}
\bigoplus_{\{\varepsilon,\mathbf{1}-\varepsilon\}\subseteq \Sigma'_{\mathbf{p-1}-\rr}} M_{\bar{\gamma}_{\varepsilon}}$$
as $\HH$-modules, where $\gamma_{\mathbf{\varepsilon}}=\gamma_{\mathbf{1-\varepsilon}}=
\{\chi\alpha^{\mathbf{\varepsilon}\centerdot(\mathbf{p-1}-\rr)},
 (\chi\alpha^{\mathbf{\varepsilon}\centerdot(\mathbf{p-1}-\rr)})^{s}\}$ and 
 $\bar{\gamma}_{\mathbf{\varepsilon}}=\bar{\gamma}_{\mathbf{1-\varepsilon}}=\{\chi^{s}\alpha^{\mathbf{\varepsilon}\centerdot \rr},
 (\chi^{s}\alpha^{\mathbf{\varepsilon}\centerdot \rr })^{s}\}$.
\end{prop}
\begin{proof}In Propositions \ref{injiwahori} and \ref{injregular} we have showed that we may extend the action of $\HH_K$ on
$(\inj \rho_{\chi,J}\oplus \inj \rho_{\chi^{s},\overline{J}})^{I_1}$
to the action of $\HH$, so that the resulting modules are isomorphic to the ones considered above. We will show that $H_0(X,\II_{\gamma})^{I_1}$
realizes this extension. 
By Proposition \ref{isores} (or alternatively Corollary \ref{thepoint}) we have
$$H_0(X,\II_{\gamma})|_{\kzr}\cong Y_{\gamma,0},\quad H_0(X,\II_{\gamma})|_{\kon}\cong Y_{\gamma,1}$$
as $\kzr$ and $\kon$-representations, respectively. Moreover, the diagram
\begin {displaymath}
\xymatrix {Y_{\gamma,_0} \ar[r]^{\cong} & H_0(X,\II_{\gamma}) \\
          Y_{\gamma, 1} \ar[u]^{\id} \ar[r]^{\cong} & H_0(X,\II_{\gamma}) \ar[u]_{\id}}
\end {displaymath}
of $F^{\times}I$-representations commutes. So 
$$(Y_{\gamma,0})^{I_1}\cong H_0(X,\II_{\gamma})^{I_1}$$
as $\HH_K$-modules. Lemma \ref{injinv} implies that 
$$H_0(X,\II_{\gamma})^{I_1}\cong (\inj \rho_{\chi,J}\oplus \inj \rho_{\chi^{s},\overline{J}})^{I_1}$$
as $\HH_K$-modules, and we know the right hand side from Propositions \ref{injiwahori} and \ref{injregular}.
It remains to determine the action of $T_{\Pi}$. Corollary \ref{T}  implies that for every $v\in H_0(X,\II_{\gamma})^{I_1}$ we have 
$$v T_{\Pi}= \Pi^{-1}v.$$
Hence the action of $T_{\Pi}$ is determined by the isomorphism
$$Y_{\gamma, 1}\cong H_0(X,\II_{\gamma})|_{\kon}.$$
Since $\BB$ is a basis of ($Y_{\gamma,0})^{I_1}$, it is enough to know how $\Pi^{-1}$ acts on the basis vectors. Let $\WW(b,b')$ be one of the 
$\kon$-invariant subspaces of $Y_{\gamma,1}$, as before.  We have extended the action of $F^{\times}I$ on $Y_{\gamma,0}|_{F^{\times}I}$ to $\kon$ so that 
$$\Pi^{-1}b=b',\quad \Pi^{-1}b'= b.$$
Hence, if we consider $\BB$ also as a basis of $H_0(X,\II_{\gamma})^{I_1}$ we have 
$$b T_{\Pi}= b',\quad b'T_{\Pi}=b.$$  
 Now the statement of the Proposition is just a realization of Propositions \ref{injiwahori} and \ref{injregular}.
\end{proof}
\subsection{Construction}
Now we will construct an embedding $D_{\gamma}\hookrightarrow Y_{\gamma}$.
Suppose that $\chi=\chi^{s}$, then we consider vectors $b_{\mathbf{0}}$ and $b_{\mathbf{0}}+b_{\mathbf{1}}$ in $(Y_{\gamma,0})^{I_1}$. 
Lemmas \ref{genezero} and \ref{gensteinberg} imply that
$$\langle K b_{\mathbf{0}}\rangle_{\Fbar}\cong \tilde{\rho}_{\chi,S},\quad \langle K(b_{\mathbf{0}}+b_{\mathbf{1}})\rangle_{\Fbar}
\cong \tilde{\rho}_{\chi,\emptyset}$$ 
as $F^{\times}K$-representations. We have constructed the action of $\kon$ on $Y_{\gamma,1}$ so that
$$\Pi^{-1}b_{\mathbf{0}}=b_{\mathbf{0}}+b_{\mathbf{1}}, \quad \Pi^{-1} (b_{\mathbf{0}}+b_{\mathbf{1}})=b_{\mathbf{0}}.$$
Suppose that $\chi\neq \chi^{s}$, then we consider vectors $b_{\mathbf{0}}$ and $\bar{b}_{\mathbf{0}}$ in $(Y_{\gamma,0})^{I_1}$. 
Lemmas \ref{genezero}  implies that
$$\langle K b_{\mathbf{0}}\rangle_{\Fbar}\cong \tilde{\rho}_{\chi,\emptyset},\quad \langle K \bar{b}_{\mathbf{0}}\rangle_{\Fbar}
\cong \tilde{\rho}_{\chi^{s},\emptyset}$$ 
as $F^{\times}K$-representations. We have constructed the action of $\kon$ on $Y_{\gamma,1}$ so that
$$\Pi^{-1}b_{\mathbf{0}}=\bar{b}_{\mathbf{0}}, \quad \Pi^{-1} \bar{b}_{\mathbf{0}}=b_{\mathbf{0}}.$$
Hence, in both cases we get an embedding $D_{\gamma}\hookrightarrow Y_{\gamma}$ in the category $\diag$. This gives us an embedding of $G$ equivariant
coefficient systems $\VV_{\gamma}\hookrightarrow \II_{\gamma}$.
\begin{defi}Let $\pi_{\gamma}$ be a representation of $G$, given by
$$\pi_{\gamma}=\Image(H_0(X,\VV_{\gamma})\rightarrow H_0(X,\II_{\gamma})).$$
\end{defi}
\begin{thm}\label{irrconst}For each $\gamma=\{\chi,\chi^{s}\}$, the representation $\pi_{\gamma}$ is irreducible and supersingular. Moreover, 
$\pi_{\gamma}^{I_1}$ 
contains an $\HH$-submodule isomorphic to $M_{\gamma}$. Further, if 
$$\pi_{\gamma}\cong \pi_{\gamma'}$$
then $\gamma=\gamma'$.
\end{thm}
\begin{proof}Lemma \ref{inj} implies that $\pi_{\gamma}$ is non-zero. So by Corollary \ref{superquot} it is enough to prove that $\pi_{\gamma}$ is 
irreducible. 
To ease the notation we identify the underlying vector spaces of $Y_{\gamma,0}$ and $H_0(X,\II_{\gamma})$. If $\chi=\chi^{s}$ then Lemma \ref{omegagener}
implies that
$$\pi_{\gamma}=\langle G b_{\mathbf{0}}\rangle_{\Fbar}=\langle G (b_{\mathbf{0}}+b_{\mathbf{1}})\rangle_{\Fbar}.$$
If $\chi\neq \chi^{s}$ then Lemma \ref{omegagener} implies that
$$\pi_{\gamma}=\langle G b_{\mathbf{0}}\rangle_{\Fbar}=\langle G \bar{b}_{\mathbf{0}}\rangle_{\Fbar}.$$
This can be rephrased in a different way. By Proposition \ref{isores} we have
$$H_0(X,\II_{\gamma})|_{K}\cong \Inj \rho_{\chi,J}\oplus \Inj \rho_{\chi^{s},\overline{J}}$$
as $K$-representations. Lemma \ref{socle} implies that
$$\rho_{\chi,J}\oplus \rho_{\chi^{s},\overline{J}}\cong \soc (H_0(X,\II_{\gamma})|_{K}).$$
Hence, if $\chi=\chi^{s}$ then 
$$(\soc(H_0(X,\II_{\gamma})|_{K}))^{I_1}=\langle b_{\mathbf{0}}, b_{\mathbf{0}}+b_{\mathbf{1}}\rangle_{\Fbar}$$
and if $\chi\neq\chi^{s}$ then
$$ (\soc(H_0(X,\II_{\gamma})|_{K}))^{I_1}=\langle b_{\mathbf{0}},\bar{b}_{\mathbf{0}}\rangle_{\Fbar}$$
and hence
$$\pi_{\gamma}=\langle G (\soc(H_0(X,\II_{\gamma})|_{K}))^{I_1}\rangle_{\Fbar}.$$
Suppose that $\pi'$ is non-zero $G$-invariant subspace of $\pi_{\gamma}$ then by Lemma \ref{propinvariants} $(\pi')^{K_1}\neq 0$, and hence
$$\soc(\pi'|_{K})\neq 0.$$
We have trivially $\soc(\pi'|_{K})\subseteq \soc(H_0(X,\II_{\gamma})|_{K})$. Hence
$$(\soc(H_0(X,\II_{\gamma})|_{K}))^{I_1}\cap (\pi')^{I_1}\neq 0.$$
The space $(\soc(H_0(X,\II_{\gamma})|_{K}))^{I_1}$ is $\HH$-invariant, and in fact isomorphic to the irreducible module $M_{\gamma}$. Hence,
$$(\soc(H_0(X,\II_{\gamma})|_{K}))^{I_1}\subseteq (\pi')^{I_1}$$
and this implies that $\pi'=\pi_{\gamma}$. Hence $\pi_{\gamma}$ is irreducible. 

Suppose that $\pi_{\gamma}\cong \pi_{\gamma'}$, then this induces an isomorphism of vector spaces
$$ \phi:(\soc(\pi_{\gamma}|_{K}))^{I_1}\cong(\soc(\pi_{\gamma'}|_{K}))^{I_1}.$$
The argument above implies that both spaces are $\HH$-invariant and  Corollary \ref{T} implies that $\phi$ is an  
isomorphism of $\HH$-modules. Hence,
$$M_{\gamma}\cong(\soc(\pi_{\gamma}|_{K}))^{I_1}\cong(\soc(\pi_{\gamma'}|_{K}))^{I_1}\cong M_{\gamma'}.$$
Lemma \ref{gammagamma} implies that $\gamma=\gamma'$.
\end{proof}
\begin{cor}\label{soclegamma} The representation $H_0(X, \II_{\gamma})$ is an essential extension of $\pi_{\gamma}$ in $\Rep_G$. In particular, 
$$\pi_{\gamma}\cong \soc(H_0(X,\II_{\gamma})),$$
where $\soc(H_0(X,\II_{\gamma}))$ is the subspace of $H_0(X,\II_{\gamma})$ generated by all the irreducible subrepresentations.
\end{cor}
\begin{proof} Let $\pi$ be a non-zero $G$-invariant subspace of $H_0(X, \II_{\gamma})$. The proof of Theorem \ref{irrconst} shows that
$(\soc(H_0(X,\II_{\gamma})|_K))^{I_1}$ is a subspace of $\pi^{I_1}$. This implies that $\pi_{\gamma}$ is a subspace of $\pi$. The last part 
is immediate. 
\end{proof}

\subsubsection{Twists by unramified quasi-characters}
Let $\lambda\in \Fbar^{\times}$, we define an unramified quasi-character $\mu_{\lambda}:F^{\times}\rightarrow \Fbar^{\times}$, by
$$\mu_{\lambda}(x)= \lambda^{\val_{F}(x)}.$$
\begin{lem}\label{mulambda} Suppose that $\pi_{\gamma}\otimes \mu_{\lambda}\circ \det \cong \pi_{\gamma'}$, then $\gamma=\gamma'$ and $\lambda=\pm 1$.
\end{lem}
\begin{proof} Our fixed uniformiser $\pif$ acts on $\pi_{\gamma}\otimes \mu_{\lambda}\circ \det $, by a scalar $\lambda^2$, and it acts trivially 
on $\pi_{\gamma'}$. Hence, $\lambda=\pm 1$.
By Lemma \ref{twistbyun} $M_{\gamma}\otimes\mu_{-1}\circ\det\cong M_{\gamma}$, and hence by the argument of \ref{irrconst} $M_{\gamma'}\cong M_{\gamma}$, 
which implies that $\gamma=\gamma'$.
\end{proof}

\begin{prop}\label{minusone} Suppose that $q=p$, then $\pi_{\gamma}\otimes(\mu_{-1}\circ\det)\cong \pi_{\gamma}$.
\end{prop}
\begin{proof} By Corollary \ref{soclegamma} it is enough to show that $Y_{\gamma}\otimes (\mu_{-1}\circ\det)\cong Y_{\gamma}$
in $\diag$. We claim that we always have 
$$Y_{\gamma,1}\cong Y_{\gamma,1}\otimes(\mu_{-1}\circ\det)$$ 
as  $\kon$-representations. Since $F^{\times}I$ is contained in the kernel of $\mu_{-1}\circ\det$, it is enough to examine the action of $\Pi$. 
We recall that the action of 
$\kon$ was defined, by fixing a certain  isomorphism $\phi:\WW(b)\cong \WW(b')$, and then letting $\Pi^{-1}$ act on $\WW(b,b')=\WW(b)\oplus\WW(b')$ by
$$\Pi^{-1}(v+w)= \phi^{-1}(w)+\phi(v).$$
Let $\iota_1$ be an $F^{\times}I$-equivariant isomorphism 
$$\iota_1: \WW(b)\oplus\WW(b')\cong \WW(b)\oplus\WW(b'),\quad v+w\mapsto v-w,$$
then, since $\mu_{-1}(\det(\Pi^{-1}))=-1$, we have
$$\Pi^{-1}\otimes\mu_{-1}(\det(\Pi^{-1}))(\iota_1(v+w))=\phi^{-1}(w)-\phi(v)=\iota_1(\Pi^{-1}(v+w)).$$ 
Hence $\WW(b,b')\cong \WW(b,b')\otimes (\mu_{-1}\circ\det)$ as $\kon$-representations and hence $Y_{\gamma,1}\cong Y_{\gamma,1}\otimes 
(\mu_{-1}\circ\det)$ as $\kon$-representations. Since $F^{\times}K$ is contained in the kernel of $\mu_{-1}\circ\det$ we also have 
$Y_{\gamma,0}\cong Y_{\gamma,0}\otimes(\mu_{-1}\circ\det)$. However, 
to define an isomorphism in $\diag$ we need to find $\iota_0:Y_{\gamma,0}\cong Y_{\gamma, 0}$, which is compatible with $\iota_1$ via the restriction 
maps. If $p=q$ this is easy, since 
if $\chi=\chi^{s}$, then
$$\WW_{\rr}\oplus\WW_{\mathbf{p-1}-\rr}=\WW(b_{\mathbf{0}})\oplus \WW(b_{\mathbf{0}}+b_{\mathbf{1}})$$
and if $\chi\neq \chi^{s}$ then 
$$\WW_{\rr}\oplus\WW_{\mathbf{p-1}-\rr}=(\WW(b_{\mathbf{0}})\oplus \WW(b_{\mathbf{1}}))\oplus 
(\WW(\bar{b}_{\mathbf{0}})\oplus\WW(\bar{b}_{\mathbf{1}}))$$
and the subspaces that $\Pi$ `swaps' come from different injective envelopes. Note, that this is not the case if $q\neq p$. Hence, if we define
$$\iota_0:\WW_{\rr}\oplus \WW_{\mathbf{p-1}-\rr}\cong\WW_{\rr}\oplus \WW_{\mathbf{p-1}-\rr},\quad v+w\mapsto v-w$$
then $\iota=(\iota_0,\iota_1)$ is an isomorphism $\iota:Y_{\gamma}\cong Y_{\gamma}\otimes(\mu_{-1}\circ\det)$.
\end{proof}

\begin{lem}\label{adm}The representations $H_0(X,\II_{\gamma})$  and $\pi_{\gamma}$ are  admissible. 
\end{lem}
\begin{proof}Proposition \ref{Injmod}, Lemma \ref{ohdear}.
\end{proof}

Our main result can be summarised as follows.
\begin{thm}\label{summary}Let $\pif$ be a fixed uniformiser, then there exists at least $q(q-1)/2$ pairwise non-isomorphic, irreducible, supersingular,
admissible representations of $G$, which admit a central character, such that $\pif$ acts trivially.
\end{thm} 
\begin{proof} There are precisely $q(q-1)/2$ orbits $\gamma=\{\chi,\chi^{s}\}$. Then the statement follows from Theorem \ref{irrconst} and 
Corollary \ref{adm}. Each $\pi_{\gamma}$ admits a central character, since $H_0(X,\VV_{\gamma})$ admits a central character. 
If $\lambda\in \oF^{\times}$, then it acts on $H_0(X,\VV_{\gamma})$ by a scalar
$$ \chi(\begin{pmatrix} \lambda & 0 \\0 & \lambda \end{pmatrix})=\chi^{s}(\begin{pmatrix} \lambda & 0 \\0 & \lambda \end{pmatrix})$$
and $\pif$ acts trivially by construction.
\end{proof}
If $F=\QQ_p$ then we may apply the results of Breuil \cite{breuil}. 
\begin{cor}Suppose that $F=\QQ_p$, then  $\pi_{\gamma}$ is independent up to isomorphism of the choices made in the construction of $Y_{\gamma}$.
Moreover, if $\pi$ is an irreducible supersingular representation of $G$, admitting a central character,  then there exists $\lambda\in \Fbar^{\times}$, 
unique up to a sign, and a unique $\gamma$, such that
$$\pi\cong \pi_{\gamma}\otimes(\mu_{\lambda}\circ \det).$$
\end{cor}
\begin{proof}In \cite{breuil} Breuil has determined all the supersingular representations, in the case of $F=\QQ_p$. As a consequence, by \cite{vig}
Theorem E.7.2, the functor of $I_1$-invariants, $\Rep_{G}\rightarrow \Mod-\HH$, $\pi\mapsto \pi^{I_1}$ induces a bijection between the isomorphism classes
of irreducible supersingular representations with a central character and isomorphism classes of supersingular right modules of $\HH$. In particular, 
there are precisely $p(p-1)/2$ isomorphism classes of supersingular representations with a central character, such that $\pif$ acts trivially. By Theorem
\ref{summary} our construction yields at least $p(p-1)/2$ such representations. Hence $\pi_{\gamma}$ does not depend up to isomorphism on the choices 
made for $Y_{\gamma}$. 

Let $\pi$ be any supersingular representation of $G$ with a central character. We may always twist $\pi$ by an unramified quasi-character, so that $\pif$ 
acts trivially. Hence by above
$$\pi\cong \pi_{\gamma}\otimes(\mu_{\lambda}\circ\det)$$
and by Lemma \ref{mulambda} and Proposition \ref{minusone}, $\gamma$ is determined uniquely and $\lambda$ up to $\pm 1$. 
\end{proof}

\end{document}